\documentclass[poms,final,nonblindrev]{poms1_V1}
\OneAndAHalfSpacedXI

\usepackage{amsmath,amsfonts,amssymb}
\usepackage{makecell,comment}
\usepackage{algorithm}
\usepackage{algorithmicx}
\usepackage{algpseudocode}

\floatname{algorithm}{Algorithm}

\allowdisplaybreaks[4]
\usepackage{esdiff}
\usepackage{mathtools}
\usepackage{ulem}
\usepackage{multirow,threeparttable,longtable,rotating,booktabs}
\usepackage{enumerate}

\usepackage{mathabx}
\usepackage{color,bm,graphicx,booktabs}
\usepackage{xspace}
\usepackage[utf8]{inputenc}
\usepackage[round]{natbib}
\usepackage{xcolor}
\definecolor{hexred}{HTML}{C91C36}

\def\Halmos{\mbox{\quad$\square$}}
\usepackage{color-edits}
\addauthor{qt}{hexred}
\addauthor{hy}{brown}
\addauthor{sy}{magenta}
\addauthor{lc}{cyan}
\addauthor{te}{red}
\usepackage{xurl}
\usepackage[hypertexnames=false]{hyperref}
\hypersetup{
	colorlinks = true,
	linkcolor=black!70!green,          
	citecolor=black!70!green,          
	filecolor=magenta,                 
	urlcolor=black!70!green,           
	linktocpage=true
}
\usepackage{cleveref}
\usepackage{anyfontsize}
\makeatletter
\@ifundefined{abovecaptionskip}{\newlength\abovecaptionskip}{}
\makeatother
\usepackage{graphicx}

\TheoremsNumberedThrough
\ECRepeatTheorems
\EquationsNumberedThrough

\theoremstyle{TH}

\RUNAUTHOR{Zhang et al.}
\RUNTITLE{Robust Inventory Planning and Advertising under Product Return}
\TITLE{Data-Driven Robust Inventory Planning and Advertising under Product Return}

\ARTICLEAUTHORS{%
\AUTHOR{Hongyan Zhang}
\AFF{University of Science and Technology of China, \EMAIL{hyzh@ustc.edu.cn}}
\AUTHOR{Yugang Yu}
\AFF{University of Science and Technology of China, \EMAIL{ygyu@ustc.edu.cn}}
\AUTHOR{Qinshen Tang}
\AFF{Nanyang Technological University, \EMAIL{qinshen.tang@ntu.edu.sg}}
\AUTHOR{Ye Shi\footnote{Corresponding author. The authors are listed in reverse alphabetical order with implied equal authorship.}}
\AFF{University of Science and Technology of China, \EMAIL{hecules@ustc.edu.cn}}
\AUTHOR{Li Chen}
\AFF{The University of Sydney, \EMAIL{li.chen@sydney.edu.au}}
}

\newcommand{\mf}[1]{\mathbf{#1}}
\newcommand{\mfu}[1]{\uline{\mkern-1mu{\mathbf{#1}}\mkern-1mu}}

\newcommand{\mfh}[1]{\hat{\mathbf{#1}}}

\newcommand{\mfhu}[1]{\hat{\uline{\mathbf{#1}}}}

\newcommand{\bmt}[1]{\tilde{\bm{#1}}}
\newcommand{\bmh}[1]{\hat{\bm{#1}}}

\newcommand{\bmd}[1]{\bm{#1}^{\dag}}

\newcommand{\bmu}[1]{\uline{\mkern-1mu{\bm{#1}}\mkern-1mu}}

\newcommand{\tu}[1]{\uline{\mkern-1mu{\tilde{#1}}\mkern-1mu}}

\newcommand{\bmtu}[1]{\uline{\tilde{\bm{#1}}}}
\newcommand{\bmhu}[1]{\uline{\hat{\bm{#1}}}}
\newcommand{\bmdu}[1]{\uline{\bm{#1}}^{\dag}}

\newcommand{\ep}[1]{\mathbb{E}_{\mathbb{P}}\left[#1\right]}
\newcommand{\eps}[1]{\mathbb{E}_{\mathbb{P}_s}\left[#1\right]}

\newcommand{\ephat}[1]{\mathbb{E}_{\hat{\mathbb{P}}}\left[#1\right]}
\newcommand{\epstar}[1]{\mathbb{E}_{\mathbb{P}^\star}\left[#1\right]}

\newcommand{\bbp}{\mathbb{P}}
\newcommand{\bbr}{\mathbb{R}}
\newcommand{\bbq}{\mathbb{Q}}

\newcommand{\cx}{\mathcal{X}}

\newcommand{\cz}{\mathcal{Z}}
\newcommand{\cy}{\mathcal{Y}}
\newcommand{\cq}{\mathcal{Q}}
\newcommand{\cf}{\mathcal{F}}
\newcommand{\ci}{\mathcal{I}}
\newcommand{\cl}{\mathcal{L}}

\newcommand{\prox}{p}

\makeatletter
\newcommand*\bigcdot{\mathpalette\bigcdot@{.75}}
\newcommand*\bigcdot@[2]{\mathbin{\vcenter{\hbox{\scalebox{#2}{$\m@th#1\bullet$}}}}}
\makeatother
\newcommand{\trace}[1]{\mathrm{tr}(#1)}

\newcommand{\comm}[1]{{\color{red} #1}}
\definecolor{revisionteal}{RGB}{0,110,125}
\newcommand{\revmonthly}[1]{\textcolor{revisionteal}{#1}}
\definecolor{revisionsigned}{RGB}{180,90,0}
\newcommand{\revgain}[1]{\textcolor{revisionsigned}{#1}}
\providecommand{\qtcomment}[1]{\textcolor{hexred}{[QT: #1]}}

\newcommand{\cred}[1]{\textcolor{red}{#1}}
\newcommand{\cgreen}[1]{\textcolor{green!45!black}{#1}}

\newcommand{\code}[1]{\texttt{\bfseries{#1}}}

\bibpunct[, ]{(}{)}{,}{a}{}{,}
\def\bibfont{\small}

\crefname{assm}{Assumption}{Assumptions}
\Crefname{assm}{Assumption}{Assumptions}
\crefname{prop}{Proposition}{Propositions}
\Crefname{prop}{Proposition}{Propositions}
\crefname{thm}{Theorem}{Theorems}
\Crefname{thm}{Theorem}{Theorems}
\crefname{coro}{Corollary}{Corollaries}
\Crefname{coro}{Corollary}{Corollaries}
\crefname{rem}{Remark}{Remarks}
\Crefname{rem}{Remark}{Remarks}
\begin{document}

\ABSTRACT{
Online advertising is a central marketing lever for e-commerce retailers, yet it complicates operational planning by simultaneously stimulating product demand and increasing product returns. Motivated by the practice of a fashion retailer on Amazon, we study a retailer's joint inventory and advertising problem across multiple products over a single season, in which advertising has both direct and cross-product spillover effects on demand and returns.  Using linear net-demand and return models that capture both effects, we formulate a contextual stochastic optimization model that jointly determines inventory levels and advertising decisions before the selling season. Because ad clicks are realized only after the advertising decision, the predictive covariates are partially observable and decision-dependent, and their joint distribution with the demand and return errors is unknown. We therefore develop an estimate-then-robustify (ETR) approach grounded in a data-driven robust satisficing framework. We establish a finite-sample guarantee on the out-of-sample profit of the ETR solution that holds for any parameter estimator, show that the solution is asymptotically optimal when the estimator is consistent and the profit target approaches the empirical optimal profit, and derive a mixed-integer conic approximation via affine recourse adaptation that is exact for a single product. On the retailer's data, ETR improves out-of-sample profit over conventional estimate-then-optimize by 6.39\% and earns the higher profit in 74.5\% of instances, with the largest gains when training data are scarce; the gain falls to 1.02\% when the decision model omits cross-product spillovers. 
}
\maketitle

\section{Introduction}\label{sec:introduction}

Retailers increasingly rely on online advertising, such as sponsored paid search ads, to boost sales and strengthen their market share in e-commerce markets \citep{ du2017bidding, kim2023coordinating}. For example, 65\% of small and mid-sized businesses run their own paid search campaigns on the Google Platform, a key contributor to Google's staggering \$264.59 billion in annual advertising revenue \citep{google2025}. On Amazon, sellers, on average, spent approximately \$298,000 per day on pay-per-click (PPC) ads in 2025 \citep{amazon12025}. 
A primary driver behind the rapid expansion of online advertising lies in its effectiveness in converting highly measurable ad clicks (through banners, search results, and social media posts) into sales in a cost-efficient manner, especially when compared to traditional advertising channels such as television, print, and billboards \citep{hu2016incentive,asadpour2019concise}.

\begin{figure}[htbp]   
    \centering  
    \includegraphics[width=0.85\textwidth]{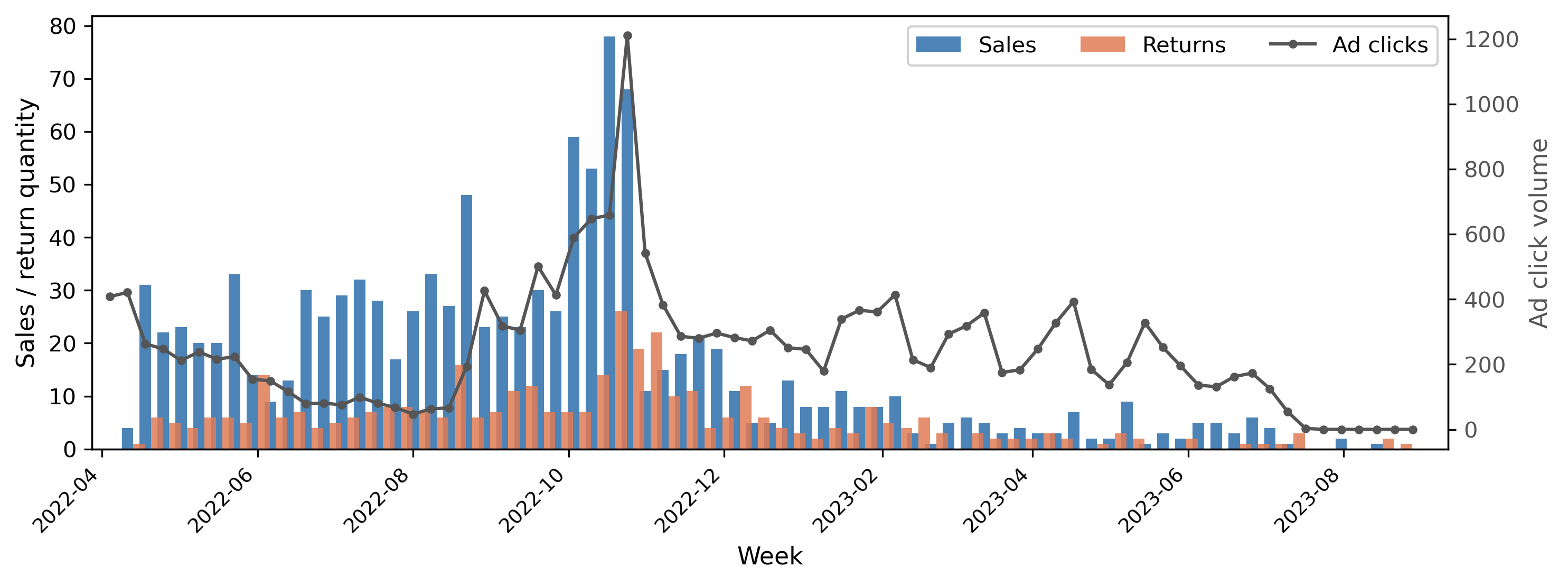}  \\
    \caption{Weekly sales and return quantities, and ad click volume of a representative product.}
    \label{fig::motivating_example}
\end{figure}

Online advertising is associated with both demand and product returns in the observed data. \Cref{fig::motivating_example} illustrates weekly sales, return quantities, and ad click volumes for a representative product sold by a
fashion retailer that extensively uses paid search ads on Amazon. The figure reveals that higher ad clicks are associated with increases in both sales and returns, suggesting that while advertising stimulates demand, it also amplifies product returns. This pattern may be explained by several consumer behavioral mechanisms. First, online advertisements often employ filters or AI-enhanced visuals that exaggerate product attributes such as fabric texture, color, or fit, particularly in apparel categories. These enhancements can widen the discrepancy between customers' pre-purchase expectations and the actual product received, leading to expectation mismatch. Second, when combined with lenient return policies such as free or hassle-free returns, online advertising may encourage impulse purchases \citep{feng2023modeling} or opportunistic returns \citep{altug2021managing} by framing transactions as low-risk trials, thereby increasing the likelihood of subsequent returns. 

The intertwined effects of online advertising on product sales and returns pose significant operational challenges for e-commerce retailers when coordinating inventory and advertising decisions. First, once the impact of advertising on product returns is accounted for, increased ad spending does not necessarily translate into higher net sales. In practice, retailers often incur significant reverse logistics costs from excessive returns. For example, under the Fulfillment by Amazon return policy, sellers are charged handling and processing fees for returned items \citep{amazonreturn2024}, which can materially erode margins and, in some cases, lead to negative profitability. Consequently, retailers must exercise greater caution in advertising decisions by proactively anticipating how advertising may inflate returns. Second, in practice, retailers typically plan their inventory levels and online advertising campaigns well ahead of selling seasons. As a result, ad clicks and the associated ad-driven demand are uncertain at the time of inventory and advertising planning. This sequencing amplifies volatility in both demand and returns, further complicating the effective coordination of marketing and operations decisions. Third, for retailers managing multiple products, advertising decisions exert dual effects on product demand and returns: \textit{direct effect}, which captures the impact of a product's own advertising on its sales and return quantities; and \textit{spillover effects}, whereby advertising for one product impacts the sales and returns of the retailer's other products \citep{kumar2015demand,zhang2023post,chen2025multichannel}. The interplay of these direct and spillover effects further complicates the relationship between online advertising, sales, and returns, substantially increasing the complexity of jointly managing advertising and inventory decisions.

Motivated by these challenges and by our industrial partner's practice, this study examines a retailer’s problem of jointly optimizing inventory and advertising decisions across multiple products prior to the selling season. The retailer decides the inventory level and whether to advertise each product, recognizing that advertising simultaneously affects demand and returns and induces both direct and cross-product spillover effects, all under substantial uncertainty. We formalize this decision problem as a contextual stochastic optimization (CSO) model, in which demand, returns, and ad clicks are stochastic and influenced by observable contextual features. Within this framework, we seek to address two core research questions:
\begin{enumerate}
    \item[(a)] How can we model product demand and returns in a way that explicitly captures both direct and spillover effects of advertising?
    \item[(b)] How can we efficiently optimize the retailer’s joint inventory and advertising decisions in the presence of these intertwined effects?
\end{enumerate}

To answer the first question, we decompose demand into net demand and product returns and propose linear models for both components. In the formulation, net demand and returns depend linearly on advertising decisions, \textit{observable} covariates such as prices and product attributes, and \textit{unobservable} ad clicks---that is, click volumes are realized only after the advertising decision is made and are therefore treated as unobservable covariates at decision time throughout the paper. The linear structure captures both direct advertising effects (own ad clicks on own demand) and cross-product spillover effects (other products' ad clicks on a focal product's demand and returns), while preserving the tractability required for the subsequent optimization. Section~\ref{subsec::empirical evidences} estimates these models on our industrial partner's data and compares them with alternative specifications. 

To address the second question, we first consider the full-information setting where the model parameters are known and ad clicks and error terms jointly follow a multivariate normal distribution. The optimal inventory is then a newsvendor quantity whose safety stock depends on which products are advertised, and the optimal advertising decision solves a supermodular maximization problem, which for two products reduces to a threshold policy in the mean ad clicks (\Cref{prop:prel_result}). Supermodular maximization is tractable without side constraints but computing an optimal advertising vector becomes challenging under a cardinality constraint \citep{karp1972reducibility}. More importantly, both the distribution and the parameters that these results presume known must in practice be inferred from a limited history of past seasons.

Neither presumption holds in practice, which raises two difficulties. First, \textit{distributional ambiguity}: the joint distribution of ad clicks and error terms is unknown. Second, \textit{decision-dependent prediction}: the advertising decision determines which potential ad clicks are realized and therefore enter the predictive model as a covariate, so the induced distribution of net demand and returns shifts with the decision. Together with the combinatorial structure identified above, these challenges preclude a direct solution of the contextual stochastic optimization model and motivate a data-driven formulation. The conventional estimate-then-optimize (ETO) approach estimates the model parameters from historical data, forms an empirical distribution from the regression residuals, and maximizes expected profit under it (\textit{e.g.}, \citealt{kannan2025data}). Its appeal is computational, but it can be optimistic for two reasons: the optimizer favors decisions flattered by sampling noise, the optimizer's curse of \citet{smith2006optimizer}, and the in-sample residuals it optimizes against can understate out-of-sample error. Both effects can be more pronounced as data grow scarce (Section~\ref{subsec:ETO}).

To address these limitations, we propose an estimate-then-robustify (ETR) approach that builds on the robust satisficing framework of \citet{long2023robust} and the residuals-based approach of \citet{kannan2023residuals}. Rather than committing to the empirical distribution, ETR seeks a decision that attains a given profit target and whose profit deteriorates as slowly as possible as the true distribution moves away from the empirical one, with the movement measured in the type-$1$ Wasserstein distance; the target replaces the ambiguity radius that distributionally robust formulations must calibrate. Both difficulties shape the analysis that follows: the performance guarantee must accommodate unobservable, decision-dependent covariates, and the reformulation must handle the bilinear coupling between advertising decisions and random ad clicks.

Our contributions are fourfold.
\begin{enumerate}
\item \textit{Modeling.} We propose linear net-demand and return models in which advertising acts through realized ad clicks, with a direct effect for each product and a spillover effect for each ordered pair of products (\Cref{assm:demand and return model}), and we embed them in a contextual stochastic optimization model,~\ref{model::CSO}, that jointly determines inventory and advertising before the selling season (Section~\ref{sec:IPAP}). Under full information, the optimal inventory is a newsvendor quantity whose safety stock depends on the advertising decision, and the optimal advertising decision solves a supermodular maximization problem that reduces to a threshold policy for two products (\Cref{prop:prel_result}). On the retailer's data, the linear specification attains the lowest out-of-sample prediction error among five candidate specifications (Section~\ref{subsec::empirical evidences}).

\item \textit{Statistical guarantees.} For the estimate-then-robustify model~\ref{model::ETR}, we establish a finite-sample guarantee on the out-of-sample expected profit that holds for any parameter estimate, with an explicit high-probability bound for least squares under strict exogeneity (\Cref{thm:out-of-sample}, \Cref{cor:hp}), and we show that the ETR solution is asymptotically optimal for~\ref{model::CSO} when the estimator is consistent and the target approaches the empirical optimal profit (\Cref{prop:asymptotic}). These results extend the residuals-based guarantees of \citet{kannan2023residuals} and \citet{zhu2024residuals} to covariates that are unobservable at the time of decision and enter the prediction through the decision itself, and to the target-based robust satisficing criterion, under which~\ref{model::ETR} minimizes fragility.

\item \textit{Computation.} We reformulate~\ref{model::ETR} exactly as a mixed-integer conic program whose size grows as $\mathcal{O}(2^N S N)$ in the number of products $N$ and sample size $S$ (\Cref{thm::P_refo}), and we derive by affine recourse adaptation an approximation,~\ref{model::bar_ETR}, of size $\mathcal{O}(S N^2)$ that is exact for a single product and inherits the finite-sample guarantee (\Cref{Theo::LDR_refo}, \Cref{prop::optimality of single product}). In our experiments, \ref{model::bar_ETR} matches the exact solution to within 0.0003\% in training profit while solving up to 38 times faster (Section~\ref{subsec:experiments}).

\item \textit{Empirical evidence.} Our case study shows that \ref{model::bar_ETR} improves out-of-sample profit over estimate-then-optimize by 6.17\% and earns the higher profit in 74.5\% of instances, with the largest gains when training data are scarce. Without spillover effects in the decision model, the gain falls to 1.02\%. 
The percentage profit gain from robustification is 5.73\% when the models are re-estimated by instrumental variables to address possible endogeneity of ad clicks (Sections~\ref{subsec:experiments} and~\ref{sec:endo_extension}).
\end{enumerate}

\subsection{Related Work}

This study relates most closely to four streams of literature: online advertising and retail operations, product returns, robust and joint inventory--advertising decisions, and residuals-based data-driven optimization. 

In the context of online retailing, a substantial body of work examines firms’ advertising decisions and their performance implications; see \cite{choi2020online} for a comprehensive review. Early studies investigate firms’ choices among advertising pricing schemes, such as cost-per-click, cost-per-impression, and cost-per-action \citep{asdemir2012pricing, hu2016incentive}, as well as the determinants of ad ranking and effectiveness in sponsored search markets \citep{ghose2009empirical, du2017bidding}. From a single-product or aggregate perspective, prior research analyzes a range of advertising decision problems, including bidding optimization \citep{ye2015sponsored, asadpour2019concise}, budget allocation \citep{tan2005allocating}, and the choice between spending-based and traffic-based decision processes \citep{liu2018advertising}. Nevertheless, this stream largely abstracts away from product returns, which erode the marginal value of advertising and distort optimal advertising levels, particularly in multiproduct settings where spillover effects further complicate decision making. 

A related stream examines cross-product advertising spillovers (see, \textit{e.g.}, \citealt{kumar2015demand, zhang2023post, chen2025multichannel}). These studies document how advertising for one product affects the demand for others but leave its effect on product returns largely unexplored. 

Product returns are a critical characteristic of modern e-commerce retailing and have attracted considerable attention in the operations management literature \citep{abdulla2019taking}. We focus on consumer-initiated returns. Existing studies typically examine how returns affect inventory decisions \citep{ketzenberg2009optimal, oh2024behavioral}, pricing \citep{shulman2009optimal, hu2019dynamic}, or assortment planning \citep{alptekinouglu2014carry}. Limited attention has been devoted to how product returns interact with advertising. 
Advertising, where it appears, is treated as exogenous rather than jointly optimized with inventory or pricing in the presence of returns. Related work also
jointly studies customer influence, pricing, and return policies \citep{huang2026influencing}. Our focus is joint inventory and advertising across products, with spillovers in demand and returns and uncertainty characterized from data.

Relatedly, the prediction of product returns has followed two main approaches. One focuses on estimating return rates using product characteristics or auxiliary information \citep{dzyabura2023leveraging}. Another directly predicts return volumes from historical data using statistical or machine-learning methods \citep{shang2020using, cui2020predicting}. However, these studies generally treat returns as outcomes that are independent of firms’ marketing actions. Our return model instead lets returns depend on advertising through both direct and spillover effects. 


Joint optimization of inventory and advertising (or sales effort) has also been studied in the literature. Classic examples include \citet{sogomonian1993modeling}, who coordinate promotion timing and production decisions within a firm, and \citet{cheng1999periodic}, who characterize optimal inventory policies when demand is influenced by the firm's promotion decisions; game-theoretic treatments of advertising in supply chains include the cooperative advertising and pricing model of \citet{he2009cooperative}. These studies adopt stylized game-theoretic or stochastic-control formulations and do not incorporate data-driven uncertainty characterization. 
This paper studies the joint inventory and advertising problem in a multiproduct setting under a data-driven robust optimization framework that explicitly accounts for advertising spillover effects across products and their impact on both demand and returns.

Methodologically, our estimate-then-robustify framework builds on two lines of work. Residuals-based distributionally robust optimization (DRO), introduced by \citet{kannan2023residuals}, centers a Wasserstein ambiguity set at the empirical distribution of regression residuals and establishes finite-sample performance guarantees under exogenous covariates. \citet{zhu2024residuals} extend it to uncertainty that depends on the decision as well as on the covariates, with decision-dependent nominal distributions and ambiguity radii, and establish asymptotic optimality, rates of convergence, and finite-sample guarantees for Wasserstein, sample-robust, and $\phi$-divergence ambiguity sets. Robust satisficing \citep{long2023robust} replaces the calibrated ambiguity radius by a profit target and minimizes the rate at which the profit guarantee deteriorates with the distance from the empirical distribution; \citet{sim2025analytics} bring covariates into this framework but treat them as observable and exogenous. Our setting differs from all four in that the ad clicks are covariates that are unobservable when the decision is made and enter the prediction through the decision itself, and our guarantees extend robust satisficing to this case of decision-dependent prediction with unobservable covariates.

\subsection{Notation}\label{subsec:organization and notation}
Let $[N]=\{1,2,\ldots,N\}$, and, for an index set $\mathcal{I} \subseteq [N]$, let $[N]\setminus\ci$ denote the indices in $[N]$ that are not in $\ci$. We use boldface symbols, such as $\bm{x}\in \mathbb{R}^N$ and $\mf{A}\in \mathbb{R}^{M\times N}$ to denote vectors (column, by default) and matrices, and we denote by $x_n$ the $n$-th element of vector~$\bm{x}$, by $\bm a_m$ the transpose of the $m$-th row vector of matrix~$\mf{A}$, by ${a}_{mn}$ the $(m,n)$-th entry of matrix~$\mf{A}$. Special vectors and matrices of the appropriate dimension include $\bm{0}$, $\bm{1}$, $\bm{e}_n$, and $\mf I$, which correspond to (respectively) the vector of all zeros, the vector of all ones, the $n$-th standard basis vector, and the identity matrix.
We use $\bmt{z} \sim \mathbb{P} \in \mathcal{M}(\mathcal{Z})$ to denote a random vector~$\bmt{z}$ governed by a probability distribution~$\mathbb{P}$, where $\mathcal{M}(\mathcal{Z})$ represents the set of all probability distributions in~$\mathcal{Z}$ with finite first moments. For a set $ \mathcal{Z} \subseteq\mathbb{R}^N$, the term $\mathbb{P}[\bmt{z}\in\mathcal{Z}]$ represents the probability of~$\bmt{z}$ lying in the set~$\mathcal{Z}$ evaluated on the probability distribution~$\mathbb{P}$. We use $\mathbb{E}_{\mathbb{P}}[\cdot]$ to signify the expectation with respect to $\mathbb{P}$, and use $\text{Var}[\cdot]$ to represent the marginal variance. Additionally, $\bm x \odot \bm y$ is a vector of the same dimension as the operands, with elements given by $(\bm x \odot \bm y)_n = x_n y_n$. 
Norm operator is denoted by $\| \cdot \| $ and $\| \cdot \|_*$ represents the dual norm. We denote by $\|\cdot\|_{\rm F}$ the Frobenius norm of a matrix. For a positive definite matrix $\bm \Sigma$, we define the matrix norm $\|\mf A\|_{\bm \Sigma} = \sqrt{\trace{\mf A \bm \Sigma \mf A^\top}}$ for a matrix $\mf A$ with as many columns as $\bm \Sigma$ has rows, and the vector norm $\|\bm v\|_{\bm \Sigma} = \sqrt{\bm v^\top \bm \Sigma \bm v}$. We denote by $\mathrm{Diag}(\bm x) \in \bbr^{N\times N}$ a diagonal matrix with entries from a vector $\bm x \in \bbr^N$ and by $\mathrm{diag}(\mf X) \in \bbr^{N}$ a vector whose entries are the diagonal entries of $\mf X \in \bbr^{N \times N}$.


\section{Inventory Planning and Advertising Problem}\label{sec:IPAP}
This section formulates the joint inventory planning and advertising problem, based on the practice from our industrial partner, \textit{i.e.}, the aforementioned retailer that sells fashion products on Amazon.

We consider a retailer selling a set of $N$ products to a (possibly overseas) market over a short selling season, typically lasting only a few weeks. Products are pre-positioned in an onshore warehouse (\textit{e.g.}, a Fulfillment-by-Amazon (FBA) warehouse in the case of our industrial partner) near the target market to enable timely order fulfillment, and each unit carries a procurement cost $c_n (>0)$. During the season, inventory replenishment is not possible due to the prolonged lead times; for example, shipping from China to the U.S. typically takes around 40 days \citep{stanley2025the}. As the selling season progresses, the stochastic demand $\tilde{d}_n$ for product $n$ is realized as $d_n$ and fulfilled at the retail price $p_n (>0)$ using the pre-positioned inventory. Leftover units at the end of the season incur a holding cost $h_n (>0)$ per unit. Conversely, if available inventory is insufficient to meet realized demand, the unmet portion is backlogged at a per-unit penalty cost $\alpha_n (>0)$. Our industrial partner routinely adopts this practice, and it is common in cross-border operations; for example, when inventory at an Amazon warehouse is depleted, backlogged orders may be fulfilled by shipping from an offshore warehouse via international express delivery \citep{shi2020analytics}.

Customers may choose to return purchased items for various reasons, including size or color mismatches, preference changes, or dissatisfaction upon receipt. A return of product $n$ requires the retailer to issue a full refund to the customer and additionally pay a return-processing fee $\gamma_n (>0)$ to the platform for reverse logistics, quality inspection, etc. Returns are operationally complex and time-consuming: data from our industrial partner reveals that it typically takes about 35 days for returned items to be transported back to Amazon warehouses. As a result, the returned unit may not be resold within the same selling season. Accounting for returns, we denote $\tu{d}_n$ as the net demand for product $n$ so that $\tu{d}_n = \tilde {d}_n - \tilde r_n$, where $\tilde{r}_n$ represents the number of returned units. 

Before the selling season starts, the retailer jointly determines inventory levels and advertising decisions. Let $x_n$ denote the inventory level selected for product $n\in[N]$, and let $\bm x \in \mathcal{X} \subseteq \mathbb{R}_+^N$ denote the inventory vector. Concurrently, the retailer chooses which products to advertise. We focus on automatic PPC advertising, a widely used mechanism on platforms such as Amazon and Google. Under this mechanism, the retailer does not need to manually select keywords and manage bids \citep{Insights2022}, allowing it to focus on the core decision of selecting which products to advertise. Let $y_n \in \{0,1\}$ denote the advertising decision for product $n$. When $y_n = 1$, the retailer promotes product $n$, generating a realized number $\xi_n$ of ad clicks with an advertising expenditure of $k_n\xi_n$, where $k_n (>0)$ denotes the average cost per click. Conversely, when $y_n = 0$, no promotion is undertaken for product $n$. For notational convenience, we denote the advertising decision vector by $\bm y \in \mathcal{Y} \subseteq \{0,1\}^N$, the realized ad clicks vector by $\bm{\xi} = (\xi_1, \ldots, \xi_N)^\top \in \bbr_+^N$, and the corresponding random vector by $\tilde{\bm{\xi}}$.

Advertising influences consumer behavior by generating ad clicks, a fraction of which convert into demand, with a portion of that demand ultimately returned. Accordingly, both net demand and returns are affected by advertising through the sequential processes of click-to-demand and demand-to-return conversion. In addition to advertising, both the retailer's net demand and product return may depend on a variety of external factors, including retail prices, product attributes, fixed time effects, and weather conditions. We represent the covariates associated with the net demand $\tu{d}_n$ for product $n$ as a vector  $\bmu v \in \uline{\mathcal{V}} \subseteq \mathbb{R}^{J_1}$. Likewise, the covariates relevant to return $\tilde{r}_n$ are denoted by a vector $\bm v \in \mathcal{V} \subseteq \mathbb{R}^{J_2}$.

To capture the combined effects of advertising and covariates on the retailer's net demand $\tu{d}_n$ and return $\tilde{r}_n$, we impose the following model assumption: 
\begin{assm}\label{assm:demand and return model}
    For each product $n \in [N]$, net demand (resp. returns) is modeled as a linear function of its associated covariates $\bmu v$ (resp. $\bm v$), and advertising exposures $\bm y \odot \bmt \xi$, given by:
    \begin{align}
        \tu{d}_n & = (\bmu a_n)^{\top} \bmu v + (\bmu b_n)^{\top} ( \bm y \odot \bmt \xi) + \tu{\epsilon}_n, \label{eq:nd}\\
        \tilde{r}_n &= (\bm a_n)^{\top} \bm v + (\bm b_n)^{\top} ( \bm y \odot \bmt \xi) + \tilde{\epsilon}_n.\label{eq:return}
    \end{align}
    Here, $\bmu a_n \in \mathbb{R}^{J_1}$ (resp. $\bm a_n \in \mathbb{R}^{J_2}$) and $\bmu b_n \in \mathbb{R}^N$ (resp. $\bm b_n \in \mathbb{R}^N$) are deterministic parameters. 
    The random vector $\bmt z := (\bmt \xi, \bmtu \epsilon, \bmt \epsilon)$ has a joint distribution $\mathbb P^\star \in \mathcal M(\cz)$---where $\cz \subseteq \mathbb R^{3N}$---that does not depend on $\bm y$. 
\end{assm}

In the net-demand function \eqref{eq:nd} (\textit{resp.}, return function \eqref{eq:return}), the parameter $\uline{b}_{nn}$ (\textit{resp.} $b_{nn}$) captures the direct effect of product $n$'s ad clicks on its own net demands (\textit{resp.} returns), while $\uline{b}_{nm}$ (\textit{resp.} $b_{nm}$) for $m\neq n$ reflects the spillover effects of product $m$'s ad clicks on product $n$'s net demand (\textit{resp.} return). We express net demand and returns in compact vector form across all products:
\begin{align}\label{model:demand vector}
    \bmtu{d} & = \mfu A\, \bmu v + \mfu B ( \bm y \odot \bmt \xi) + \bmtu{\epsilon},\\
    \label{model:return vector}
    \bmt{r} & = \mf A\, \bm v + \mf B ( \bm y \odot \bmt \xi) + \bmt{\epsilon},
\end{align}
where $\mfu A \in \uline{\mathcal{A}} \subseteq \mathbb R^{N \times J_1}$, $\mf A \in \mathcal{A} \subseteq \mathbb R^{N \times J_2}$, $\mfu B \in \uline{\mathcal{B}} \subseteq \mathbb R^{N \times N}$, $\mf B \in \mathcal{B} \subseteq \mathbb R^{N \times N}$ are the parameter matrices associated with covariates and advertising effects. The vectors $\bmtu{\epsilon},\bmt{\epsilon}$ denote the random errors in net demand and return, respectively. Notably, the random ad clicks $\tilde{\bm \xi}$ can also be regarded as covariates. In this sense, the prediction of net-demand and return involves both \textit{observable} covariates ($\bmu v$ and $\bm v$) and \textit{unobservable} covariates ($\tilde{\bm \xi}$), which substantially complicates the statistical analysis and computational approximation as we discuss in detail in Section~\ref{sec:ETR}. We empirically validate the linear specification of the net-demand and return models using real-world data in Section~\ref{subsec::empirical evidences}. Moreover, the invariance of $\mathbb P^\star$ to $\bm y$ is motivated by the hypothesis that, under automatic PPC, impressions are allocated through platform-wide keyword auctions in which the retailer's own products are a negligible fraction of the competing advertising inventory, so activating one more product draws impressions from competing sellers rather than from the retailer's other campaigns, and cross-product interactions arise downstream, through the spillover coefficients $\uline{b}_{nm}$ and $b_{nm}$. As a descriptive check, conditional on own campaign spend, a product's weekly clicks are essentially unrelated to the number of other focal products advertised in the same week (see Appendix~\ref{app:click_exogeneity} for a test of this invariance).

Accordingly, the total demand is given by $\bmt d = \bmtu{d} + \bmt{r}$ and total profit can thus be expressed as:
\begin{align*}
    \bm p^\top \bmt d- \big(\bm h^\top(\bm x-\bmt d)^+ + \bm \alpha^\top(\bmt d-\bm x)^+\big) - \bm c^\top \bm x - (\bm p + \bm \gamma)^\top \bmt r - \bm k^\top (\bm y \odot \bmt\xi),
\end{align*}
in which the term $\bm p^\top \bmt d$ is the revenue from product sales, $\bm h^\top(\bm x-\bmt d)^+ + \bm \alpha^\top(\bmt d-\bm x)^+$ accounts for inventory holding costs and backorder penalties,  $\bm c^\top \bm x$ is the inventory procurement cost,   $(\bm p + \bm \gamma)^\top \bmt r$ covers full refund and per-unit return-processing fee for each returned unit, and $\bm k^\top (\bm y \odot \bmt\xi)$ is the advertising expenditure. 

For any realization $(\bm \xi,\bmu \epsilon,\bm \epsilon)$, the profit function can be written explicitly as:
\begin{equation}\label{eq: profit function def}
    \begin{aligned}
        \Pi(\mfu A, \mfu B, & \mf A, \mf B,  \bmu v, \bm v, \bm x, \bm y, \bm \xi, \bmu \epsilon, \bm \epsilon) 
        :=(\bm p - \bm \alpha)^\top \big(\mfu A\, \bmu v + \mfu B ( \bm y \odot \bm \xi) + \bmu{\epsilon}\big) 
        - (\bm \alpha + \bm \gamma)^\top \big(\mf A\, \bm v + \mf B ( \bm y \odot \bm \xi) + \bm{\epsilon}\big) \\
        &- (\bm h + \bm \alpha)^\top \big(\bm x - \mfu A\, \bmu v - \mf A\, \bm v - (\mfu B + \mf B) ( \bm y \odot \bm \xi) - \bmu{\epsilon} - \bm{\epsilon}\big)^+ + (\bm \alpha - \bm c)^\top \bm x - \bm k^\top(\bm y \odot \bm \xi).
    \end{aligned}
\end{equation}

Observing the covariates $(\bmu v, \bm v)$, the retailer chooses the two here-and-now decision variables---the inventory vector $\bm x \in \mathcal{X}$ and the binary advertising-activation vector $\bm y \in \mathcal{Y}$---before the ad clicks and error terms $(\bmt \xi, \bmtu \epsilon, \bmt \epsilon)$ are realized. The joint inventory planning and advertising problem can thus be formulated as the following mixed-integer contextual stochastic optimization (CSO) problem:
\begin{equation}\label{model::CSO}\tag{$\code{CSO}$}
    \begin{array}{c@{\;}l@{\;}l@{\;}}
        \max~& \mathbb{E}_{\mathbb P^\star} 
        \Big[
        \Pi(\mfu A, \mfu B, \mf A, \mf B, \bmu v, \bm v, \bm x, \bm y, \bmt \xi, \bmtu \epsilon, \bmt \epsilon) 
        \Big] \\
        ~	{\rm s.t.} ~ &\bm x \in \mathcal{X}, \; \bm y \in \mathcal{Y}.
    \end{array}
\end{equation}

In~\ref{model::CSO}, inventory and advertising are modeled as concurrent decisions. This modeling choice reflects three features of the retailer's operations. First, inventory must be committed before the selling season because of long replenishment lead times. Advertising decisions are also planned in advance because they affect product-level demand, returns, and cross-product spillovers. Inventory and advertising decisions are therefore interdependent and are optimized jointly. Second, adjusting advertising based on inventory levels may not be desirable. For example, a product with high inventory may not be a suitable advertising target if its ads generate demand with high return rates, while advertising a product with favorable cross-product spillovers may be beneficial even when its inventory is limited. Joint optimization accounts for these effects. Third, our industrial partner manages automatic PPC advertising at the campaign level rather than through continuous bid adjustments. Revising campaigns during the selling season also requires additional setup and management effort, which limits the frequency of such adjustments. 

\section{Preliminaries}\label{sec:vali and prel}

\subsection{Full Information Setting}\label{sec:preliminary}

We first analyze the~\ref{model::CSO} model under full information, assuming that the model parameters $(\mfu A, \mfu B, \mf A, \mf B)$ and the joint distribution $\mathbb{P}^\star$ are known. We decompose the random ad clicks as $\bmt \xi=\bm \mu+\bmt \zeta$, where $\bm \mu = \mathbb{E}_{\mathbb{P}^\star}[\bmt \xi] > \bm 0$ denotes the expected number of ad clicks, and $\bmt\zeta$ is a zero-mean random error term, \textit{i.e.}, $\mathbb{E}_{\mathbb{P}^\star}[\bmt \zeta] = \bm 0$. For analytical tractability, we assume  $\mathcal{X} = \mathbb{R}_+^N$ and that the joint distribution $\bar{\bbp}^\star$ of $(\bmt\zeta,\bmtu \epsilon,\bmt \epsilon)$ is known to follow a multivariate normal distribution with independent components and zero-mean residual blocks, $\mathbb{E}_{\mathbb{P}^\star}[\uline{\tilde{\bm \epsilon}}] = \mathbb{E}_{\mathbb{P}^\star}[\tilde{\bm \epsilon}] = \bm 0$---meaning all off-diagonal entries of the correlation matrix are zero. We further assume $\bm \alpha > \bm c$ componentwise, ensuring that expediting a unit costs more than pre-positioning it and ruling out the degenerate regime in which the retailer optimally stocks nothing and backlogs all demand \citep{zipkin2000foundations, bertsimas2006robust}, and that mean demand is large enough for the inventory in~\Cref{prop:prel_result}(i) below to be nonnegative for every $\bm y \in \cy$. Let $\Phi(\cdot)$, $\phi(\cdot)$, and $\Phi^{-1}(\cdot)$ denote the standard normal CDF, PDF, and inverse CDF, respectively, and they are applied element-wise when their arguments are vectors. We further define the marginal variances $\text{Var}[\bmt \zeta]=\uwave{\bm\sigma}^2,\; \text{Var}[\bmtu\epsilon]=\bmu \sigma^2$, and $\text{Var}[\bmt\epsilon]=\bm \sigma^2$. We also define $\bm{\Upsilon}(\bm y)$ as the aggregated stochastic deviation, with each component given by:
$\Upsilon_n(\bm y) := \sqrt{(\bmu b_n + \bm b_n)^\top \mathrm{Diag}( \bm y \odot \uwave{\bm\sigma}^2) (\bmu b_n + \bm b_n) + \uline{\sigma}_n^2+\sigma_n^2}$ for $n \in [N]$. 
Based on these definitions and assumptions, we derive the following results. 

\begin{prop}\label{prop:prel_result}
    In the full information setting,
    
    (i) Given $\bm y \in \mathcal{Y}$, the retailer's optimal inventory decision is
    $$
    \bm x^\star(\bm y)=\mfu A\, \bmu v +\mf A\, \bm v + (\mfu B + \mf B) ( \bm y \odot \bm \mu)  +\bm{\Upsilon}(\bm y) \odot \Phi^{-1}\Big(\frac{\bm \alpha-\bm c}{\bm \alpha+ \bm h}\Big),
    $$
    where the last term $\bm{\Upsilon}(\bm y) \odot \Phi^{-1}\Big(\frac{\bm \alpha-\bm c}{\bm \alpha+ \bm h}\Big)$ is the safety stock.
    
    (ii) The optimal advertising decision $\bm y^\star$ is a maximizer of a supermodular optimization problem:
    \begin{align}\label{prob:ad}
        \max_{\bm y \in \mathcal{Y}}\left\{\left((\bm p- \bm c)^\top \mfu B-(\bm\gamma+\bm c)^\top \mf B - \bm k^\top\right) (\bm y \odot \bm \mu)-\bm \Upsilon(\bm y)^\top \left((\bm \alpha+\bm h) \odot \phi\left(\Phi^{-1}\left(\frac{\bm \alpha-\bm c}{\bm \alpha+\bm h}\right)\right) \right) \right\},
    \end{align}
    where $(\bm \alpha+\bm h) \odot \phi\big(\Phi^{-1}(\frac{\bm \alpha-\bm c}{\bm \alpha+\bm h})\big) > \bm 0$ is the marginal inventory cost factor.
    
    (iii) In the special case with $N=2$ and $\mathcal{Y} = \{0,1\}^2$, the optimal advertising decision $\bm y^\star=(y^\star_1,y^\star_2)$ follows the threshold-based policy:
    $\bm y^\star=(1,1)$ if $\iota_1\mu_1\geq \psi_{01}$, $\iota_2\mu_2\geq \psi_{10}$, and $\iota_1\mu_1+\iota_2\mu_2\geq \psi_{10}+\hat\psi_{10}$; $\bm y^\star=(1,0)$ if $\iota_1\mu_1\geq \hat\psi_{10}$, $\iota_2\mu_2\leq \psi_{10}$, and $\iota_1\mu_1-\iota_2\mu_2\geq \psi_{01}-\psi_{10}$; $\bm y^\star=(0,1)$ if $\iota_1\mu_1\leq \psi_{01}$, $\iota_2\mu_2\geq \hat\psi_{01}$, and $\iota_1\mu_1-\iota_2\mu_2\leq \psi_{01}-\psi_{10}$; and $\bm y^\star=(0,0)$ if $\iota_1\mu_1\leq \hat\psi_{10}$, $\iota_2\mu_2\leq \hat\psi_{01}$, and $\iota_1\mu_1+\iota_2\mu_2\leq \psi_{10}+\hat\psi_{10}$.
    Here, $\iota_n:=\sum_{m\in\{1,2\}}[(p_m-c_m)\uline b_{mn}-(\gamma_m+c_m)b_{mn}]-k_n$ for $n\in\{1,2\}$, and $(\psi_{y_1,y_2},\hat\psi_{y_1,y_2})$ are parameters defined in Appendix \ref{app:prel-col}. 
\end{prop}
\textit{Proof.} We relegate the proofs of all analytical results to the Appendix~\ref{appendix:Proofs}. $\hfill\square$

Part (\textit{i}) of \Cref{prop:prel_result} provides a closed-form expression for the optimal inventory decision for any advertising decision $\bm y \in \mathcal{Y}$. The first three terms represent the mean demand, incorporating both covariate effects and expected ad-driven demand, while the final term is the safety stock, which balances the risk of overstocking and understocking. This safety stock is jointly determined by the (element-wise) targeted service level $(\bm \alpha - \bm c)/(\bm \alpha + \bm h)$ and the aggregated deviation $\bm \Upsilon(\bm y)$. This structure generalizes the classical newsvendor problem to a multivariate setting where demand is endogenously dependent on advertising decisions.

Part (\textit{ii}) shows that the \ref{model::CSO} model reduces to a supermodular maximization problem over the binary feasible region. The objective captures a trade-off between expected incremental net profit from advertising (via direct and spillover effects encoded in $\mfu B$ and $\mf B$) and a penalty associated with increased inventory risk (via $\bm \Upsilon(\bm y)$). The supermodular nature of the problem reflects complementarities among advertising decisions: promoting one product may mitigate the marginal increment in inventory risk caused by promoting other products (reflected by the submodularity of $\bm \Upsilon(\bm y)$), thereby enhancing the overall marginal benefit of joint promotion. However, the computational complexity of Problem~\eqref{prob:ad} depends critically on the presence of additional combinatorial constraints. Standard practical constraints such as ``a constraint on the number of advertised products'' render the problem generally NP-hard (see e.g., \citealt{karp1972reducibility, shu2005stochastic}). 

\begin{figure}[ht]
    \centering
    \vspace{-0.3cm}
    \includegraphics[width=0.9\textwidth]{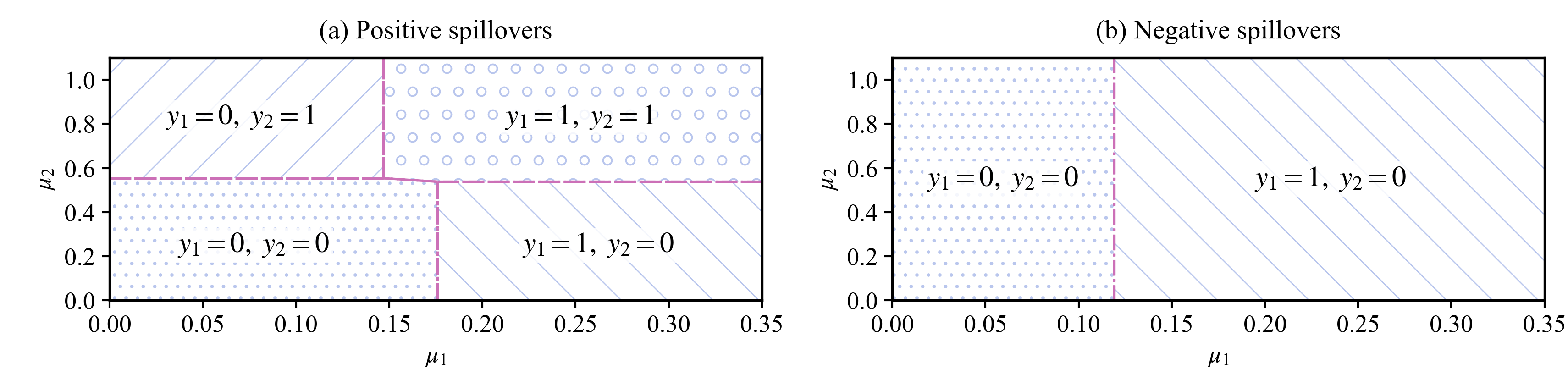}
    \caption{Optimal advertising strategy for two products in the $(\mu_1,\mu_2)$ parameter space, where $\bm p = (3,4)^\top, \bm c = 0.25 \bm p, \bm h = 0.06 \bm p, \bm \alpha = 0.30 \bm p, \bm \gamma = 0.20 \bm p, \bm k = (0.004,0.005)^\top, \uwave{\bm \sigma} = (2,3)^\top, \bmu \sigma = (1,2)^\top, \bm \sigma = (2,1)^\top$. The left panel depicts positive spillover effects with $\mfu B = \big(\begin{smallmatrix} 0.2 & 0.3 \\ 0.1 & 0.2 \end{smallmatrix}\big)$ and $\mf B = \big(\begin{smallmatrix} 0.1 & 0.2 \\ 0.2 & 0.3 \end{smallmatrix}\big)$; the right one depicts negative spillover effects with $\mfu B = \big(\begin{smallmatrix} 0.2 & -0.3 \\ -0.1 & 0.2 \end{smallmatrix}\big)$ and $\mf B = \big(\begin{smallmatrix} 0.1 & -0.2 \\ -0.2 & 0.3 \end{smallmatrix}\big)$.}
    \label{fig:stra}
    \vspace{-0.3cm}
\end{figure}

By contrast, in the unconstrained case, the problem can be solved in polynomial time by combinatorial algorithms (\textit{e.g.}, \citealt{iwata2001combinatorial}). Part (\textit{iii}) illustrates this structure in the two-product case by providing an explicit threshold-based binary advertising policy $\bm y^\star = (y_1^\star, y_2^\star)$. The policy relies on adjusted metrics $\iota_1\mu_1, \iota_2\mu_2$ (net marginal benefits of advertising each product, accounting for costs, revenue, and returns) compared against thresholds $\psi_{10}, \hat\psi_{10},\psi_{01}, \hat\psi_{01}$, partitioning the decision space into four optimality regions that may overlap at ties, that correspond to the four binary configurations: advertise both, neither, or a single product. These regions are illustrated in \Cref{fig:stra}.

\subsection{Estimate-then-Optimize Approach}\label{subsec:ETO}

Notably, \Cref{prop:prel_result} relies on an idealized setting in which the retailer has complete knowledge of both the model parameters $(\mfu A, \mfu B, \mf A, \mf B)$ and the joint distribution $\mathbb{P}^\star$ of the~\ref{model::CSO} model. In practice, however, 
even with well-estimated parameters, the retailer faces substantial uncertainty about future realizations of ad clicks, demand shocks, and return shocks. Additionally, the problem involves \textit{decision-dependent prediction}, as both net demand $\uline{\tilde{\bm d}}$ and return $\tilde{\bm r}$ depend on the advertising decision $\bm y$; this characteristic may give rise to complex bilinear or non-convex optimization problems \citep[see, \textit{e.g.},][]{ryu2025nurse} that are computationally challenging to solve. 
    
In practice, retailers typically have access to a historical dataset $\mathcal{D}$, consisting of $S$ samples. Each sample includes covariates $\bmhu v_s$ and $\bmh v_s$, advertising decision $\bmh y_s$, observed potential ad clicks $\bmh \xi_s$, net demand $\bmhu d_s$, and return $\bmh r_s$. We denote this dataset as $\mathcal{D} = \{ (\bmhu v_s, \bmhu d_s, \bmh v_s, \bmh r_s, \bmh y_s, \bmh \xi_s)\}_{s \in [S]}$. To address the challenges of unknown parameters and distribution inherent in solving~\ref{model::CSO}, the retailer can employ the classical two-step ``estimate-then-optimize" approach, detailed as follows:

\noindent\textbf{Estimation.} The parameters $(\mfu A, \mfu B, \mf A, \mf B)$ can be estimated using the least squares method: 
\begin{align*}
    (\mfhu A,\mfhu B):=\argmin_{\mfu A, \mfu B}~ \sum_{s \in [S]} \| \bmhu d_s - \mfu A \, \bmhu v_s - \mfu B (\bmh y_s \odot \bmh \xi_s) \|_2^2,\; 
    (\mfh A ,\mfh B):=\argmin_{\mf A, \mf B}~ \sum_{s \in [S]} \| \bmh r_s - \mf A \, \bmh v_s - \mf B (\bmh y_s \odot \bmh \xi_s) \|_2^2.
\end{align*}
The resulting regression residuals are given by:  
$$
\bmhu \epsilon_{s} \coloneqq \bmhu d_{s} - \mfhu A \, \bmhu v_s - \mfhu B (\bmh y_s \odot \bmh \xi_s),~ \bmh \epsilon_{s} \coloneqq \bmh r_{s} - \mfh A \, \bmh v_s - \mfh B (\bmh y_s \odot \bmh \xi_s).
$$
These residuals capture the discrepancies between the observed data and the fitted model predictions and are used to construct an empirical distribution of $(\bmt \xi, \bmtu{\epsilon}, \bmt{\epsilon})$ as follows:
\[
\hat{\mathbb{P}} := \frac{1}{S}\sum_{s\in[S]} \delta_{(\bmh \xi_s, \bmhu{\epsilon}_s, \bmh{\epsilon}_s)}.
\] 
Here $\delta_{\bm z}$ denotes the Dirac probability measure at $\bm z$, and the support set is chosen so that $\hat{\mathbb{P}}\in \mathcal{M}(\cz)$.

\noindent\textbf{Optimization.} Given the parameter estimates $(\mfhu{A}, \mfhu{B}, \mfh A, \mfh B)$ and the empirical distribution $\hat{\mathbb{P}}$, the retailer then solves the following optimization problem:
\begin{equation}\label{model::ETO}\tag{$\code{ETO}$}
    \hat{\pi}^\star_{\code{ETO}} = \left[
    \begin{array}{c@{\;}l@{\;}l@{\;}}
        \max~& \mathbb{E}_{\hat{\mathbb P}} 
        \Big[
        \Pi(\mfhu A, \mfhu B, \mfh A, \mfh B, \bmu v, \bm v, \bm x, \bm y, \bmt \xi, \bmtu \epsilon, \bmt \epsilon) 
        \Big] \\
        {\rm s.t.} ~ &\bm x \in \mathcal{X}, \; \bm y \in \mathcal{Y}
    \end{array}
    \right].
\end{equation}
Provided $\mathcal{X}$ and $\mathcal{Y}$ admit mixed-integer representations, Problem~\ref{model::ETO} can be reformulated as a mixed-integer linear program (MILP), making the estimate-then-optimize approach implementable with standard MILP solvers. However, this approach can suffer from the \textit{optimizer's curse} \citep{smith2006optimizer}. Because model~\ref{model::ETO} selects the decision with the highest estimated profit under the empirical distribution, the optimization can favor decisions whose estimated performance is inflated by sampling noise. As a result, the empirical optimal profit $\hat{\pi}^\star_{\code{ETO}}$ could be optimistically biased relative to the true expected profit. This issue is particularly pronounced when the available sample size $S$ is small, which is common in operational settings. 

\section{Estimate-then-Robustify Model}\label{sec:ETR}

Motivated by the limitations of the estimate-then-optimize approach discussed above, and building on the insights of \citet{long2023robust} and \citet{kannan2023residuals}, we propose an estimate-then-robustify framework. This framework not only addresses the challenges inherent in the~\ref{model::CSO} model but also mitigates the drawbacks arising from optimizing over empirical distributions in the~\ref{model::ETO} model. The estimate-then-robustify model is formally defined as follows:
\begin{equation}\label{model::ETR}\tag{$\code{ETR}$}
    \hat{\kappa}^{\star}_{\tau} = \left[
    \begin{array}{c@{\;}l@{\;}l@{\;}}
        \min~& \kappa \\
        {\rm s.t.}~& \ep{\Pi(\mfhu A, \mfhu B, \mfh A, \mfh B, \bmu v, \bm v, \bm x, \bm y, \bmt \xi, \bmtu \epsilon, \bmt \epsilon)} \geq \tau - \kappa \Delta (\bbp, \hat{\bbp})  & ~ \forall \bbp \in \mathcal{M}(\cz) \\
        & \bm x \in \cx, \; \bm y \in \cy, \; \kappa \in \bbr_+
    \end{array}
    \right].
\end{equation}
Here, the parameter $\tau$ denotes the retailer's given profit target; the distance metric $\Delta (\bbp, \hat{\bbp})$ is a type-$1$ Wasserstein distance that quantifies the distance between $\mathbb{P}$ and  $\hat{\mathbb{P}}$, and is defined as follows: 
\[
\Delta(\mathbb{P}, \hat{\mathbb{P}}) := \inf_{\mathbb{Q} \in \mathcal{M}(\cz \times \cz)} \; \Big\{\mathbb{E}_{\mathbb{Q}}\, [\|\bmt z - \bmt \nu \|]~\big|~ (\bmt z, \bmt \nu) \sim \mathbb{Q}, \; \bmt z \sim \mathbb{P}, \; \bmt \nu \sim \hat{\mathbb{P}}
\Big\}.
\]

Unless stated otherwise, we use the block-separable ground norm
\[
\|(\bm \xi,\bmu\epsilon,\bm\epsilon)\|_{\oplus}
:=
\|\bm \xi\|_2+\|\bmu\epsilon\|_2+\|\bm\epsilon\|_2
\]
in the Wasserstein distance $\Delta$. This choice leads to the blockwise lifted support set and the corresponding dual-norm constraints in the conic reformulations below.

The profit target $\tau$ specifies the minimum profit level that the retailer seeks to maintain relative to a reference benchmark, such as the empirical optimal profit $\hat{\pi}^\star_{\code{ETO}}$. Following \citet{long2023robust}, we parameterize the target as $\tau=\hat{\pi}^\star_{\code{ETO}} - (1-\vartheta)|\hat{\pi}^\star_{\code{ETO}}|$, where $\vartheta\in[0,1]$ controls the target level and, consequently, the conservatism of the resulting decision. When $\vartheta=1$, \ref{model::ETR} selects, among the optimal solutions of~\ref{model::ETO}, one with the smallest fragility. While existing studies typically treat $\vartheta$ as an exogenous tuning parameter, we follow the idea of \citet{kannan2023residuals} and develop a cross-validation procedure that selects $\vartheta$ to maximize out-of-sample profit  based on realized covariates. Appendix~\ref{app:protocol} provides the details in Algorithm~\ref{alg::tau}.

Model~\ref{model::ETR} differs from model~\ref{model::ETO}, which maximizes expected profit under the empirical distribution. Instead, model~\ref{model::ETR} selects inventory and advertising decisions $(\bm{x},\bm{y})$ that achieve the profit target $\tau$ while limiting their sensitivity to distributional shifts away from $\hat{\mathbb{P}}$. Specifically, $\kappa$ measures the rate at which the profit guarantee deteriorates as a candidate distribution $\mathbb{P}$ moves away from $\hat{\mathbb{P}}$. The constraint requires the expected profit under any $\mathbb{P}\in\mathcal{M}(\mathcal{Z})$ to be at least $\tau-\kappa\Delta(\mathbb{P},\hat{\mathbb{P}})$. Thus, a smaller $\kappa$ provides a stronger robustness guarantee: the lower bound on expected profit decreases more slowly as the distribution deviates from the empirical distribution.

In summary, like model~\ref{model::ETO}, model~\ref{model::ETR} first estimates the net-demand and return parameters from historical data and then treats these estimates as fixed during optimization. Model~\ref{model::ETR} accounts for residual distributional ambiguity by requiring the expected profit under every $\mathbb{P}\in\mathcal{M}(\cz)$ to remain above $\tau-\kappa\Delta(\mathbb{P},\hat{\mathbb{P}})$. Rather than introducing a separate uncertainty set for the estimated parameters, we account for parameter estimation error indirectly through its effect on the distributional discrepancy $\Delta(\mathbb{P},\hat{\mathbb{P}})$, as formalized by the finite-sample performance guarantee. Explicitly incorporating parameter uncertainty into the optimization model would instead lead to an intractable bilinear formulation.

We first verify that model~\ref{model::ETR} is feasible with bounded infimal fragility when the empirical optimum is attained.

\begin{prop}\label{prop:bound}
Suppose model~\ref{model::ETR} attains a finite optimal value. For any $\tau \le \hat{\pi}^\star_{\code{ETO}}$, model~\ref{model::ETR} is feasible and its infimal value satisfies $\hat{\kappa}^{\star}_{\tau} \le C_{\code{ETO}} < \infty$, where $C_{\code{ETO}}$ is a Lipschitz constant of $\Pi(\mfhu A, \mfhu B, \mfh A, \mfh B, \bmu v, \bm v, \hat{\bm x}_{\code{ETO}}, \hat{\bm y}_{\code{ETO}}, \cdot, \cdot, \cdot)$ with respect to $\|\cdot\|_\oplus$ and $(\hat{\bm x}_{\code{ETO}},\hat{\bm y}_{\code{ETO}})$ is an optimal solution of~\ref{model::ETO}.
\end{prop}

Under the parameterization $\tau=\hat{\pi}^\star_{\code{ETO}} - (1-\vartheta)|\hat{\pi}^\star_{\code{ETO}}|$, the condition $\tau \le \hat{\pi}^\star_{\code{ETO}}$ holds for either sign of $\hat{\pi}^\star_{\code{ETO}}$ whenever $\vartheta \in [0,1]$, so \Cref{prop:bound} applies to every target the retailer selects.

The remainder of this section analyzes model~\ref{model::ETR} from two perspectives. Section~\ref{subsec:guarantee} studies its statistical performance: we establish a finite-sample guarantee on the out-of-sample expected profit of the~\ref{model::ETR} solution, and we show that the solution is asymptotically optimal for~\ref{model::CSO} when the parameter estimator is consistent and the target approaches the empirical optimal profit. Sections~\ref{subsec:exact} and~\ref{subsec:approximation} address computation: we reformulate~\ref{model::ETR} as a mixed-integer conic program whose size grows exponentially in the number of products, and we then develop an affine recourse approximation whose size is polynomial in the number of products and which is exact for a single product.

\subsection{Finite-Sample Guarantee and Asymptotic Optimality}\label{subsec:guarantee}
Solving model~\ref{model::ETR} using the historical dataset $\mathcal{D}$ yields an optimal solution $(\hat{\bm x},\hat{\bm y})$. We next study how this data-driven solution performs out of sample under the true but unknown parameters $(\mfu A,\mfu B,\mf A,\mf B)$ and distribution $\mathbb{P}^\star$. A standard approach in data-driven optimization is to establish finite-sample guarantees on the out-of-sample performance of the prescribed decisions (see, \textit{e.g.}, \citealt{mohajerin2018data, kannan2023residuals}). In our setting, we seek a probabilistic guarantee of the form
\begin{equation}\label{eq: finite sample performance}
    \mathbb{P}^S\Big[\mathbb{E}_{\mathbb{P}^\star}[\Pi(\mfu A, \mfu B, \mf A, \mf B, \bmu v, \bm v, \hat{\bm x}, \hat{\bm y}, \bmt \xi, \bmtu \epsilon, \bmt \epsilon)] \geq \tau - \hat{\kappa}^{\star}_{\tau}\cdot \Gamma \Big] \geq 1 - \delta,
\end{equation}
where $\delta\in(0,1)$ is the prescribed significance level, $\Gamma$ quantifies the discrepancy between the empirical model and the true data-generating process, and $\mathbb{P}^S$ denotes the joint distribution of the $S$ samples used to construct $(\hat{\bm x},\hat{\bm y})$.

The analysis in this subsection has two objectives:
\begin{enumerate}
    \item[(a)] establish, for a given sample size $S$, a $(1-\delta)$-confidence guarantee on the out-of-sample expected profit of $(\hat{\bm x},\hat{\bm y})$, as in~\eqref{eq: finite sample performance}; and
    \item[(b)] show that the out-of-sample expected profit of $(\hat{\bm x},\hat{\bm y})$ converges to the optimal value of~\ref{model::CSO} as $S \to \infty$.
\end{enumerate}

To establish a finite-sample performance guarantee for model~\ref{model::ETR}, we impose an assumption on the data-generating process and a light-tail condition on the distribution of ad clicks and error terms.

\begin{assm}[Data-Generating Process]\label{assm: DGP}

(i) The historical ad clicks and error terms $\bmh z_s:=(\bmh \xi_s, \bmu \epsilon_s^\dag, \bm \epsilon_s^\dag)$, $s \in [S]$, are independent and identically distributed samples of $\bmt z =(\bmt \xi, \bmtu \epsilon, \bmt \epsilon) \sim \bbp^\star$, and the historical net demand and returns are generated by the true parameters according to 
$$
\begin{array}{lll}
	\bmhu d_s  = \mfu A \, \bmhu v_s + \mfu B (\bmh y_s \odot \bmh \xi_s) + \bmdu \epsilon_s, \quad
	\bmh r_s  = \mf A \, \bmh v_s + \mf B (\bmh y_s \odot \bmh \xi_s) + \bmd \epsilon_s, \qquad s \in [S],
\end{array}
$$
where the historical covariates and advertising decisions $(\bmhu v_s, \bmh v_s, \bmh y_s)$, $s \in [S]$, are treated as given.
(ii) The sample second-moment matrices
$$
\bmu \Sigma := \frac{1}{S} \sum_{s \in [S]} (\bmhu v_s, \bmh y_s \odot \bmh \xi_s) (\bmhu v_s, \bmh y_s \odot \bmh \xi_s)^\top , \quad \bm \Sigma := \frac{1}{S} \sum_{s \in [S]} (\bmh v_s, \bmh y_s \odot \bmh \xi_s) (\bmh v_s, \bmh y_s \odot \bmh \xi_s)^\top 
$$
are positive definite.
(iii) The random vector $\bmt z \sim \bbp^\star$ faced by the retailer in the selling season is independent of the historical data $\{\bmh z_s\}_{s \in [S]}$.

\end{assm}

Part (i) of \Cref{assm: DGP} states that the historical data are generated by the model of \Cref{assm:demand and return model} under the same distribution $\bbp^\star$ that governs the selling season, with the historical covariates and advertising decisions held fixed. Here $\bmh \xi_s$ denotes the potential click volume of week $s$, and the realized ad exposure is $\bmh y_s \odot \bmh \xi_s$. The model therefore imposes no mechanical dependence between $\bmh \xi_s$ and $\bmh y_s$: when $\bmh y_{s,n} = 0$, the potential clicks of product $n$ are not revealed by the realized click records and do not enter the net-demand and return equations, and product $n$ enters the estimation with zero ad exposure, which is the exposure term $\bmh y_s \odot \bmh \xi_s$ specified by the model. 
Part (ii) requires the regressors to have full column rank, which presumes $S \ge \max\{J_1, J_2\} + N$ and sufficient variation in $(\bmh y_s, \bmh \xi_s)$ across weeks; both the on/off variation in the campaigns and the variation in click volumes across active weeks contribute to identifying the advertising-effect matrices. In our data, campaigns are active in 82.2\% of product-weeks, per-product activation shares range from 51\% to 100\%, and click volumes vary widely across active weeks (\Cref{table:weekly_stat_main}). Part (iii) separates the historical data from the selling season.

\begin{assm}[Light Tails]\label{assm: light tail}
The random vector $\bmt z
\sim \bbp^\star$ is light-tailed: there exists $q > 1$ such that $\varpi \coloneqq \epstar{\exp(\| \bmt z \|_{\oplus}^q)} < \infty$. 
\end{assm}

\Cref{assm: light tail} is the standard condition under which the empirical distribution of independent samples concentrates around $\bbp^\star$ in the type-$1$ Wasserstein distance \citep{fournier2015rate}. It is satisfied by random vectors with bounded support, which covers click volumes bounded by the number of impressions, and by sub-Gaussian random vectors for any $q \in (1,2)$.
		
\begin{thm}[Out-of-Sample Guarantee]\label{thm:out-of-sample}

Suppose Assumptions~\ref{assm:demand and return model}, \ref{assm: DGP}, and \ref{assm: light tail} hold, and that $\cz = \cz_\xi \times \bbr^{N} \times \bbr^{N}$ for some $\cz_\xi \subseteq \bbr^N_+$ with $\bbp^\star[\bmt\xi \in \cz_\xi] = 1$. 
Given $(\bmu v, \bm v)$ and any estimate $(\mfhu A, \mfhu B, \mfh A, \mfh B)$ of the model parameters, let $\hat{\bbp}$ be the empirical distribution of the residuals $\{(\bmh \xi_s, \bmhu \epsilon_s, \bmh \epsilon_s)\}_{s \in [S]}$ induced by this estimate as in Section~\ref{subsec:ETO}, and let $(\hat{\bm x}, \hat{\bm y},\hat{\kappa}^{\star}_{\tau})$ be an optimal solution of~\ref{model::ETR}. Then for any $\delta \in (0,1)$, with probability at least $1-\delta$,
\begin{align*}
	\epstar{\Pi(\mfu A, \mfu B, \mf A, \mf B, \bmu v, \bm v, \hat{\bm x},  \hat{\bm y}, \bmt \xi, \bmtu \epsilon, \bmt \epsilon)} \ge \tau - \hat{\kappa}^{\star}_{\tau} \cdot \Big(D + \big(\sqrt{\|\bmu \Sigma^{-1}\|} + \sqrt{\|\bm \Sigma^{-1}\|}\big)\Gamma_1\Big)\big(R  \Gamma_1 + D  \Gamma_2\big),
\end{align*}
where $D:=1 + \|\mfu B\|_{\rm F} + \|\mf B\|_{\rm F}$,
$
R:=  \max\limits_{s \in [S]}\left\{ \left\| \left[\begin{array}{c}
	\bmu v - \bmhu v_s \\
	( \hat{\bm y} - \bmh y_s) \odot \bmh \xi_s
\end{array} \right]\right\|_{\bmu \Sigma^{-1}}  + \left\|\left[\begin{array}{c}
	\bm v - \bmh v_s \\
	(\hat{\bm y} - \bmh y_s) \odot \bmh \xi_s
\end{array} \right]\right\|_{\bm \Sigma^{-1}} \right\},
$\\
$
\Gamma_1 := \max\big\{\|(\mfu A-\mfhu A,\mfu B-\mfhu B)\|_{\bmu \Sigma} ,\; \|(\mf A-\mfh A,\mf B-\mfh B)\|_{\bm \Sigma}\big\},
$
and 
$
\Gamma_2 := \left\{\begin{array}{ll}
	\left(\frac{\log (2\beta_1/\delta)}{\beta_2 S} \right)^{1/(3N)} & \text{ if } S \ge \frac{\log (2\beta_1/\delta)}{\beta_2}, \\[4pt]
	\left(\frac{\log (2\beta_1/\delta)}{\beta_2 S} \right)^{1/q} & \text{ if } S < \frac{\log (2\beta_1/\delta)}{\beta_2},
\end{array} \right.\\
$
$\beta_1,\beta_2$ are positive constants that depend only on $q$, $\varpi$, and $N$, and the matrix norm applied to $\bmu \Sigma^{-1}$ and $\bm \Sigma^{-1}$ is the spectral norm.
\end{thm}


\Cref{thm:out-of-sample} bounds the out-of-sample expected profit of the~\ref{model::ETR} solution from below by the target $\tau$ minus $\hat{\kappa}^\star_\tau \Gamma$, where $\Gamma:= (D + (\sqrt{\|\bmu \Sigma^{-1}\|}+ \sqrt{\|\bm \Sigma^{-1}\|})\Gamma_1)(R  \Gamma_1 + D  \Gamma_2)$. The fragility $\hat{\kappa}^\star_\tau$ is the quantity that model~\ref{model::ETR} minimizes, so the model minimizes the coefficient of the deviation term in its own guarantee. The term $\Gamma_1$ is the estimation error of the model parameters, and $\Gamma_2$ bounds, with probability at least $1-\delta/2$,  the Wasserstein distance between the empirical and the true distributions of ad clicks and error terms. The constant $D$ measures how the advertising-effect matrices transmit perturbations in ad clicks to net demand and returns, and $R$ measures the distance, in the norms $\|\cdot\|_{\bmu \Sigma^{-1}}$ and $\|\cdot\|_{\bm \Sigma^{-1}}$, between the new covariates and advertising decision and the historical ones. The term $\Gamma_2$ vanishes as $S \to \infty$ for any fixed $\delta$, and $\Gamma_1$ vanishes when the estimator is consistent.

The bound holds for any estimate $(\mfhu A, \mfhu B, \mfh A, \mfh B)$, because $\Gamma_1$ is defined as the realized estimation error of the given estimate; the probability refers to the sampling variability of the historical ad clicks and error terms, and $\hat{\bbp}$ is built from the residuals of the given estimate, whether obtained by ordinary least squares (OLS), by instrumental variables (IV), or by another method. The proof uses only the feasibility of $(\hat{\bm x}, \hat{\bm y},\hat{\kappa}^{\star}_{\tau})$, so the same bound holds for any feasible solution of~\ref{model::ETR} with its own $\kappa$, which allows the affine recourse approximation in Section~\ref{subsec:approximation} to preserve the guarantee. Because $D$ and $\Gamma_1$ depend on the true parameters, $\Gamma$ is not computable from the data. \Cref{cor:hp} below replaces $\Gamma_1$ by an explicit bound in terms of the sub-Gaussian variance proxy for the OLS estimator. 

For the OLS estimator under strict exogeneity, $\Gamma_1$ admits an explicit high-probability bound.

\begin{coro}[OLS estimates]\label{cor:hp}
Suppose the assumptions of \Cref{thm:out-of-sample} hold. Let $\mfu W \in \mathbb{R}^{S \times (J_1+N)}$ and $\mf W \in \mathbb{R}^{S \times (J_2+N)}$ denote the regressor matrices whose $s$-th rows are $(\bmhu v_s, \bmh y_s \odot \bmh \xi_s)^\top$ and $(\bmh v_s, \bmh y_s \odot \bmh \xi_s)^\top$, and suppose that, conditional on $(\mfu W, \mf W)$, the error vectors $\bmdu \epsilon_s$ and $\bmd \epsilon_s$, $s \in [S]$, are independent across $s$, have mean zero, and are sub-Gaussian with variance proxy $\sigma^2_\prox$, \textit{i.e.}, $\mathbb{E}[\exp(\bm a^\top \bmdu \epsilon_s) \mid \mfu W] \le \exp(\sigma^2_\prox \|\bm a\|_2^2/2)$ for all $\bm a \in \mathbb{R}^N$, and likewise for $\bmd \epsilon_s$ given $\mf W$. Let $(\mfhu A, \mfhu B, \mfh A, \mfh B)$ be the OLS estimate and $\bar m := \max\{J_1 + N,\, J_2 + N\}$. Then, for any $\delta \in (0,1)$, with probability at least $1-\delta$,
\[
\Gamma_1 \le \bar \Gamma_1 := \sigma_\prox \sqrt{\frac{2 N \bar m + 3 \log(4/\delta)}{S}}
\]
and
\[
\epstar{\Pi(\mfu A, \mfu B, \mf A, \mf B, \bmu v, \bm v, \hat{\bm x},  \hat{\bm y}, \bmt \xi, \bmtu \epsilon, \bmt \epsilon)} \ge \tau - \hat{\kappa}^{\star}_{\tau} \cdot \Big(D + \big(\sqrt{\|\bmu \Sigma^{-1}\|} + \sqrt{\|\bm \Sigma^{-1}\|}\big)\, \bar \Gamma_1\Big)\big(R\, \bar \Gamma_1 + D\, \Gamma_2\big).
\]
\end{coro}

The conditional mean-zero condition is the strict exogeneity condition under which the OLS estimator is unbiased. Under this condition, the estimation error is of order $\sqrt{N\bar m/S}$ for fixed $\delta$, which decays faster than the Wasserstein term $\Gamma_2$. Hence, when $R$ and the inverse second-moment matrices remain bounded as $S$ grows (Assumption~\ref{assm:design} below) the deviation term in \Cref{cor:hp} is dominated by $D^2 \Gamma_2$ for large $S$. If ad clicks are correlated with the error terms, the condition fails, the OLS estimate converges to a pseudo-true parameter, and $\Gamma_1$ does not vanish; \Cref{thm:out-of-sample} still holds with the realized $\Gamma_1$, but the guarantee no longer tightens as $S$ grows.

We next show that, when the estimator is consistent, the solution of~\ref{model::ETR} is asymptotically optimal for~\ref{model::CSO}. Let
\[
\pi^\star_{\code{CSO}} := \sup_{\bm x \in \cx,\, \bm y \in \cy} \epstar{\Pi(\mfu A, \mfu B, \mf A, \mf B, \bmu v, \bm v, \bm x, \bm y, \bmt \xi, \bmtu \epsilon, \bmt \epsilon)}
\]
denote the optimal value of~\ref{model::CSO}, which is finite under \Cref{assm: light tail}. For a dataset of size $S$, we index by $S$ the quantities that depend on it: the estimate $(\mfhu A_S, \mfhu B_S, \mfh A_S, \mfh B_S)$, the empirical distribution $\hat{\bbp}_S$, the empirical optimal profit $\hat{\pi}^\star_{\code{ETO},S}$, the target $\tau_S$, the optimal solution $(\hat{\bm x}_S, \hat{\bm y}_S, \hat{\kappa}^\star_{\tau,S})$ of~\ref{model::ETR}, the sample second-moment matrices $\bmu \Sigma_S$ and $\bm \Sigma_S$, and the estimation error $\Gamma_{1,S}$. We also define
\[
\bar R_S := \max_{\bm y \in \cy}\, \max_{s \in [S]}\left\{ \left\| \left[\begin{array}{c} \bmu v - \bmhu v_s \\ ( \bm y - \bmh y_s) \odot \bmh \xi_s \end{array} \right]\right\|_{\bmu \Sigma_S^{-1}}  + \left\|\left[\begin{array}{c} \bm v - \bmh v_s \\ (\bm y - \bmh y_s) \odot \bmh \xi_s \end{array} \right]\right\|_{\bm \Sigma_S^{-1}} \right\},
\]
which is the constant $R$ of \Cref{thm:out-of-sample} maximized over all candidate advertising decisions. The asymptotic analysis requires the historical regressors to remain well conditioned as $S$ grows.

\begin{assm}[Design Regularity]\label{assm:design}
There exists a constant $\lambda > 0$ such that, with probability tending to one as $S \to \infty$, (i) the smallest eigenvalues of $\bmu \Sigma_S$ and $\bm \Sigma_S$ are at least $\lambda$; and (ii) $\bar R_S \le 1/\lambda$.
\end{assm}

Condition (i) holds by the law of large numbers when the sample second-moment matrices converge to positive-definite population second-moment matrices, which requires non-collinear covariates and each product to be advertised with positive frequency in the historical data. Condition (ii) states that the new covariates and every candidate advertising decision remain at uniformly bounded distances from the historical regressors, measured in the norms $\|\cdot\|_{\bmu \Sigma_S^{-1}}$ and $\|\cdot\|_{\bm \Sigma_S^{-1}}$; given (i), it holds whenever covariates and clicks have uniformly bounded support, as click volumes do in practice.

\begin{prop}[Asymptotic Optimality of ETR]\label{prop:asymptotic}
Suppose Assumptions~\ref{assm:demand and return model}--\ref{assm:design} hold, and let $(\hat{\bm x}_S, \hat{\bm y}_S)$ be an optimal solution of~\ref{model::ETR} with a target $\tau_S \le \hat{\pi}^\star_{\code{ETO},S}$. If $\Gamma_{1,S} \to 0$ in probability, then $\hat{\pi}^\star_{\code{ETO},S} \to \pi^\star_{\code{CSO}}$ in probability as $S \to \infty$. If, in addition, $\hat{\pi}^\star_{\code{ETO},S} - \tau_S \to 0$ in probability, then
\[
\epstar{\Pi(\mfu A, \mfu B, \mf A, \mf B, \bmu v, \bm v, \hat{\bm x}_S,  \hat{\bm y}_S, \bmt \xi, \bmtu \epsilon, \bmt \epsilon)} \to \pi^\star_{\code{CSO}} \quad \text{in probability as } S \to \infty.
\]
\end{prop}

The proof bounds the out-of-sample suboptimality of the~\ref{model::ETR} solution by three terms and shows that each vanishes. The first term is the slack $\hat{\pi}^\star_{\code{ETO},S} - \tau_S$ between the empirical optimal profit and the target, which the retailer controls through $\vartheta$. The second term bounds the distance between the empirical optimal profit and $\pi^\star_{\code{CSO}}$; any method that optimizes an estimated model incurs this term, and it vanishes at the rate of $\Gamma_{1,S}$ and $\Gamma_2$. The third term is the deviation term of \Cref{thm:out-of-sample} with $\hat{\kappa}^\star_{\tau,S}$ replaced by a uniform bound obtained from \Cref{prop:bound}; this is the term that model~\ref{model::ETR} minimizes.

The proposition requires two conditions beyond \Cref{assm:design}. The condition $\Gamma_{1,S} \to 0$ is consistency of the estimator in the norms $\|\cdot\|_{\bmu \Sigma_S}$ and $\|\cdot\|_{\bm \Sigma_S}$; it holds for OLS under strict exogeneity (\Cref{cor:hp}) and requires suitable instrument and regularization conditions for the IV estimator (\Cref{rem:iv} below). Estimator consistency is thus needed for asymptotic optimality but not for the finite-sample guarantee. The condition $\hat{\pi}^\star_{\code{ETO},S} - \tau_S \to 0$ concerns the target. Under the parameterization $\tau_S = \hat{\pi}^\star_{\code{ETO},S} - (1 - \vartheta_S) |\hat{\pi}^\star_{\code{ETO},S}|$, the slack equals $(1-\vartheta_S)|\hat{\pi}^\star_{\code{ETO},S}|$, and it vanishes whenever $\vartheta_S \to 1$, because $\hat{\pi}^\star_{\code{ETO},S}$ converges to the finite limit $\pi^\star_{\code{CSO}}$ by the first statement of the proposition. For a fixed $\vartheta < 1$, the bound in the proof shows that the positive part of the suboptimality minus $(1-\vartheta)|\pi^\star_{\code{CSO}}|$ converges to zero in probability. This is an upper bound on the loss permitted by the target slack, and the cross-validation procedure in Algorithm~\ref{alg::tau} selects $\vartheta$ to trade it off against the finite-sample protection. 

\begin{rem}[IV estimates]\label{rem:iv}
\Cref{thm:out-of-sample} also holds for the product-specific regularized IV estimator of Section~\ref{sec:endo_extension}, with $\hat{\bbp}$ formed from its residuals, because the finite-sample result allows any parameter estimator. Applying~\Cref{prop:asymptotic} additionally requires consistency of that estimator; for regularized IV, this depends on instrument relevance and exclusion and on the asymptotic behavior of the penalty. Weak instruments can increase parameter uncertainty \citep{stock2005testing}, which enters the finite-sample bound through $\Gamma_{1,S}$.
\end{rem}

\Cref{thm:out-of-sample} and \Cref{prop:asymptotic} extend existing guarantees in two directions. \citet{kannan2023residuals} establish finite-sample guarantees and asymptotic optimality for residuals-based DRO with exogenous, decision-independent, and observable covariates, and \citet{zhu2024residuals} extend these results to decision-dependent uncertainty. Our results accommodate the unobservable covariate $\bmt \xi$, whose entry into the net-demand and return models depends on the advertising decision, and they apply to the target-based robust satisficing formulation, in which the fragility $\hat{\kappa}^\star_\tau$ rather than a prescribed ambiguity radius governs the degree of protection. In the robust satisficing literature, \citet{long2023robust} omit covariates and \citet{sim2025analytics} treat covariates as observable and exogenous.

\subsection{Exact Reformulations}\label{subsec:exact}
Having established the performance guarantee of the \ref{model::ETR} model, it is important to note that the model is not readily solvable using off-the-shelf solvers. The computational challenge stems from its intrinsic complexity, which becomes evident after some reformulations presented below.

Observe that~\ref{model::ETR} under the Wasserstein distance $\Delta$ can
be equivalently written as (see for example, \citealt{Chen_Sim_Xiong_2020robust}; \citealt{long2023robust}): 
\begin{equation}\label{model::rs_lift}
    \hat{\kappa}^{\star}_{\tau} = \left[
    \begin{array}{c@{\;}l@{\;}l@{\;}l@{\;}}
        &\min ~& \kappa \\
        &{\rm s.t.}~& \tau - \mathbb{E}_{\mathbb{Q}}\Big[\Pi(\mfhu A, \mfhu B, \mfh A, \mfh B, \bmu v, \bm v, \bm x, \bm y, \bmt \xi, \bmtu \epsilon, \bmt \epsilon)\Big] \leq \kappa \mathbb{E}_{\mathbb{Q}}[\tilde{w}]  & ~ \forall \mathbb{Q} \in \bar{\cq} \\
        && \bm x \in \cx, \; \bm y \in \cy, \; \kappa \in \bbr_+.
    \end{array}
    \right]
\end{equation}
where we introduce a lifted ambiguity set $\bar{\cq}$ incorporating both the primary random variables $(\bmt{\xi}, \bmtu{\epsilon}, \bmt{\epsilon})$ and lifted random variables $\tilde{w}$ and $\tilde{s}$: 
\begin{equation*}
    \Bar{\cq} := \left\{ 
    \bbp \in \mathcal{M}(\mathbb{R}^{3N+1} \times [S]) \left| 
    \begin{array}{ll}
        (\bmt\xi, \bmtu\epsilon, \bmt \epsilon, \tilde{w}, \tilde{s}) \sim \bbp &  \\
        \bbp[\tilde{s} = s] = 1/S & \forall s \in [S]\\
        \bbp[(\bmt\xi, \bmtu\epsilon, \bmt \epsilon, \tilde{w}) \in \widebar{\mathcal{Z}}_s \mid \tilde{s} = s ]= 1 & \forall s \in [S] \\
    \end{array}
    \right.
    \right\}
\end{equation*}
in which the lifted support set is defined as 
\begin{align*}
    \widebar{\mathcal{Z}}_s := \Big\{(\bm\xi, \bmu\epsilon, \bm \epsilon, w) \in \cz \times \bbr  \mid &~ \|\bm \xi - \bmh{\xi}_{s}\| + \|\bmu \epsilon - \bmhu \epsilon_{s}\| + \|\bm \epsilon - \bmh \epsilon_{s}\| \le w \Big\}.
\end{align*}

\begin{prop}\label{prop::P_rsr_refo}
    The \ref{model::ETR} model is equivalent to the following robust optimization problem:
    \begin{equation}\label{model::robust_opt_m1}
        \hat{\kappa}^{\star}_{\tau} = \left[
        \begin{array}{c@{\;}l@{\;}l@{\;}l@{\;}}
            \min~& \kappa \\
            {\rm s.t.}~& \bm t^\top \bm 1 \ge \tau S \\
            & (\bm p - \bm \alpha)^\top \big(\mfhu A\, \bmu v + \mfhu B ( \bm y \odot \bm \xi) + \bmu{\epsilon}\big) - (\bm \alpha + \bm \gamma)^\top \big(\mfh A\, \bm v + \mfh B ( \bm y \odot \bm \xi) + \bm{\epsilon}\big)   \\
            & ~~~~- (\bm h + \bm \alpha)^\top (\bm x - \mfhu A\, \bmu v - \mfh A\, \bm v - (\mfhu B + \mfh B) ( \bm y \odot \bm \xi) - \bmu{\epsilon} - \bm{\epsilon})^+ \\
            & ~~~~+ (\bm \alpha - \bm c)^\top \bm x 
            - \bm k^\top(\bm y \odot \bm \xi) + \kappa w \ge t_s \qquad  \forall (\bm \xi, \bmu \epsilon, \bm \epsilon, w) \in \widebar{\mathcal{Z}}_s, s \in [S]\\
            & \bm x \in \cx, \; \bm y \in \cy, \; \kappa \in \bbr_+,\; \bm t \in \bbr^S
        \end{array}
        \right].
    \end{equation} 
\end{prop}
This formulation, however, remains NP-hard even for $S=1$ due to its multiproduct nature (see, \textit{e.g.}, \citealt{ardestani2016robust,Hanasusanto_Kuhn_2018conic}). To illustrate,  note that 
\begin{equation}
    \label{eq: expand exponential pieces}
    \begin{aligned}
        & (\bm h + \bm \alpha)^\top (\bm x - \mfhu A\, \bmu v - \mfh A\, \bm v - (\mfhu B + \mfh B) ( \bm y \odot \bm \xi) - \bmu{\epsilon} - \bm{\epsilon})^+\\
        & = \max_{\ci \subseteq [N]} \;\left\{ \sum\limits_{n \in \ci} \; (h_n + \alpha_n) (x_n - \bmhu a_n^\top \bmu v - \bmh a_n^\top \bm v - (\bmhu b_n + \bmh b_n)^\top (\bm y \odot \bm \xi) - \uline{\epsilon}_n - \epsilon_n)\right\}
    \end{aligned}
\end{equation}
is a piecewise affine function involving $2^N$ pieces. Hence, there would be an exponential number of robust constraints in the formulation \eqref{model::robust_opt_m1}, making it difficult to solve. Nevertheless, the following theorem gives an equivalent conic reformulation, which we refer to as the \code{EXACT} model.

\begin{thm}\label{thm::P_refo}
    Suppose $\mathcal{Z} = \big\{(\bm\xi, \bmu\epsilon, \bm \epsilon) \in \bbr^{3N} : \bm \xi \ge \bm 0, \bmu \epsilon \in [\bmhu\epsilon^{b},\bmhu\epsilon^{u}], \bm \epsilon \in [\bmh\epsilon^{b},\bmh\epsilon^{u}] \big\}$\, where $\bmhu\epsilon^{b}$ and $\bmhu\epsilon^{u}$ (resp.\ $\bmh\epsilon^{b}$ and $\bmh\epsilon^{u}$) are componentwise lower and upper bounds on the net-demand (resp.\ return) errors. Provided that $\cx$ and $\cy$ admit mixed-integer conic representations, Model~\eqref{model::robust_opt_m1}, and hence the \ref{model::ETR} model, can be equivalently reformulated as a mixed-integer conic optimization problem with $\mathcal{O}(2^N S N)$ number of decision variables and constraints, in addition to those required in $\mathcal{X}$ and $\mathcal{Y}$, where the~\ref{model::exact} model is given in Appendix~\ref{app:P_refo}.
\end{thm}

In the implementation, the lower and upper bounds of the support of $\bmu \epsilon$ (resp. $\bm \epsilon$) can be set to the componentwise minima and maxima of the training residuals $\{\bmhu\epsilon_s\}_{s \in [S]}$ (resp.\ $\{\bmh\epsilon_s\}_{s \in [S]}$). Restricting the support to this box removes distributions from the ambiguity set of~\ref{model::ETR}, so the fragility computed by~\ref{model::exact},  is a lower bound on the exact fragility over the unrestricted support for which \Cref{thm:out-of-sample} is stated. Although~\ref{model::ETR} can be equivalently reformulated as a mixed-integer conic optimization problem, \Cref{thm::P_refo} demonstrates that both the number of decision variables and constraints of the resulting~\ref{model::exact} model grow exponentially with the number of products $N$, rendering the model computationally intractable even for moderate problem sizes. To overcome this challenge, we reinterpret \ref{model::ETR} as a two-stage robust optimization problem and develop a computationally tractable approximation based on affine recourse adaptation. This approximation significantly reduces complexity while preserving the essential structure and robustness of the~\ref{model::ETR} model.

\subsection{Tractable Approximation}\label{subsec:approximation}

Recall that Model~\eqref{model::robust_opt_m1} is, in general, NP-hard. Nevertheless, it admits a natural interpretation as a two-stage robust optimization problem with recourse. In this formulation, the variables $\bm{x}$ and $\bm{y}$ are the \textit{here-and-now} decisions, while the expression $ (\bm x - \mfhu A\, \bmu v - \mfh A\, \bm v - (\mfhu B + \mfh B) ( \bm y \odot \bm \xi) - \bmu{\epsilon} - \bm{\epsilon})^+$ represents the \textit{wait-and-see} recourse decision, which we denote by $\bm{\ell} \in \mathbb{R}^N$. 
In other words, the variable $\bm \ell$ acts as an adaptive decision that is an arbitrary functional of the (infinitely many) realizations of the random variables associated with net demand, ad clicks, and return, making the problem even more challenging to optimize. 

A common approach in the robust optimization literature is to approximate such arbitrary functionals by tractable classes of functions—typically (piecewise) static or affine mappings, via an \textit{affine recourse adaptation} (see, \textit{e.g.}, \citealt{bertsimas2019adaptive}, \citealt{Chen_Sim_Xiong_2020robust}, \citealt{han2023finite}). That is, the recourse decision $\bm{\ell}$ is assumed to be an affine mapping of the random variables $(\bmt{\xi}, \bmtu{\epsilon}, \bmt{\epsilon}, \tilde{w})$. Specifically, we define the set of affine mappings:
\begin{align*}
    \cl := \left\{\bm{\ell}  \in \mathcal{R}^{3N+1, N}  ~\left\vert~
    \begin{aligned}
        & \exists \bm{\ell}^0,\bm{\ell}^{4} \in \mathbb{R}^N,  \mf{L}^{1}, \mf{L}^{2}, \mf{L}^{3} \in \mathbb{R}^{N \times N} \\
        &\bm \ell(\bm \xi, \bmu \epsilon, \bm \epsilon, w)  = \bm{\ell}^0 + \mf{L}^{1} \bm\xi + \mf{L}^{2} \bmu\epsilon + \mf{L}^{3} \bm\epsilon + \bm{\ell}^4 w
    \end{aligned}
    \right.
    \right\}
\end{align*}
where $\mathcal{R}^{3N+1, N}$ is the set of all measurable mappings from $\mathbb{R}^{3N+1}$ to $\mathbb{R}^{N}$. Correspondingly, we approximate Model~\eqref{model::robust_opt_m1} by the following robust optimization problem:
\begin{equation}\label{model::bar_ETR}\tag{$\widebar{\code{ETR}}$}
    \scalebox{0.92}{$
        \bar{\kappa}^\star_\tau :=\left[
        \begin{array}{c@{\;}l@{\;}l@{\;}l@{\;}l@{\;}}
            &\min~& \kappa \\
            &{\rm s.t.}~& \bm t^\top \bm 1 \ge \tau S \\
            && (\bm p - \bm \alpha)^\top \big(\mfhu A\, \bmu v + \mfhu B ( \bm y \odot \bm \xi) + \bmu{\epsilon}\big)  + (\bm \alpha - \bm c)^\top \bm x 
            \\
            &&~~~~ -(\bm \alpha + \bm \gamma)^\top \big(\mfh A\, \bm v + \mfh B ( \bm y \odot \bm \xi) + \bm{\epsilon}\big) - \bm k^\top(\bm y \odot \bm \xi)\\
            &&~~~~ -(\bm h+\bm \alpha)^\top \bm{\ell}_s(\bm\xi, \bmu\epsilon, \bm\epsilon, w) + \kappa w \ge t_s  && \forall (\bm \xi, \bmu \epsilon, \bm \epsilon, w) \in \widebar{\mathcal{Z}}_s, s \in [S]\\
            && \bm{\ell}_s(\bm\xi, \bmu\epsilon, \bm\epsilon, w) \ge \bm x - \mfhu A\, \bmu v - \mfh A\, \bm v - (\mfhu B + \mfh B) ( \bm y \odot \bm \xi) - \bmu{\epsilon} - \bm{\epsilon} &&  \forall (\bm \xi, \bmu \epsilon, \bm \epsilon, w) \in \widebar{\mathcal{Z}}_s, s \in [S] \\
            && \bm{\ell}_s(\bm\xi, \bmu\epsilon, \bm\epsilon, w) \ge \bm 0 &&   \forall (\bm \xi, \bmu \epsilon, \bm \epsilon, w) \in \widebar{\mathcal{Z}}_s, s \in [S] \\
            && \bm{\ell}_s \in \cl && \forall s \in [S] \\
            && \bm x \in \cx, \; \bm y \in \cy, \; \kappa \in \bbr_+, \; \bm t \in \mathbb{R}^S
\end{array}
\right].$}
\end{equation}
		
It follows that, on the same support $\mathcal{Z}$,  $\bar{\kappa}^\star_\tau \geq \hat{\kappa}^{\star}_{\tau}$. In other words, restricting $\bm{\ell}(\cdot)$ to an affine decision class $\mathcal{L}$ yields an upper bound $\bar{\kappa}^\star_\tau \geq \hat{\kappa}^{\star}_{\tau}$ on the optimal fragility of the~\ref{model::ETR} problem. Importantly, the~\ref{model::bar_ETR} model can be further reformulated as a mixed-integer conic optimization problem that can be handled by off-the-shelf solvers such as \textsc{Gurobi}.  
		
\begin{thm}\label{Theo::LDR_refo}
	Suppose $\cz$ as in \Cref{thm::P_refo}, and that $\cx$ and $\cy$ admit mixed-integer conic representations. The~\ref{model::bar_ETR} model can be reformulated as a mixed-integer conic optimization problem with $\mathcal{O}(SN^2)$ variables and constraints, in addition to those required in $\mathcal{X}$ and $\mathcal{Y}$. 
\end{thm}
		
\Cref{Theo::LDR_refo} establishes that, through the affine recourse adaptation,  the~\ref{model::bar_ETR} model has significantly fewer variables and constraints compared to the~\ref{model::exact} model. We present the final reformulation~\eqref{model::RS_L_CP} in the Appendix~\ref{app:LDR_refo}. The following proposition further guarantees the optimality of the affine recourse formulation in the single-product case.
\begin{prop}\label{prop::optimality of single product}
	When $N=1$, we have $\hat{\kappa}^{\star}_{\tau} = \bar{\kappa}^\star_\tau$.
\end{prop}
		
We remark that the optimality of the affine recourse adaptation approximation has only been established in limited contexts, such as the moment-based ambiguity set \citep{Perakis2022robust,bertsimas2019adaptive}, or the scenario-wise ambiguity set \citep{Hao_He_Hu_Jiang_2020robust} under certain separability assumptions. In contrast, our study demonstrates that optimality can be achieved within an estimate-then-robustify framework that incorporates a scenario-wise lifted ambiguity set in the single-product case without separability assumptions. Although the theoretical results do not establish optimality for the multiproduct case ($N>1$), our numerical case study indicates that the approximate solutions of~\ref{model::bar_ETR} are nearly identical to the solution obtained from~\ref{model::exact}. 
		
Before moving to the numerical section, we remark that the affine recourse adaptation is conservative in a precise and reassuring sense. Any optimal solution $(\bar{\bm x}, \bar{\bm y}, \bar{\kappa}^\star_\tau)$ of~\ref{model::bar_ETR} is feasible for~\ref{model::ETR} on the same support $\mathcal{Z}$: since $\bm \ell_s \ge \max\{\bm x - \mfhu A \bmu v - \mfh A \bm v - (\mfhu B + \mfh B)(\bm y \odot \bm \xi) - \bmu \epsilon - \bm \epsilon, \bm 0\}$ pointwise on $\widebar{\cz}_s$, the affine recourse under-estimates profit, so satisfying the~\ref{model::bar_ETR} constraints implies satisfying the robust satisficing constraints of~\ref{model::ETR} with the same $\kappa$. Consequently, the finite-sample guarantee of \Cref{thm:out-of-sample} extends to the affine solution, with $\hat{\kappa}^\star_\tau$ replaced by the optimal fragility under affine recourse approximation. The tractable approximation therefore does not forfeit the out-of-sample performance guarantee; it only (weakly) inflates the fragility certificate.

\section{Model Validation and Numerical Experiments}\label{sec:case study}

Section~\ref{subsec::empirical evidences} estimates the linear net-demand and return models of \Cref{assm:demand and return model} on a proprietary dataset from our industrial partner and compares their predictive accuracy with that of alternative specifications. 
Section~\ref{subsec:experiments} then evaluates the computational efficiency of~\ref{model::bar_ETR} against~\ref{model::exact} and out-of-sample performance of~\ref{model::bar_ETR} against the~\ref{model::ETO} benchmark, as well as a spillover-ignorant variant. Section~\ref{sec:endo_extension} re-estimates the net-demand and return models using IV to address potential endogeneity in ad clicks, and repeats the experiment under the IV estimates.

\subsection{Data and Model Validation}\label{subsec::empirical evidences}

The dataset records sales, returns, and PPC advertising for 84 products over 76 weeks from March 2022 to August 2023, in 61,932 transaction records. Sales records contain the order time, order and product IDs, retail price, sales volume, ad clicks, and cost per click; return records contain the return time, order and product IDs, return volume, and refund amount. 

When inventory is unavailable at Amazon, customers can opt for shipment from the offshore warehouse, so orders are recorded even in stock-out weeks; the retailer also keeps enough stock that only 2.84\% of orders were fulfilled this way. We therefore treat observed sales as realized demand. The ten best-selling products account for 94.08\% of total sales and 90.19\% of total revenue, so we focus on ten selected products, labeled P1 through P10, and aggregate the data to the product-week level. After an isolation-forest outlier screen with a 1\% contamination rate and the removal of product-weeks with negative net demand, which arise when returns of earlier sales exceed the week's sales, the panel contains 695 product-week observations; \Cref{table:weekly_stat_main} in Appendix~\ref{appendix:addi_comp_resu} reports descriptive statistics.

We estimate the net-demand and return models of \Cref{assm:demand and return model} by ridge regression, with a distinct direct effect for each product and a distinct spillover effect for each ordered pair of products. For product $n$ in week $t$, the regression equations are
\begin{align}
\uline{d}_{nt} &= \uline{a}_{n0} + \bmu a^{\top} \bmu v_{nt} + \uline{b}_{nn}\,\xi_{nt} + \sum_{m \ne n} \uline{b}_{nm}\,\xi_{mt} + \uline{\epsilon}_{nt}, \label{eq:ridge_nd}\\
r_{nt} &= a_{n0} + \bm a^{\top} \bm v_{nt} + b_{nn}\,\xi_{nt} + \sum_{m \ne n} b_{nm}\,\xi_{mt} + \epsilon_{nt}, \label{eq:ridge_r}
\end{align}
where $\bmu v_{nt} = \bm v_{nt}$ collects the retail price $p_{nt}$, the week index $t$, and eleven indicators of the calendar month of the week's starting date (January is the reference month), and $\xi_{mt}$ is the recorded weekly ad clicks of product $m$, equal to zero in weeks in which $m$ is not advertised, so that it is the realized ad exposure $\bmh y_s \odot \bmh \xi_s$ of Section~\ref{subsec:ETO}. The intercepts $(\uline{a}_{n0}, a_{n0})$ are product-specific, $(\bmu a, \bm a)$ are shared across products within each outcome, and the advertising coefficients $\uline{b}_{nm}$ and $b_{nm}$ are the unrestricted entries of $\mfu B$ and $\mf B$ in \eqref{model:demand vector}--\eqref{model:return vector}, with rows indexing the receiving product and columns the advertised product, so each matrix contains $N$ direct and $N(N-1)$ spillover effects with no symmetry or sign restrictions. In the notation of Assumption~\ref{assm:demand and return model}, $\bmu v=\bm v$ collects a constant and the covariates of all products, and \eqref{eq:ridge_nd}--\eqref{eq:ridge_r} restrict the row $\bm a_n$ of $\mf A$ to the product-specific intercept $a_{n0}$ and the shared coefficients $\bm a$ on product $n$'s own covariates $\bm v_{nt}$, and likewise $\bmu a_n$ to $(\uline{a}_{n0},\bmu a)$. Each product-week is one regression observation, whereas a sample $s$ of Section~\ref{subsec:ETO} collects the joint observations across products in one week; here and in Appendix~\ref{app:prediction_models} weeks are indexed by $t$ and the hats that mark observed data are omitted. Appendix~\ref{app:prediction_models} gives the fitting procedure and reports the coefficient matrices in \Cref{table:product_specific_coefficients}.

We compare the linear specification with four alternatives: a log-linear specification, a return-rate specification, and two mixed specifications that pair a linear equation for one outcome with a log-linear equation for the other. All include product-specific direct effects and pair-specific spillover effects and are fitted by ridge; the return-rate specification fits a sales equation and splits its nonnegative sales prediction into net demand and returns using product-by-calendar-month return rates. Coefficients, penalties, and return rates are estimated from training observations only, and the five specifications are evaluated on the same 200 random 80/20 splits of the product-week observations. See Appendix~\ref{app:prediction_models} for more details. \Cref{table:product_specific_prediction} reports the prediction performance of all the specifications. 

\begin{table}[htbp]
    \centering
    \vspace{-0.2cm}
    \caption{Prediction performance of five different specifications}
    \label{table:product_specific_prediction}
    \begin{threeparttable}
        \renewcommand{\arraystretch}{1}
        \setlength{\tabcolsep}{3pt}
        \fontsize{9}{10}\selectfont
        \begin{tabular}{lccccc}
            \toprule
            Specification & In-sample MSE & Out-of-sample MSE & Net-demand MSE & Return MSE & 90\% range \\
            \midrule
            Linear & 39.83 & \textbf{54.67} & \textbf{99.35} & \textbf{9.99} & [32.75, 82.85] \\
            Log-linear & 68.52 & 110.23 & 207.34 & 13.12 & [38.05, 336.64] \\
            Return rate & \textbf{39.12} & 56.73 & 99.95 & 13.50 & [33.70, 84.06] \\
            Linear net/log return & 40.59 & 56.24 & \textbf{99.35} & 13.12 & [34.87, 83.39] \\
            Log net/linear return & 67.75 & 108.67 & 207.34 & \textbf{9.99} & [38.30, 324.62] \\
            \bottomrule
        \end{tabular}
        \begin{tablenotes}
            \item Notes. 
            In-sample mean squared error (MSE) is the full-panel fitted error; test MSE averages 200 splits, each with 556 training and 139 test observations. Joint MSE averages net-demand and return MSEs in original units. Outcome-specific columns report test errors; the range gives the 5th and 95th percentiles of split-level joint test MSE, not an interval for its mean.
        \end{tablenotes}
    \end{threeparttable}
    \vspace{-0.4cm}
\end{table}

\Cref{table:product_specific_prediction} shows that the linear model has the lowest joint test MSE (54.67) among the five specifications, compared with 56.73 for the return-rate model and 110.23 for the log-linear model. It also has the lowest return test MSE (9.99, tied with log net/linear return) and net-demand test MSE (99.35, tied with linear net/log return). We therefore use the linear equations in the numerical experiments, preserving the affine outcome structure of Section~\ref{sec:ETR}. Note that ridge regularizes the estimated advertising effects but does not resolve possible endogeneity of clicks; Section~\ref{sec:endo_extension} examines that issue under a separate IV specification.

\subsection{Numerical Experiments}\label{subsec:experiments}

We now evaluate the proposed approach on the data of Section~\ref{subsec::empirical evidences}, addressing three questions: (i) how the computation time and solution quality of the affine approximation~\ref{model::bar_ETR} compare with those of the exact reformulation~\ref{model::exact} as the number of products $N$ and the sample size $S$ grow; (ii) how much out-of-sample profit robustification adds over the~\ref{model::ETO} benchmark, and how this value depends on $N$, $S$, and the profit target; and (iii) how the gains differ with and without modeling cross-product spillovers. 
The experiments are implemented in Python and solved with Gurobi. Computation times are measured on a MacBook Pro with an Apple M3 Pro processor and 36~GB of RAM.

An instance consists of a number of products $N\in\{2,4,6,8,10\}$, a training sample of $S\in\{10,15,20,30\}$ weeks, and the associated test set. Following the fixed product order P1--P10 of Section~\ref{subsec::empirical evidences}, the instances with $N$ products use the first $N$ products. For each $(N,S)$ pair, we evaluate 20 training/test splits of the 76 weeks. Each split provides a 30-week training pool, whose first $S$ weeks form the training sample, and a disjoint 46-week test set, giving 400 instances in total. On each training sample, we estimate the net-demand and return models of Section~\ref{subsec::empirical evidences}, form the empirical distribution $\hat{\bbp}$ from their residuals, and solve~\ref{model::ETO} and~\ref{model::bar_ETR}, the latter with $\tau=\hat{\pi}^\star_{\code{ETO}} - (1-\vartheta)|\hat{\pi}^\star_{\code{ETO}}|$ and $\vartheta$ selected by the cross-validation (CV) procedure of Algorithm~\ref{alg::tau} on the training weeks only. The resulting~\ref{model::ETO} and~\ref{model::bar_ETR} policies are evaluated on the test weeks by residual replay under a full-panel calibration of the two models, which retains each week's covariates and baseline variation and replaces its ad exposure with that implied by the policy.  
We observe that, in the retailer's dataset, weekly profit is negative in 10 of the 76 weeks (13.2\%). Negative profits also arise under~\ref{model::ETO} in some instances, and a percentage profit gain computed for a single instance is then misleading: it has the wrong sign when the benchmark profit is negative and is arbitrarily large when the benchmark profit is close to zero. We therefore define the percentage profit gain of policy $A$ over its benchmark $B$ across matched instances $\mathcal G$ as
\begin{equation}
100\,\cdot \frac{\sum_{i\in\mathcal G}(\pi^A_i-\pi^B_i)}{\sum_{i\in\mathcal G}\pi^B_i}\%,
\label{eq:aggregate_profit_gain}
\end{equation}
where $\pi^A_i$ and $\pi^B_i$ are the average weekly out-of-sample profits in instance $i$. 
Appendix~\ref{app:protocol} gives the parameters and the remaining details of the design.

\textbf{Computational efficiency and performance efficacy.} 
We solve both~\ref{model::exact} and~\ref{model::bar_ETR} on 400 training instances, 20 per $(N,S)$ cell, at a common fixed target $\vartheta=0.98$. \Cref{table:time_ave} reports the mean total computation time and empirical profit gap, defined in Appendix~\ref{app:protocol}, of the two models. 

\begin{table}[!htbp]
\vspace{-0.2cm}
\centering
\caption{{Mean total computation time (seconds) and empirical profit gap (\%) of~\ref{model::exact} and~\ref{model::bar_ETR}}}
\label{table:time_ave}
\begin{threeparttable}
\renewcommand{\baselinestretch}{1}
\fontsize{8}{9.6}\selectfont
\setlength{\tabcolsep}{2pt}
\renewcommand{\arraystretch}{1.2}
\begin{tabular}{ccrrrrrrrrrrrrrrr}
\toprule
\multicolumn{2}{c}{$N$} & \multicolumn{3}{c}{2} & \multicolumn{3}{c}{4} & \multicolumn{3}{c}{6} & \multicolumn{3}{c}{8} & \multicolumn{3}{c}{10} \\
\cmidrule(lr){3-5}\cmidrule(lr){6-8}\cmidrule(lr){9-11}\cmidrule(lr){12-14}\cmidrule(lr){15-17}
\multicolumn{2}{c}{Model} & \ref{model::exact} & \ref{model::bar_ETR} & Gap (\%) & \ref{model::exact} & \ref{model::bar_ETR} & Gap (\%) & \ref{model::exact} & \ref{model::bar_ETR} & Gap (\%) & \ref{model::exact} & \ref{model::bar_ETR} & Gap (\%) & \ref{model::exact} & \ref{model::bar_ETR} & Gap (\%) \\
\midrule
\multirow{4}{*}{$S$} & 10 & 0.050 & 0.052 & 0.000000 & 0.448 & 0.297 & 0.000271 & 3.127 & 0.557 & 0.000000 & 20.310 & 1.535 & 0.000014 & 114.616 & 2.979 & 0.000018 \\
 & 15 & 0.072 & 0.081 & 0.000000 & 0.681 & 0.606 & 0.000000 & 4.765 & 1.129 & 0.000000 & 30.373 & 3.040 & 0.000000 & 182.601 & 6.422 & 0.000000 \\
 & 20 & 0.097 & 0.114 & 0.000000 & 0.894 & 0.549 & 0.000000 & 6.414 & 1.880 & 0.000000 & 40.678 & 4.997 & 0.000000 & 235.738 & 11.415 & 0.000000 \\
 & 30 & 0.144 & 0.207 & 0.000000 & 1.374 & 1.014 & 0.000000 & 10.109 & 4.328 & 0.000000 & 63.808 & 11.353 & 0.000000 & 384.992 & 25.671 & 0.000000 \\
\bottomrule
\end{tabular}
\end{threeparttable}
\end{table}

As shown in~\Cref{table:time_ave}, the computational time of~\ref{model::exact} increases sharply with the number of products $N$ and roughly linearly with the sample size $S$, whereas the time required for solving~\ref{model::bar_ETR} grows much more slowly, typically by only a few seconds. Despite this substantial improvement in efficiency, the two models exhibit nearly identical performance: the absolute average profit gap is consistently below 0.0003\%, with many cases showing a gap that rounds to zero at six decimal places. These findings demonstrate that although~\ref{model::bar_ETR} is an approximation, it delivers solutions that are almost the same as those of~\ref{model::exact}. Given its scalability and accuracy, we therefore focus on~\ref{model::bar_ETR} in the remainder of our analysis and exclude~\ref{model::exact} from subsequent experiments.

\begin{figure}[htbp]
    \vspace{-0.2cm}
    \centering
    \includegraphics[width=0.52\textwidth]{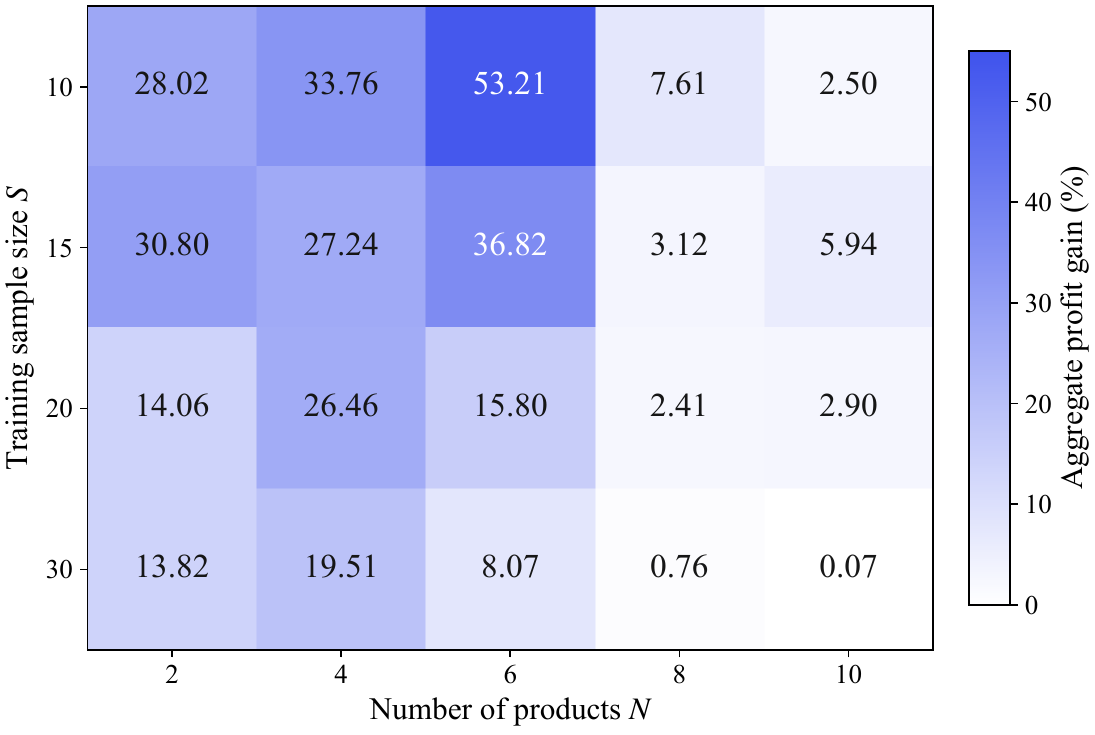}
    \caption{Percentage profit gains of~\ref{model::bar_ETR} over~\ref{model::ETO} with respect to training sample size $S$ and number of products $N$.} \label{fig:per_comp}
    \vspace{-0.2cm}
\end{figure}

\textbf{The value of robustification.} To evaluate the value of robustification, \textit{i.e.}, the benefit gained from addressing distributional ambiguity, we compare the out-of-sample average profits of~\ref{model::bar_ETR} with~\ref{model::ETO}. Figure~\ref{fig:per_comp} presents the percentage profit gain of~\ref{model::bar_ETR} relative to~\ref{model::ETO} across different training sample sizes $S$ and numbers of products $N$. 

Figure~\ref{fig:per_comp} shows positive percentage profit gains in all 20 cells, ranging from 0.07\% ($N=10$, $S=30$) to 53.21\% ($N=6$, $S=10$). Averaging all instances, \ref{model::bar_ETR} improves on~\ref{model::ETO} by 6.39\% and earns the higher profit in 74.5\% of instances. Pooling over $N$, the percentage profit gain decreases from 11.33\% at $S=10$ to 9.61\%, 5.59\%, and 2.88\% at $S=15,20,30$, respectively. Pooling over $S$, the gains are 19.82\%, 25.67\%, 21.36\%, 2.76\%, and 2.52\% for $N=2,4,6,8,10$. Robustification therefore adds the most value when training data are scarce and for portfolios with fewer products, although the gain is not monotone in portfolio size. This pattern is consistent with the optimizer's curse discussed in Section~\ref{subsec:ETO}, which can be more pronounced with limited historical data.

\begin{figure}[htbp]
    \vspace{-0.2cm}
    \centering
    \includegraphics[width=0.92\textwidth]{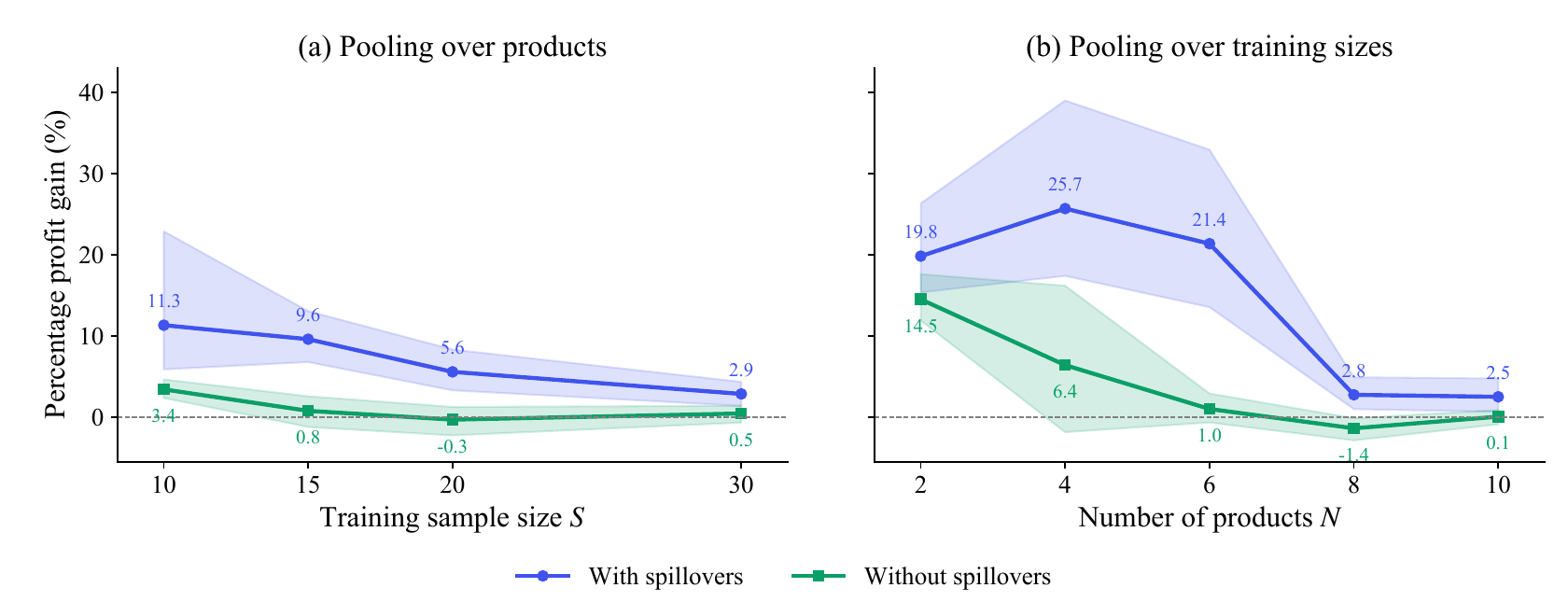}
    \caption{Percentage profit gain of~\ref{model::bar_ETR} over~\ref{model::ETO} with and without spillovers, (a) against the training sample size $S$, pooling over $N$, and (b) against the number of products $N$, pooling over $S$. 
    Shaded bands are 95\% bootstrap intervals. 
    }
    \label{fig:robustification_spillover_density}
\end{figure}

\textbf{On Spillover Effects.} In practice, retailers often overlook the spillover effects of online advertisements across different products. To examine the implications of this omission, we compare the out-of-sample performance of~\ref{model::bar_ETR} relative to~\ref{model::ETO} under two settings: one that accounts for spillover effects and one that ignores them. \Cref{fig:robustification_spillover_density} plots the resulting percentage profit gain against the training sample size $S$ (left panel) and the number of products $N$ (right panel). Across all instances, the percentage profit gain is 6.39\% with spillover effects and 1.02\% without.

\Cref{fig:robustification_spillover_density}(a) shows that the value of robustification declines with training sample size. Without spillover effects (green), the gain is 3.4\% at $S=10$ and its 95\% bootstrap interval includes zero for $S\ge15$. With spillover effects (blue), the gain decreases from 11.3\% at $S=10$ to 9.6\%, 5.6\%, and 2.9\% at $S=15,20,30$, respectively. The full specification exceeds the diagonal-only specification by 2.4 to 8.8 percentage points, and the paired 95\% bootstrap interval for this difference is positive at every $S$.

\Cref{fig:robustification_spillover_density}(b) shows that the full specification also produces a larger robustification gain at every $N$, with positive paired 95\% bootstrap intervals for the differences. At $N=2$, the gains are 19.8\% with spillovers and 14.5\% without. At $N=4$ and $6$, the gains are 25.7\% and 21.4\% with spillovers, compared with 6.4\% and 1.0\% without. For $N\in\{8,10\}$, ETO earns substantially higher profit and both specifications yield gains below 2.8\%. Without spillovers, the gain at $N=8$ is $-1.4$\%, with a 95\% bootstrap interval of $[-2.84,-0.11]$\%. The relative value of robustification is therefore greatest for the smaller portfolios in this experiment, particularly when cross-product effects are modeled.

These results show that, first, robustification pays most when data are scarce: a retailer with only a few weeks of relevant history stands to gain most from~\ref{model::bar_ETR}, whereas one with a long, stable history loses little by relying on~\ref{model::ETO}. Second, the relative robustification gain is larger under the specification that includes spillover effects. Robustification guards against estimation error within the assumed model, not against the omission of spillover effects, so robustifying a model that ignores them yields little and can even reduce profit. Conversely, the full specification must estimate cross-product effects from the same data, which plausibly makes its estimates noisier and leaves more room for robustification to protect profit. A retailer that advertises across related products should therefore pair a spillover-aware model with~\ref{model::bar_ETR} rather than plug its estimates directly into the optimization, especially when data are limited.

\subsection{Potential Endogeneity of Ad Clicks}\label{sec:endo_extension}

The estimates of Section~\ref{subsec::empirical evidences} treat ad clicks as exogenous conditional on the covariates. Clicks may, however, be correlated with unobserved demand shocks: the platform's auction algorithm can allocate impressions in response to real-time conversion rates, so weeks with high latent demand may also attract more clicks \citep{ghose2009empirical,du2017bidding}. A standard approach to addressing endogeneity is instrumental-variable (IV) estimation. Appendix~\ref{app:pooled_iv_evidence} reports first-stage evidence of instrument relevance and compares pooled OLS and IV estimates. Building on this evidence, we examine the operational value of robustification under IV estimation.

Retaining the product-specific structure of the main experiment, we use regularized two-stage least squares (2SLS) for net demand and returns, with product intercepts, shared controls, and all $N^2$ direct and product-pair advertising coefficients per outcome. The excluded instruments are the previous-calendar-week clicks of all $N$ products and market advertising expenditure excluding the receiving product, interacted with receiving-product indicators. This expanded set supplies $N(N+1)$ candidate instrument columns for the $N^2$ endogenous advertising regressors. Its use requires relevance and exclusion conditional on the controls, including the absence of direct lagged-advertising carryover and correlation between lagged clicks and current outcome shocks.

We use $S=30$, the largest training sample in the main design, to support the larger set of product-specific coefficients through the nested estimation steps. Five-fold target selection leaves 24 weeks for each fit, and penalty selection uses a further five-fold procedure within those weeks. Shorter samples leave less information for the first stages and product-pair effects.  The profit target is selected by the same training-only procedure as in the main experiment. For each $N\in\{2,4,6,8,10\}$, we evaluate 100 instances, each with 30 training weeks and 46 disjoint test weeks, yielding 500 matched~\ref{model::ETO}--\ref{model::bar_ETR} comparisons. The cost parameters are unchanged. Both policies are evaluated using a common full-panel product-specific regularized-IV calibration. We compute the percentage profit gain with the total~\ref{model::ETO} profit being positive for every reported group.

\begin{figure}[htbp]
    \vspace{-0.2cm}
    \centering
    \includegraphics[width=0.50\textwidth]{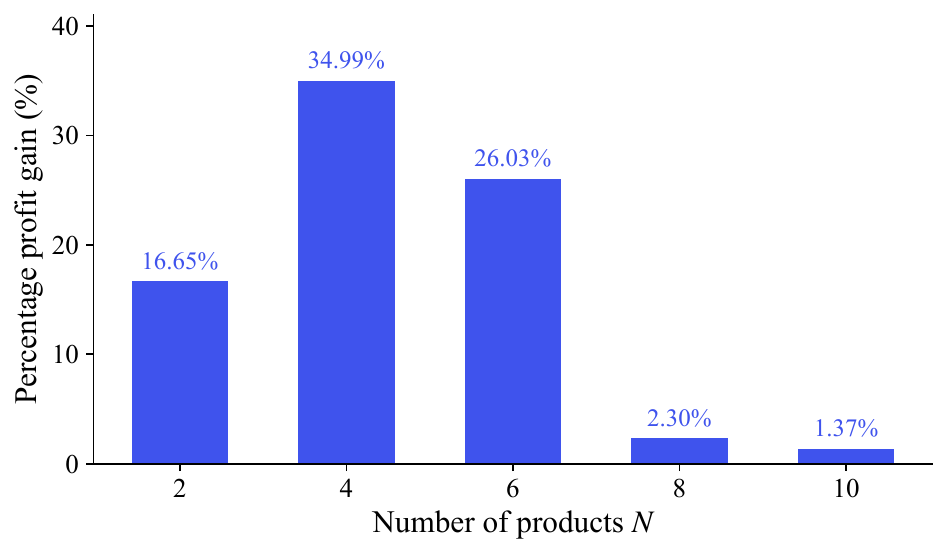}
    \caption{Percentage profit gain of~\ref{model::bar_ETR} over~\ref{model::ETO} under product-specific regularized IV estimation at $S=30$ for number of products $N \in \{2,4,6,8,10\}$.}
    \label{fig:iv_heatmap}
    \vspace{-0.3cm}
\end{figure}

Across all instances, the percentage profit gain is 5.73\%. The gains are 16.65\%, 34.99\%, 26.03\%, 2.30\%, and 1.37\% for $N=2,4,6,8,10$, respectively. Robustification thus adds value under the product-specific IV specification at every portfolio size. As in the main experiment, relative gains are larger for portfolios with fewer products and smaller for $N=8,10$. 

\section{Concluding Remarks} \label{section::conclusion}

We study a retailer's joint inventory and advertising decisions across multiple products when advertising affects both demand and returns and spills over across products. Using weekly data from a fashion retailer on Amazon, we estimate linear net-demand and return models with direct and spillover advertising effects; the direct effects are positive for both outcomes, while spillover effects vary in sign and magnitude. We formulate the joint decision as a contextual stochastic optimization problem in which the covariates are partially observable and enter the prediction through the advertising decision. To address distributional ambiguity, we propose an estimate-then-robustify approach---the~\ref{model::ETR} model---that minimizes the fragility of a specified profit target to deviations from the empirical distribution of clicks and residuals.

Under mild assumptions, we establish a finite-sample guarantee on the out-of-sample expected profit that explicitly accounts for parameter estimation error and distributional ambiguity. We further establish asymptotic optimality when the parameter estimates are consistent and the gap between the profit target and the empirical optimal profit vanishes. For specified uncertainty support, we derive an exact mixed-integer conic program and a conservative affine recourse approximation. The approximation has a polynomial formulation size, is exact for a single product, and substantially reduces computation time in our numerical experiments.

On these data, we find that the linear specification achieves the lowest joint prediction error among the five specifications evaluated. Building on the linear specification, in the case study, \ref{model::bar_ETR} improves on~\ref{model::ETO} by 6.39\% across all instances and earns higher profit in 74.5\% of them, with the largest gains when training data are scarce and for fewer products. The relative gain from robustification is also larger under the specification that includes spillovers than under the specification that omits them. These findings show that the value of robustification depends on both data availability and the predictive specification. 
With product-specific regularized IV estimates at $S=30$, robustification yields a 5.73\% percentage profit gain across 500 instances and a positive percentage profit gain at every portfolio size.

There are several potential directions for future research. First, the joint distribution of potential ad clicks and error terms is assumed not to depend on the advertising decision. Allowing it to depend on which products are advertised would capture interactions among campaigns at the click-generation stage in addition to the downstream spillover effects considered here. Second, the optimization treats estimated coefficients as fixed, with parameter estimation error entering the performance guarantee. A tractable formulation that jointly addresses parameter uncertainty and distributional ambiguity would extend the framework's protection. Finally, the numerical evidence comes from calibrated replay for one retailer; further studies could examine how the gains vary across product categories, advertising mechanisms, and operating environments.

\bibliographystyle{pomsref}
\let\oldbibliography\thebibliography
\renewcommand{\thebibliography}[1]{%
  \oldbibliography{#1}%
  \baselineskip14pt
  \setlength{\itemsep}{10pt}%
}
\bibliography{References}

\ECSwitch
\renewcommand{\theHsection}{EC.\arabic{section}}
\renewcommand{\theHequation}{EC.\arabic{equation}}
\renewcommand{\theHfigure}{EC.\arabic{figure}}
\renewcommand{\theHtable}{EC.\arabic{table}}

\ECHead{E-Companion of ``Data-Driven Robust Inventory Planning and Advertising under Product Return"}


\section{Proofs for Propositions and Theorems} \label{appendix:Proofs}

\subsection{Proof of Proposition \ref{prop:prel_result}\label{app:prel}}

For notation simplicity, we define $\bm \mu := \ep{\bmt \xi} > \bm 0$, $\bmt\zeta:= \bmt \xi - \bm \mu$, $\bmt \varepsilon \coloneqq (\mfu B+ \mf B) (\bm y \odot \bmt \zeta) + \bmtu\epsilon+\bmt\epsilon$,  and we denote by $\Phi(\cdot)$, $\phi(\cdot)$,  $\Phi^{-1}(\cdot)$ the CDF, PDF, and inverse CDF of the standard normal distribution. Then, the total demand for all products is expressed by: 
$\bmt d=\bmtu d+\bmt r=\mfu A \, \bmu v+\mf A \, \bm v+\mfu B (\bm y \odot \bmt \xi) + \bmtu\epsilon+\mf B (\bm y \odot \bmt \xi) + \bmt\epsilon=\mfu A \, \bmu v+\mf A \, \bm v+(\mfu B + \mf B) (\bm y \odot \bm \mu) + \bmt\varepsilon$.

Given the advertising decisions $\bm y$, we first consider the retailer's optimal inventory decisions for all products. Notably,  the retailer's inventory problem can be decomposed as $N$ separable newsvendor problems, and the first-order condition of the $n$-th problem is $\mathbb{P}[\tilde d_n \le x_n] = (\alpha_n - c_n)/(\alpha_n + h_n)$, which lies in $(0,1)$ because $\bm \alpha > \bm c$. Since $\tilde d_n$ is normal with mean $\bmu a_n^\top \bmu v + \bm a_n^\top \bm v + (\bmu b_n + \bm b_n)^\top(\bm y \odot \bm \mu)$ and variance $\Upsilon_n(\bm y)^2$, this yields the following optimal inventory decision: $\bm x^\star(\bm y)=\mfu A \, \bmu v+\mf A \, \bm v+(\mfu B+\mf B)(\bm y \odot \bm \mu) + {\bm \varkappa}$,
in which ${\bm \varkappa}$ represents the all products' safety stocks. 
Given the marginal variances $\text{Var}[\bmt\zeta]=\uwave{\bm\sigma}^2$, $\text{Var}[\bmtu \epsilon]=\bmu \sigma^2$ and $\text{Var}[\bmt \epsilon]=\bm\sigma^2$,  we have the following closed form expression:
$
\bm \varkappa = \bm\Upsilon(\bm y) \odot \Phi^{-1}\left(\frac{\bm \alpha-\bm c}{\bm \alpha+\bm h}\right),
$
where ${\bm \Upsilon(\bm y)}:=\sqrt{\mathrm{diag}((\mfu B +\mf B ) \mathrm{Diag}( \bm y \odot \uwave{\bm \sigma}^2)(\mfu B +\mf B )^\top) + \bmu\sigma^2+\bm\sigma^2 }$. 

We next substitute $\bm x^\star(\bm y)$ into the expected profit. Let $\tilde Z$ be a standard normal random variable and $z_n := \Phi^{-1}((\alpha_n - c_n)/(\alpha_n + h_n))$. Using $\mathbb{E}[(z_n - \tilde Z)^+] = z_n\Phi(z_n) + \phi(z_n)$ and $\mathbb{E}[(\tilde Z - z_n)^+] = \phi(z_n) - z_n(1-\Phi(z_n))$, the expected holding, backorder, and procurement costs at $\bm x^\star(\bm y)$ sum to $\bm c^\top \mathbb{E}[\bmt d] + \bm \Upsilon(\bm y)^\top((\bm \alpha + \bm h) \odot \phi(\Phi^{-1}((\bm \alpha - \bm c)/(\bm \alpha + \bm h))))$. The expected profit therefore equals
\begin{align*}
    &\bm p^\top(\mfu A \, \bmu v+\mf A \, \bm v+(\mfu B + \mf B) (\bm y \odot \bm \mu))-\bm \Upsilon(\bm y)^\top \left((\bm\alpha + \bm h) \odot \phi\Big(\Phi^{-1}\Big(\frac{\bm\alpha-\bm c}{\bm\alpha+\bm h}\Big)\Big) \right)\\
    &~~~~ - \bm c^\top (\mfu A \, \bmu v+\mf A \, \bm v+(\mfu B +\mf B)(\bm y \odot \bm \mu))-(\bm p + \bm\gamma)^\top(\mf A \, \bm v+ \mf B (\bm y \odot \bm \mu)) - \bm k^\top(\bm y \odot \bm \mu)\\
    =&(\bm p-\bm c)^\top(\mfu A \, \bmu v)-(\bm\gamma+\bm c)^\top(\mf A \, \bm v)+((\bm p-\bm c)^\top\mfu B-(\bm \gamma+\bm c)^\top \mf B -\bm k^\top) ( \bm y \odot \bm \mu)  -\bm \Upsilon(\bm y)^\top \left((\bm\alpha + \bm h) \odot \phi\left(\Phi^{-1}\left(\frac{\bm\alpha-\bm c}{\bm\alpha+\bm h}\right)\right) \right).
\end{align*}
Since the term $(\bm p-\bm c)^\top (\mfu A \, \bmu v) - (\bm\gamma+\bm c)^\top(\mf A \, \bm v)$ is a constant, the retailer's problem of optimizing the advertising decisions is equivalent to
\begin{align*}
    &\max_{\bm y\in \cy}\left\{ \left((\bm p- \bm c)^\top\mfu B-(\bm\gamma+\bm c)^\top \mf B -\bm k^\top\right) (\bm y \odot \bm \mu) -\bm \Upsilon(\bm y)^\top \left((\bm\alpha + \bm h) \odot \phi\left(\Phi^{-1}\left(\frac{\bm\alpha-\bm c}{\bm\alpha+\bm h}\right)\right) \right)\right\},
\end{align*}
It remains to verify supermodularity. The linear term in $\bm y$ is modular. For each product $n$, the nonlinear term can be written as
$-(\alpha_n+h_n)\phi(z_n) \sqrt{\upsilon_{n0}+\sum_{m\in[N]} \upsilon_{nm}y_m}$, $(\alpha_n+h_n)\phi(z_n)>0$, $\upsilon_{n0}=\uline{\sigma}_n^2+\sigma_n^2\ge0$, and $\upsilon_{nm}=(\uline{b}_{nm}+b_{nm})^2\uwave{\sigma}_m^2\ge0$, so that $\Upsilon_n(\bm y)^2=\upsilon_{n0}+\sum_{m\in[N]}\upsilon_{nm}y_m$.
For any $i\neq j$ and any $\bm y$ with $y_i=y_j=0$, letting $A=\upsilon_{n0}+\sum_{m\neq i,j}\upsilon_{nm}y_m$, we have
\[
\begin{aligned}
     &-(\alpha_n+h_n)\phi(z_n)\Big[\sqrt{A+\upsilon_{ni}+\upsilon_{nj}}-\sqrt{A+\upsilon_{ni}}-\sqrt{A+\upsilon_{nj}}+\sqrt{A}\Big] \\
    &=(\alpha_n+h_n)\phi(z_n)\left[\sqrt{A+\upsilon_{ni}}-\sqrt{A}-(\sqrt{A+\upsilon_{ni}+\upsilon_{nj}}-\sqrt{A+\upsilon_{nj}})\right]\ge0,
\end{aligned}
\]
because the function $t\mapsto \sqrt{t+\upsilon_{ni}}-\sqrt{t}$ is nonincreasing on $\mathbb{R}_+$. Thus each nonlinear term has nonnegative discrete cross-differences and is supermodular. Summing supermodular functions and a modular function preserves supermodularity, so the advertising objective is supermodular in $\bm y$. $\hfill\Halmos$

\subsection{Proof of the Two-Product Special Case in Proposition~\ref{prop:prel_result}}\label{app:prel-col}
When $N=2$, we write the objective function \eqref{prob:ad} as $g(y_1,y_2)$, where
$$
\begin{aligned} 
    g(y_1,y_2)=&((p_1- c_1)\uline b_{11}-(\gamma_1+c_1)b_{11})y_1\mu_1+((p_1- c_1)\uline b_{12}-(\gamma_1+c_1)b_{12})y_2\mu_2 \\
    &- (\alpha_1+h_1)\phi\left(\Phi^{-1}\left(\frac{\alpha_1-c_1}{\alpha_1 + h_1}\right)\right)\sqrt{y_1(\uline b_{11}+b_{11})^2\uwave\sigma^2_{1}+y_2(\uline b_{12}+b_{12})^2\uwave{\sigma}^2_{2} + \uline\sigma_{1}^2+\sigma_{1}^2 }\\
    &+ ((p_2- c_2)\uline b_{21}-(\gamma_2+c_2)b_{21})y_1\mu_1+((p_2- c_2)\uline b_{22}-(\gamma_2+c_2) b_{22})y_2\mu_2 \\
    &- (\alpha_2+h_2)\phi\left(\Phi^{-1}\left(\frac{\alpha_2-c_2}{\alpha_2 + h_2}\right)\right)\sqrt{y_1(\uline b_{21}+b_{21})^2\uwave\sigma^2_{1}+y_2(\uline b_{22}+b_{22})^2\uwave\sigma^2_{2} + \uline\sigma_{2}^2+\sigma_{2}^2 }-y_1 k_1\mu_1-y_2 k_2\mu_2.
\end{aligned}
$$
For ease of exposition, we define 
$$
\begin{aligned}
    \iota_1 &\coloneqq ((p_1- c_1)\uline b_{11}-(\gamma_1+c_1)b_{11}+(p_2- c_2)\uline b_{21}-(\gamma_2+c_2)b_{21}-k_1) \\
    \iota_2&\coloneqq((p_1- c_1)\uline b_{12}-(\gamma_1+c_1)b_{12}+(p_2- c_2)\uline b_{22}-(\gamma_2+c_2)b_{22}-k_2) \\
    \varphi_{y_1,y_2} &\coloneqq (\alpha_1+h_1)\phi\left(\Phi^{-1}\left(\frac{\alpha_1-c_1}{\alpha_1 + h_1}\right)\right)\sqrt{y_1(\uline b_{11}+b_{11})^2\uwave\sigma^2_{1}+y_2(\uline b_{12}+b_{12})^2\uwave\sigma^2_{2} + \uline\sigma_{1}^2+\sigma_{1}^2 }\\
    \hat \varphi_{y_1,y_2} & \coloneqq (\alpha_2+h_2)\phi\left(\Phi^{-1}\left(\frac{\alpha_2-c_2}{\alpha_2 + h_2}\right)\right)\sqrt{y_1(\uline b_{21}+b_{21})^2\uwave\sigma^2_{1}+y_2(\uline b_{22}+b_{22})^2\uwave\sigma^2_{2} + \uline\sigma_{2}^2+\sigma_{2}^2 }.
\end{aligned}
$$  
Then, the objective function can be rewritten as follows:
$$
g(y_1,y_2)=\iota_1 y_1\mu_1+\iota_2 y_2\mu_2-\varphi_{y_1,y_2}-\hat \varphi_{y_1,y_2}.
$$
Next, we discuss the optimal advertising decisions for the case of $N=2$ as follows:

(\textit{i}) If $g(1,1)\ge \max\{g(1,0),g(0,1),g(0,0)\}$, then $\bm y^\star=(1,1)$ and the condition is equivalent to $\iota_1\mu_1\geq -\varphi_{01}+\varphi_{11}-\hat{\varphi}_{01}+\hat{\varphi}_{11},\iota_2 \mu_2\ge   -\varphi_{10}+\varphi_{11}-\hat{\varphi}_{10}+\hat{\varphi}_{11},\iota_1\mu_1+\iota_2 \mu_2\geq  -\varphi_{00}+\varphi_{11}-\hat{\varphi}_{00}+\hat{\varphi}_{11}$;

(\textit{ii}) If $g(1,0)\ge \max\{g(1,1),g(0,1),g(0,0)\}$, then $\bm y^\star=(1,0)$ and the condition is equivalent to $ -\varphi_{10}+\varphi_{11}-\hat{\varphi}_{10}+\hat{\varphi}_{11}\geq \iota_2 \mu_2, \iota_1\mu_1\geq  -\varphi_{00}+\varphi_{10}-\hat{\varphi}_{00}+\hat{\varphi}_{10},\iota_1\mu_1-\iota_2 \mu_2\geq  -\varphi_{01}+\varphi_{10}-\hat{\varphi}_{01}+\hat{\varphi}_{10}$; 

(\textit{iii}) If $g(0,1)\ge \max\{g(1,1),g(1,0),g(0,0)\}$, then $\bm y^\star=(0,1)$ and the condition is equivalent to $-\varphi_{01}+\varphi_{11}-\hat{\varphi}_{01}+\hat{\varphi}_{11}\geq \iota_1\mu_1,\iota_2 \mu_2\geq -\varphi_{00}+\varphi_{01}-\hat{\varphi}_{00}+\hat{\varphi}_{01},\iota_1\mu_1-\iota_2 \mu_2\le  -\varphi_{01}+\varphi_{10}-\hat{\varphi}_{01}+\hat{\varphi}_{10}$;

(\textit{iv}) if $g(0,0)\ge \max\{g(1,1),g(1,0),g(0,1)\}$, then $\bm y^\star=(0,0)$ and the condition is equivalent to $-\varphi_{00}+\varphi_{10}-\hat{\varphi}_{00}+\hat{\varphi}_{10}\geq \iota_1\mu_1,-\varphi_{00}+\varphi_{01}-\hat{\varphi}_{00}+\hat{\varphi}_{01} \geq \iota_2\mu_2,\iota_1\mu_1+\iota_2 \mu_2 \leq  -\varphi_{00}+\varphi_{11}-\hat{\varphi}_{00}+\hat{\varphi}_{11}$.

By defining the following parameters: 
\begin{align*}
    \psi_{10}:=\varphi_{11}-\varphi_{10}-\hat{\varphi}_{10}+\hat{\varphi}_{11}, \, &\hat\psi_{10}:=\varphi_{10} -\varphi_{00}-\hat{\varphi}_{00}+\hat{\varphi}_{10}\\
    \psi_{01}:=\varphi_{11}-\varphi_{01}-\hat{\varphi}_{01}+\hat{\varphi}_{11}, \, &\hat\psi_{01}:=\varphi_{01} -\varphi_{00}-\hat{\varphi}_{00}+\hat{\varphi}_{01}
\end{align*}
we can prove the result by substituting these parameters into the conditions (\textit{i})$\sim$(\textit{iv}). 
$\hfill\Halmos$

\subsection{Proof of \texorpdfstring{\Cref{prop:bound}}{Proposition~\ref{prop:bound}}}\label{app:bound}

Let $(\bmh x_{\code{ETO}}, \bmh y_{\code{ETO}})$ be an optimal solution of~\ref{model::ETO}, and for $\bm z = (\bm \xi, \bmu \epsilon, \bm \epsilon) \in \mathbb{R}^{3N}$ write
$$
\Pi_{\code{ETO}}(\bm z) := \Pi(\mfhu A, \mfhu B, \mfh A, \mfh B, \bmu v, \bm v, \bmh x_{\code{ETO}}, \bmh y_{\code{ETO}}, \bm \xi, \bmu \epsilon, \bm \epsilon).
$$
By~\eqref{eq: profit function def}, $\Pi_{\code{ETO}}$ is the sum of an affine function of $\bm z$ and the function $-(\bm h + \bm \alpha)^\top \bm u(\bm z)^+$, where $\bm u(\bm z) := \bmh x_{\code{ETO}} - \mfhu A \bmu v - \mfh A \bm v - (\mfhu B + \mfh B)(\bmh y_{\code{ETO}} \odot \bm \xi) - \bmu \epsilon - \bm \epsilon$ is affine in $\bm z$ and $t \mapsto t^+$ is $1$-Lipschitz. Hence $\Pi_{\code{ETO}}$ is Lipschitz continuous with respect to $\|\cdot\|_\oplus$, and we let $C_{\code{ETO}}$ denote a Lipschitz constant, so that
$$
|\Pi_{\code{ETO}}(\bm z) - \Pi_{\code{ETO}}(\bm z')| \le C_{\code{ETO}}\, \|\bm z - \bm z'\|_\oplus \qquad \forall\, \bm z, \bm z' \in \mathbb{R}^{3N}.
$$
Since a Lipschitz constant can be enlarged, we take $C_{\code{ETO}} > 0$ without loss of generality.

Let $\cf_\oplus$ denote the set of functions $f: \mathbb{R}^{3N} \to \mathbb{R}$ with $|f(\bm z) - f(\bm z')| \le \|\bm z - \bm z'\|_\oplus$ for all $\bm z, \bm z' \in \mathbb{R}^{3N}$. For any two distributions $\bbq_1, \bbq_2 \in \mathcal{M}(\mathbb{R}^{3N})$ with finite first moments, the type-$1$ Wasserstein distance with ground norm $\|\cdot\|_\oplus$ admits the dual representation \citep{kantorovich1958space}
$$
\Delta(\bbq_1, \bbq_2) = \sup_{f \in \cf_\oplus} \big\{ \mathbb{E}_{\bbq_1}[f(\bmt z)] - \mathbb{E}_{\bbq_2}[f(\bmt z)] \big\}.
$$
Because $\Pi_{\code{ETO}} / C_{\code{ETO}} \in \cf_\oplus$, this representation yields
\begin{equation}\label{eq: KR Lipschitz}
\mathbb{E}_{\bbq_1}[\Pi_{\code{ETO}}(\bmt z)] - \mathbb{E}_{\bbq_2}[\Pi_{\code{ETO}}(\bmt z)] \le C_{\code{ETO}} \cdot \Delta(\bbq_1, \bbq_2).
\end{equation}

Now fix $\tau \le \hat{\pi}^\star_{\code{ETO}}$ and $\bbp \in \mathcal{M}(\cz)$, we have
\begin{equation}\label{eq: Lipschitz feasibility}
\begin{aligned}
\tau - \ep{\Pi_{\code{ETO}}(\bmt z)}
& \le \hat{\pi}^\star_{\code{ETO}} - \ep{\Pi_{\code{ETO}}(\bmt z)} \\
& = \ephat{\Pi_{\code{ETO}}(\bmt z)} - \ep{\Pi_{\code{ETO}}(\bmt z)} \\
& \le C_{\code{ETO}} \cdot \Delta(\hat{\bbp}, \bbp) = C_{\code{ETO}} \cdot \Delta(\bbp, \hat{\bbp}),
\end{aligned}
\end{equation}
where the first inequality is from $\tau \le \hat{\pi}^\star_{\code{ETO}}$, the equality is from the optimality of $(\bmh x_{\code{ETO}}, \bmh y_{\code{ETO}})$ for~\ref{model::ETO}, which gives $\hat{\pi}^\star_{\code{ETO}} = \ephat{\Pi_{\code{ETO}}(\bmt z)}$, the second inequality is from~\eqref{eq: KR Lipschitz} with $\bbq_1 = \hat{\bbp}$ and $\bbq_2 = \bbp$, both of which have finite first moments, and the last equality is from the symmetry of $\Delta$. Rearranging~\eqref{eq: Lipschitz feasibility}, we obtain
$$
\ep{\Pi(\mfhu A, \mfhu B, \mfh A, \mfh B, \bmu v, \bm v, \bmh x_{\code{ETO}}, \bmh y_{\code{ETO}}, \bmt \xi, \bmtu \epsilon, \bmt \epsilon)} \ge \tau - C_{\code{ETO}} \cdot \Delta(\bbp, \hat{\bbp}) \qquad \forall\, \bbp \in \mathcal{M}(\cz).
$$
Together with $\bmh x_{\code{ETO}} \in \cx$, $\bmh y_{\code{ETO}} \in \cy$, and $C_{\code{ETO}} \in \mathbb{R}_+$, this shows that $(\bmh x_{\code{ETO}}, \bmh y_{\code{ETO}}, C_{\code{ETO}})$ is feasible for~\ref{model::ETR}. Hence~\ref{model::ETR} is feasible, and its infimal value, the infimum of feasible $\kappa$ values, satisfies $\hat{\kappa}^\star_\tau \le C_{\code{ETO}} < \infty$.
$\hfill\Halmos$
		
\subsection{Proof of \texorpdfstring{\Cref{thm:out-of-sample}}{Theorem~\ref{thm:out-of-sample}}\label{app:out-of-sample}}
To prove the result, it suffices to show: suppose Assumptions~\ref{assm:demand and return model}, \ref{assm: DGP}, and \ref{assm: light tail} hold, then for optimal solution $(\hat{\bm x}, \hat{\bm y},\hat{\kappa}^{\star}_{\tau})$ of the \ref{model::ETR} model and any $\delta \in (0,1)$,  we have the following out-of-sample guarantee with probability at least $1-\delta$: 
\begin{align*}
    \epstar{\Pi(\mfu A, \mfu B, \mf A, \mf B, \bmu v, \bm v, \bmh x, \bmh y, \bmt \xi, \bmtu \epsilon, \bmt \epsilon)} \ge \tau - \hat{\kappa}^{\star}_{\tau} \cdot (D + (\sqrt{\|\bmu \Sigma^{-1}\|} + \sqrt{\|\bm \Sigma^{-1}\|})\Gamma_1)(R  \Gamma_1 + D  \Gamma_2)
\end{align*}

Given $(\bmu v, \bm v, \bmh y)$, we denote by  $\bbq^\star$ the distribution of $( \bmt \xi, \bmtu d, \bmt r) $, \textit{i.e.},  $$\left(\bmt \xi, \; \mfu A \, \bmu v + \mfu B (\bmh y \odot \bmt \xi) + \bmtu \epsilon, \; \mf A \, \bm v + \mf B (\bmh y \odot \bmt \xi) + \bmt \epsilon\right) \sim \bbq^\star, $$ by $\bbq^\dag$ the empirical distribution over the samples $$
\left( \bmh \xi_s, \; \mfu A \, \bmu v + \mfu B (\bmh y \odot \bmh \xi_s) + \bmdu \epsilon_s, \; \mf A \, \bm v + \mf B (\bmh y \odot \bmh \xi_s) + \bmd \epsilon_s\right), \quad s \in [S], 
$$
by $\hat \bbq$ the empirical distribution over  the samples $$
\left(\bmh \xi_s, \; \mfhu A \, \bmu v + \mfhu B (\bmh y \odot \bmh \xi_s) + \bmhu \epsilon_s, \; \mfh A \, \bm v + \mfh B (\bmh y \odot \bmh \xi_s) + \bmh \epsilon_s\right), \quad s \in [S], 
$$
and by $\bbp^\dag$ the empirical distribution of $(\bmt \xi, \bmtu \epsilon, \bmt \epsilon)$ over the samples $(\bmh \xi_s,   \bmdu \epsilon_s,  \bmd \epsilon_s)$, $s \in [S]$. 

We claim
\begin{equation}\label{eq: bound qstar and qhat}
    \Delta(\bbq^\star, \hat \bbq) \le R \Gamma_1 + D \Gamma_2
\end{equation}
with probability at least $1-\delta/2$.  
The proof is as follows. Note that by triangle inequality, we have $$ \Delta(\bbq^\star, \hat \bbq)  
\le    \Delta(\bbq^\star, \bbq^\dag) + \Delta(\bbq^\dag, \hat \bbq) . $$

We bound the term $\Delta(\bbq^\star, \bbq^\dag)$ as follows. By observing that $(\bmt \xi, \bmtu d, \bmt r)$ is an affine transformation of $(\bmt \xi, \bmtu \epsilon, \bmt \epsilon)$, we have 
\begin{equation}\label{eq: bound qstar and qdag}
    \Delta(\bbq^\star, \bbq^\dag) \le \left\|\left[ \begin{array}{ccc}
        \mf I & \mf 0 & \mf 0 \\
        \mfu B\mathrm{Diag}(\bmh y) & \mf I & \mf 0 \\
        \mf B\mathrm{Diag}(\bmh y) &  \mf 0 & \mf I
    \end{array}\right] \right\| \cdot   \Delta(\bbp^\star, \bbp^\dag)  \le D \cdot   \Delta(\bbp^\star, \bbp^\dag) \le D \cdot \Gamma_2
\end{equation}
where the matrix norm $\|\cdot\|$ denotes the operator norm induced by the ground norm $\|\cdot\|_{\oplus}$ on both domain and range (recall that the Wasserstein distance $\Delta$ is defined with $\|\cdot\|_{\oplus}$). The first inequality holds because the affine map transports any coupling of $\bbp^\star$ and $\bbp^\dag$ to a coupling of $\bbq^\star$ and $\bbq^\dag$ whose transport cost is at most the operator norm of the map times the original cost. The second inequality holds because, for any $\bm z:=(\bm \xi, \bmu \epsilon, \bm \epsilon) \in \bbr^{3N}$,
$$
\begin{aligned}
    \|\bm \xi\|_2 + \|\mfu B\mathrm{Diag}(\bmh y)\, \bm \xi + \bmu \epsilon\|_2 + \|\mf B\mathrm{Diag}(\bmh y)\, \bm \xi + \bm \epsilon\|_2
    &\le \|\bm z\|_{\oplus} + \big(\|\mfu B\mathrm{Diag}(\bmh y)\|_{\rm F} + \|\mf B\mathrm{Diag}(\bmh y)\|_{\rm F}\big)\|\bm \xi\|_2 \\
    &\le \big(1 + \|\mfu B\|_{\rm F} + \|\mf B\|_{\rm F}\big) \|\bm z\|_{\oplus} = D \|\bm z\|_{\oplus},
\end{aligned}
$$
using triangle inequality, $\|\mf M \bm u\|_2 \le \|\mf M\|_{\rm F} \|\bm u\|_2$ and $\|\mf B\mathrm{Diag}(\bmh y)\|_{\rm F} \le \|\mf B\|_{\rm F}$ for binary $\bmh y$. The last inequality follows from
$\Delta(\bbp^\star, \bbp^\dag) \le \Gamma_2 $ holds with probability at least $1-\delta/2$  under Assumption~\ref{assm: light tail} \citep{fournier2015rate}. 

We next bound the term $\Delta(\bbq^\dag, \hat \bbq)$. Note that 
\begin{equation}\label{eq: bound qdag and qhat}
    \begin{array}{r@{\;}l}
        & \Delta(\bbq^\dag, \hat \bbq)  \\
        \le &  \displaystyle \frac{1}{S} \sum_{s \in [S]} \left\|\left[\begin{array}{c}
            \bmh \xi_s - \bmh \xi_s  \\
            (\mfu A - \mfhu A) \bmu v  +  (\mfu B - \mfhu B) (\bmh y \odot \bmh \xi_s) + \bmu \epsilon_s^\dag - \bmhu \epsilon_s \\
            (\mf A - \mfh A) \bm v + (\mf B - \mfh B) (\bmh y \odot \bmh \xi_s) + \bm \epsilon_s^\dag - \bmh \epsilon_s
        \end{array} \right] \right\|_{\oplus} \\
        = &  \displaystyle \frac{1}{S} \sum_{s \in [S]} \left\|\left[\begin{array}{c}
            \bm 0 \\
            (\mfu A - \mfhu A) (\bmu v - \bmhu v_s)  +  (\mfu B - \mfhu B) ((\bmh y - \bmh y_s) \odot \bmh \xi_s) \\
            (\mf A - \mfh A) (\bm v - \bmh v_s) + (\mf B - \mfh B) ((\bmh y - \bmh y_s) \odot \bmh \xi_s)
        \end{array} \right] \right\|_{\oplus} \\
        \le &  \displaystyle \frac{1}{S} \sum_{s \in [S]} \left[\left\|
        (\mfu A - \mfhu A) (\bmu v - \bmhu v_s)  +  (\mfu B - \mfhu B) ((\bmh y - \bmh y_s) \odot \bmh \xi_s) \right\|_2 + \left\|
        (\mf A - \mfh A) (\bm v - \bmh v_s) + (\mf B - \mfh B) ((\bmh y - \bmh y_s) \odot \bmh \xi_s)
        \right\|_2\right] \\
        
        \le & \displaystyle \frac{1}{S} \sum_{s \in [S]}\left[\left\|
        (\mfu A - \mfhu A, \mfu B - \mfhu B)\right\|_{\bmu \Sigma}  \cdot \left\| \left[\begin{array}{c}
            \bmu v - \bmhu v_s \\
            (\bmh y - \bmh y_s) \odot \bmh \xi_s
        \end{array} \right]\right\|_{\bmu \Sigma^{-1}}  +  \left\|
        (\mf A - \mfh A, \mf B - \mfh B) \right\|_{\bm \Sigma} \cdot \left\|\left[\begin{array}{c}
            \bm v - \bmh v_s \\
            (\bmh y - \bmh y_s) \odot \bmh \xi_s
        \end{array} \right]\right\|_{\bm \Sigma^{-1}} \right]  \\
        
        \le & \displaystyle \frac{1}{S} \sum_{s \in [S]} \left[\Gamma_1  \cdot \left\| \left[\begin{array}{c}
            \bmu v - \bmhu v_s \\
            (\bmh y - \bmh y_s) \odot \bmh \xi_s
        \end{array} \right]\right\|_{\bmu \Sigma^{-1}}  +  \Gamma_1 \cdot \left\|\left[\begin{array}{c}
            \bm v - \bmh v_s \\
            (\bmh y - \bmh y_s) \odot \bmh \xi_s
        \end{array} \right]\right\|_{\bm \Sigma^{-1}}\right]  \\
        
        \le & \displaystyle \Gamma_1 \cdot R \\
    \end{array}  
\end{equation}
where the first inequality is from the definition of the Wasserstein distance, $\bbq^\dag$, and $\hat \bbq$, the equality is from the definition of $\bmhu \epsilon_s$ and $\bmh \epsilon_s$, the third inequality is from Assumption~\ref{assm: DGP} that $\bmu \Sigma$ and $\bm \Sigma$ are positive definite and the fact that $\|\mf Q \bm u\|_2 \le \sqrt{\mathrm{tr}(\mf Q \bm \Sigma \mf Q^\top)} \sqrt{\bm u^\top \bm \Sigma^{-1} \bm u}$ for any matrix $\mf Q$, vector $\bm u$ and positive definite $\bm \Sigma$, the last two inequalities are by the definition of $\Gamma_1$ and $R$. 

Combining the bounds~\eqref{eq: bound qstar and qdag} and \eqref{eq: bound qdag and qhat}, the second of which holds surely, we obtain $
\Delta(\bbq^\star, \hat \bbq)  
\le  R \Gamma_1 + D \Gamma_2 
$
with probability at least $1-\delta/2 \ge 1-\delta$.

Let  $\pi(\bm x, \bm y,  \bmt \xi, \bmtu d, \bmt r) := (\bm p - \bm \alpha)^\top \bmtu d - (\bm \alpha + \bm \gamma)^\top \bmt r-(\bm h + \bm \alpha)^\top(\bm x-\bmtu d - \bmt r)^+ + (\bm \alpha - \bm c)^\top \bm x - \bm k^\top (\bm y \odot \bmt\xi).$ Let $\bbp^\circ$ denote the distribution of $\big(\bmt \xi,\; \bmtu d - \mfhu A\,\bmu v - \mfhu B(\bmh y \odot \bmt \xi),\; \bmt r - \mfh A\,\bm v - \mfh B(\bmh y \odot \bmt \xi)\big)$ for $(\bmt \xi, \bmtu d, \bmt r) \sim \bbq^\star$, that is, the image of $\bbq^\star$ under the affine mapping $(\bmt \xi, \bmtu d, \bmt r) \mapsto (\bmt \xi, \bmtu d - \mfhu A\,\bmu v - \mfhu B(\bmh y \odot \bmt \xi), \bmt r - \mfh A\,\bm v - \mfh B(\bmh y \odot \bmt \xi))$; the same mapping sends $\hat\bbq$ to $\hat\bbp$. The $\bmt\xi$-marginal of $\bbp^\circ$ is that of $\bbp^\star$ and its error coordinates are unrestricted in $\cz$, so $\bbp^\circ \in \mathcal{M}(\cz)$; it has finite first moments under \Cref{assm: light tail}. Then, with probability at least $1-\delta$,
$$
\begin{array}{r@{\;}l}
    &   \mathbb{E}_{\mathbb P^\star}\Big[\Pi(\mfu A, \mfu B, \mf A, \mf B, \bmu v, \bm v, \bmh x, \bmh y, \bmt \xi, \bmtu \epsilon, \bmt \epsilon)\Big]
    = \mathbb{E}_{\mathbb Q^\star}\Big[\pi(\bmh x, \bmh y,  \bmt \xi, \bmtu d, \bmt r)\Big]
    = \mathbb{E}_{\bbp^\circ}\Big[\Pi(\mfhu A, \mfhu B, \mfh A, \mfh B, \bmu v, \bm v, \bmh x, \bmh y, \bmt \xi, \bmtu \epsilon, \bmt \epsilon)\Big] \\[4pt]
    \ge & \tau - \hat{\kappa}^{\star}_{\tau} \cdot \Delta(\bbp^\circ, \hat\bbp)
    \;\ge\; \tau - \hat{\kappa}^{\star}_{\tau} \cdot  (D + (\sqrt{\|\bmu \Sigma^{-1}\|} + \sqrt{\|\bm \Sigma^{-1}\|})\Gamma_1)(R \Gamma_1 + D \Gamma_2),
\end{array}
$$
where the first equality is the definition of $\bbq^\star$, the second equality is the identity $\pi(\bmh x, \bmh y, \bm \xi, \bmu d, \bm r) = \Pi(\mfhu A, \mfhu B, \mfh A, \mfh B, \bmu v, \bm v, \bmh x, \bmh y, \bm \xi, \bmu d - \mfhu A \, \bmu v - \mfhu B (\bmh y \odot \bm \xi) , \bm r- \mfh A \, \bm v - \mfh B (\bmh y \odot \bm \xi))$, the first inequality is the robust satisficing constraint of~\ref{model::ETR} at $\bbp = \bbp^\circ$, and the last inequality combines~\eqref{eq: bound qstar and qhat} with the bound
$$\begin{array}{r@{\;}l}
    \Delta(\bbp^\circ , \hat \bbp) \le & \left\|\left[ \begin{array}{ccc}
        \mf I & \mf 0 & \mf 0 \\
        -\mfhu B\mathrm{Diag}(\bmh y) & \mf I & \mf 0 \\
        -\mfh B\mathrm{Diag}(\bmh y) &  \mf 0 & \mf I
    \end{array}\right] \right\| \cdot    \Delta(\bbq^\star, \hat \bbq) \\
    \le &  \left( 1 + \|\mfhu B\mathrm{Diag}(\bmh y)\|_{\rm F} + \| \mfh B\mathrm{Diag}(\bmh y) \|_{\rm F} \right) \cdot \Delta(\bbq^\star, \hat \bbq) \\
    \le & \left( D + \|\mfu B - \mfhu B\|_{\rm F} + \|\mf B - \mfh B \|_{\rm F} \right) \cdot \Delta(\bbq^\star, \hat \bbq)\\
    \le & \left( D + (\sqrt{\|\bmu \Sigma^{-1}\|} + \sqrt{\|\bm \Sigma^{-1}\|}) \cdot \Gamma_1 \right) \cdot \Delta(\bbq^\star, \hat \bbq).
\end{array} $$
where the first inequality is from the definition of the Wasserstein distance and the matrix norm is the $\oplus\!\to\!\oplus$ operator norm, the second inequality follows from the same componentwise argument as in equation~\eqref{eq: bound qstar and qdag} (with $\mfhu B, \mfh B$ in place of $\mfu B, \mf B$), the second-to-last inequality is from the triangle inequalities $\|\mfhu B\|_{\rm F} \le \|\mfu B\|_{\rm F} + \|\mfu B - \mfhu B\|_{\rm F}$ and $\|\mfh B\|_{\rm F} \le \|\mf B\|_{\rm F} + \|\mf B - \mfh B\|_{\rm F}$, and the last inequality is because $\|\mf X\|_{\rm F} \le \|\mf X\|_{\bm \Sigma} \sqrt{\|\bm \Sigma^{-1}\|}$ for any matrix $\mf X$ and positive definite matrix $\bm \Sigma$, whence 
$$
\begin{array}{r@{\;}l}
    \|\mfu B - \mfhu B\|_{\rm F} + \|\mf B - \mfh B\|_{\rm F} & \le \|(\mfu A-\mfhu A,\mfu B-\mfhu B)\|_{\rm F} + \|(\mf A-\mfh A,\mf B-\mfh B)\|_{\rm F} \\
    & \le  \|(\mfu A-\mfhu A,\mfu B-\mfhu B)\|_{\bmu \Sigma} \cdot \sqrt{\|\bmu \Sigma^{-1}\|} +  \|(\mf A-\mfh A,\mf B-\mfh B)\|_{\bm \Sigma} \cdot \sqrt{\|\bm \Sigma^{-1}\|}\\
    & \le (\sqrt{\|\bmu \Sigma^{-1}\|} + \sqrt{\|\bm \Sigma^{-1}\|})\cdot \Gamma_1 
\end{array}
$$ by definition of $\Gamma_1$. 
\hfill\Halmos
			
\subsection{Proof of \texorpdfstring{\Cref{cor:hp}}{Corollary~\ref{cor:hp}}}\label{app:cor_hp}

We first bound the estimation error. Consider the net-demand regression; the return regression is treated identically. Let $\mfu Y := (\bmhu d_1, \ldots, \bmhu d_S)^\top \in \mathbb{R}^{S \times N}$, $\mfu E := (\bmdu \epsilon_1, \ldots, \bmdu \epsilon_S)^\top \in \mathbb{R}^{S \times N}$, and $\bmu \Theta := (\mfu A, \mfu B) \in \mathbb{R}^{N \times (J_1 + N)}$, so that \Cref{assm: DGP}(i) reads $\mfu Y = \mfu W \bmu \Theta^\top + \mfu E$ and $\bmu \Sigma = \mfu W^\top \mfu W/S$. The OLS estimate is $\bmhu \Theta^\top = (\mfu W^\top \mfu W)^{-1} \mfu W^\top \mfu Y$, whence $(\bmhu \Theta - \bmu \Theta)^\top = (\mfu W^\top \mfu W)^{-1} \mfu W^\top \mfu E$ and
$$
\|\bmhu \Theta - \bmu \Theta\|_{\bmu \Sigma}^2 = \mathrm{tr}\big( (\bmhu \Theta - \bmu \Theta) \bmu \Sigma (\bmhu \Theta - \bmu \Theta)^\top \big)
= \frac{1}{S}\, \mathrm{tr}\big( \mfu E^\top \mfu H \mfu E \big), \qquad \mfu H := \mfu W (\mfu W^\top \mfu W)^{-1} \mfu W^\top ,
$$
where $\mfu H$ is the orthogonal projection onto the column space of $\mfu W$, so that $\mfu H^2 = \mfu H$, $\|\mfu H\| = 1$, and $\mathrm{tr}(\mfu H) = J_1 + N =: \uline{m}$. Let $\bm e := \mathrm{vec}(\mfu E) \in \mathbb{R}^{SN}$ be the vector that stacks the columns of $\mfu E$, and let $\mf G := \mf I_N \otimes \mfu H \in \mathbb{R}^{SN \times SN}$ be the block-diagonal matrix with $N$ copies of $\mfu H$ on its diagonal, where $\otimes$ denotes the Kronecker product; then $\mathrm{tr}(\mfu E^\top \mfu H \mfu E) = \bm e^\top \mf G \bm e$, and $\mf G$ is an orthogonal projection with $\mathrm{tr}(\mf G) = \mathrm{tr}(\mf G^2) = N \uline{m}$ and $\|\mf G\| = 1$. Conditional on $\mfu W$, the assumption of the corollary and independence across $s$ imply that $\bm e$ is sub-Gaussian with variance proxy $\sigma^2_\prox$, since $\mathbb{E}[\exp(\sum_{s} \bm a_s^\top \bmdu \epsilon_s) \mid \mfu W] = \prod_s \mathbb{E}[\exp(\bm a_s^\top \bmdu \epsilon_s) \mid \mfu W] \le \exp(\sigma^2_\prox \sum_s \|\bm a_s\|_2^2 / 2)$. The tail inequality for quadratic forms of sub-Gaussian vectors \citep[Theorem~2.1]{hsu2012tail} then gives, for every $t > 0$,
$$
\mathbb{P}\Big[ \bm e^\top \mf G \bm e > \sigma^2_\prox \big( \mathrm{tr}(\mf G) + 2\sqrt{\mathrm{tr}(\mf G^2)\, t} + 2 \|\mf G\| t \big) \;\Big|\; \mfu W \Big] \le e^{-t} .
$$
Choosing $t = \log(4/\delta)$ and taking expectations over $\mfu W$, we obtain with probability at least $1 - \delta/4$
$$
\|\bmhu \Theta - \bmu \Theta\|_{\bmu \Sigma} \le \sigma_\prox \sqrt{\frac{N \uline{m} + 2\sqrt{N \uline{m} \log(4/\delta)} + 2\log(4/\delta)}{S}} \le \sigma_\prox \sqrt{\frac{2 N \uline{m} + 3\log(4/\delta)}{S}} ,
$$
where the second inequality uses $2\sqrt{ab} \le a + b$ with $a = N \uline{m}$ and $b = \log(4/\delta)$. The same bound holds for the return regression with $m := J_2 + N$ in place of $\uline{m}$. Since the right-hand side is nondecreasing in the number of regressors, both bounds are at most $\bar \Gamma_1$, and the union bound yields $\mathbb{P}^S[\mathcal{E}_1] \ge 1 - \delta/2$ for $\mathcal{E}_1 := \{\Gamma_1 \le \bar \Gamma_1\}$.

We now combine this bound with \Cref{thm:out-of-sample}. Let $\bbp^\dag$ denote the empirical distribution of $(\bmh \xi_s, \bmdu \epsilon_s, \bmd \epsilon_s)$, $s \in [S]$, as in the proof of \Cref{thm:out-of-sample}, and let $\mathcal{E}_2 := \{\Delta(\bbp^\star, \bbp^\dag) \le \Gamma_2\}$. The proof of \Cref{thm:out-of-sample} shows that $\mathbb{P}^S[\mathcal{E}_2] \ge 1 - \delta/2$ and that the inequality of \Cref{thm:out-of-sample} holds on $\mathcal{E}_2$. By the union bound, $\mathbb{P}^S[\mathcal{E}_1 \cap \mathcal{E}_2] \ge 1 - \delta$. On $\mathcal{E}_1 \cap \mathcal{E}_2$, both claims of the corollary hold: the first is the definition of $\mathcal{E}_1$, and the second follows from the inequality of \Cref{thm:out-of-sample} because its deviation term
$
\Big(D + \big(\sqrt{\|\bmu \Sigma^{-1}\|} + \sqrt{\|\bm \Sigma^{-1}\|}\big)\Gamma_1\Big)\big(R \Gamma_1 + D \Gamma_2\big)
$
is a product of two nonnegative factors, each nondecreasing in $\Gamma_1$, so replacing $\Gamma_1$ by $\bar \Gamma_1 \ge \Gamma_1$ can only increase it.
$\hfill\Halmos$

\subsection{Proof of \texorpdfstring{\Cref{prop:asymptotic}}{Proposition~\ref{prop:asymptotic}}}\label{app:asymptotic}

We use the notation of the proof of \Cref{thm:out-of-sample}, write $\bm \Theta := (\mfu A, \mfu B, \mf A, \mf B)$ for the true parameters and $\bmh \Theta := (\mfhu A, \mfhu B, \mfh A, \mfh B)$ for the given estimate, and let $\mathcal{E}_2 := \{\Delta(\bbp^\star, \bbp^\dag) \le \Gamma_2\}$ be the concentration event of the proof of \Cref{cor:hp}, so that $\mathbb{P}^S[\mathcal{E}_2] \ge 1 - \delta/2$ and the inequality of \Cref{thm:out-of-sample} holds on $\mathcal{E}_2$. Until the final part of the proof, we fix the sample size $S$ and the confidence level $\delta \in (0,1)$, suppress the index $S$, and argue on $\mathcal{E}_2$.

We first establish two Lipschitz constants. Define
$$
\bar L := \max\big\{ \|\bm k\|_2,\; \|\bm p - \bm \alpha\|_2 + \|\bm h + \bm \alpha\|_2,\; \|\bm \alpha + \bm \gamma\|_2 + \|\bm h + \bm \alpha\|_2 \big\}
$$
and, for matrices $\mfu B', \mf B' \in \mathbb{R}^{N \times N}$,
$$
\begin{aligned}
\bar C(\mfu B', \mf B') := \max\Big\{ & \|\bm p - \bm \alpha\|_2 + \|\bm h + \bm \alpha\|_2,\; \|\bm \alpha + \bm \gamma\|_2 + \|\bm h + \bm \alpha\|_2, \\
& \|\bm k\|_2 + \|\mfu B'\|_{\rm F}\big(\|\bm p - \bm \alpha\|_2 + \|\bm h + \bm \alpha\|_2\big) + \|\mf B'\|_{\rm F}\big(\|\bm \alpha + \bm \gamma\|_2 + \|\bm h + \bm \alpha\|_2\big) \Big\}.
\end{aligned}
$$
Fix $(\bm x, \bm y) \in \cx \times \cy$ and abbreviate $\Pi(\bm z) := \Pi(\bmh \Theta, \bmu v, \bm v, \bm x, \bm y, \bm z)$. For $\bm z = (\bm \xi, \bmu \epsilon, \bm \epsilon)$ and $\bm z' = (\bm \xi', \bmu \epsilon', \bm \epsilon')$ in $\mathbb{R}^{3N}$, let $\bm u(\bm z) := \bm x - \mfhu A \bmu v - \mfh A \bm v - (\mfhu B + \mfh B)(\bm y \odot \bm \xi) - \bmu \epsilon - \bm \epsilon$. From the definition~\eqref{eq: profit function def},
$$
\begin{aligned}
\Pi(\bm z) - \Pi(\bm z')
= \;& (\bm p - \bm \alpha)^\top \big( \mfhu B \mathrm{Diag}(\bm y)(\bm \xi - \bm \xi') + (\bmu \epsilon - \bmu \epsilon') \big)
- (\bm \alpha + \bm \gamma)^\top \big( \mfh B \mathrm{Diag}(\bm y)(\bm \xi - \bm \xi') + (\bm \epsilon - \bm \epsilon') \big) \\
& - (\bm h + \bm \alpha)^\top \big( \bm u(\bm z)^+ - \bm u(\bm z')^+ \big) - \bm k^\top \big(\bm y \odot (\bm \xi - \bm \xi')\big).
\end{aligned}
$$
We use the elementary bounds $|\bm a^\top \bm b| \le \|\bm a\|_2 \|\bm b\|_2$, $\|\mf M \mathrm{Diag}(\bm y) \bm w\|_2 \le \|\mf M\|_{\rm F} \|\bm w\|_2$ and $\|\bm y \odot \bm k\|_2 \le \|\bm k\|_2$ for binary $\bm y$, and $\|\bm u^+ - \bm u'^+\|_2 \le \|\bm u - \bm u'\|_2$, which holds because $t \mapsto t^+$ is $1$-Lipschitz componentwise. Since $\|\bm u(\bm z) - \bm u(\bm z')\|_2 \le (\|\mfhu B\|_{\rm F} + \|\mfh B\|_{\rm F})\|\bm \xi - \bm \xi'\|_2 + \|\bmu \epsilon - \bmu \epsilon'\|_2 + \|\bm \epsilon - \bm \epsilon'\|_2$, we obtain
$$
\begin{aligned}
\big|\Pi(\bm z) - \Pi(\bm z')\big|
\le\;& \Big( \|\bm k\|_2 + \|\mfhu B\|_{\rm F}\big(\|\bm p - \bm \alpha\|_2 + \|\bm h + \bm \alpha\|_2\big) + \|\mfh B\|_{\rm F}\big(\|\bm \alpha + \bm \gamma\|_2 + \|\bm h + \bm \alpha\|_2\big) \Big) \|\bm \xi - \bm \xi'\|_2 \\
& + \big(\|\bm p - \bm \alpha\|_2 + \|\bm h + \bm \alpha\|_2\big) \|\bmu \epsilon - \bmu \epsilon'\|_2
+ \big(\|\bm \alpha + \bm \gamma\|_2 + \|\bm h + \bm \alpha\|_2\big) \|\bm \epsilon - \bm \epsilon'\|_2 \\
\le\;& \bar C(\mfhu B, \mfh B)\, \|\bm z - \bm z'\|_\oplus,
\end{aligned}
$$
where the last inequality uses that a function whose increments are bounded blockwise by $L_\xi \|\bm \xi - \bm \xi'\|_2 + L_{\uline{\epsilon}} \|\bmu \epsilon - \bmu \epsilon'\|_2 + L_\epsilon \|\bm \epsilon - \bm \epsilon'\|_2$ is $\max\{L_\xi, L_{\uline{\epsilon}}, L_\epsilon\}$-Lipschitz with respect to $\|\cdot\|_\oplus$. Hence $\bar C(\mfhu B, \mfh B)$ is a Lipschitz constant of $\Pi(\cdot) = \Pi(\bmh \Theta, \bmu v, \bm v, \bm x, \bm y, \cdot)$ for every $(\bm x, \bm y)$. The same argument applied to the function $\pi(\bm x, \bm y, \bm \xi, \bmu d, \bm r)$ defined in the proof of \Cref{thm:out-of-sample}, whose increments satisfy
$$
\begin{aligned}
&\big|\pi(\bm x, \bm y, \bm \xi, \bmu d, \bm r) - \pi(\bm x, \bm y, \bm \xi', \bmu d', \bm r')\big| \\
&\quad \le \|\bm k\|_2 \|\bm \xi - \bm \xi'\|_2 + \big(\|\bm p - \bm \alpha\|_2 + \|\bm h + \bm \alpha\|_2\big)\|\bmu d - \bmu d'\|_2 \\
&\qquad + \big(\|\bm \alpha + \bm \gamma\|_2 + \|\bm h + \bm \alpha\|_2\big)\|\bm r - \bm r'\|_2,
\end{aligned}
$$
shows that $\bar L$ is a Lipschitz constant of $\pi(\bm x, \bm y, \cdot)$ with respect to $\|\cdot\|_\oplus$ on $(\bm \xi, \bmu d, \bm r)$, for every $(\bm x, \bm y)$. Finally, the third entry in the maximum defining $\bar C(\mfhu B, \mfh B)$ is bounded, by the triangle inequality for the Frobenius norm, by the corresponding entry of $\bar C(\mfu B, \mf B)$ plus $\big(\|\bm p - \bm \alpha\|_2 + \|\bm h + \bm \alpha\|_2\big)\|\mfu B - \mfhu B\|_{\rm F} + \big(\|\bm \alpha + \bm \gamma\|_2 + \|\bm h + \bm \alpha\|_2\big)\|\mf B - \mfh B\|_{\rm F} \le \bar L \big(\|\mfu B - \mfhu B\|_{\rm F} + \|\mf B - \mfh B\|_{\rm F}\big)$, while the first two entries do not depend on the estimate. Using the bound $\|\mfu B - \mfhu B\|_{\rm F} + \|\mf B - \mfh B\|_{\rm F} \le (\sqrt{\|\bmu \Sigma^{-1}\|} + \sqrt{\|\bm \Sigma^{-1}\|}) \Gamma_1$ established at the end of the proof of \Cref{thm:out-of-sample}, we conclude
\begin{equation}\label{eq:Cbar_bound}
\bar C(\mfhu B, \mfh B) \le \bar C(\mfu B, \mf B) + \bar L \big(\sqrt{\|\bmu \Sigma^{-1}\|} + \sqrt{\|\bm \Sigma^{-1}\|}\big) \Gamma_1 .
\end{equation}

We next compare the empirical and the true objectives uniformly over $(\bm x, \bm y)$. For $\bm y \in \cy$, let $\bbq^\star(\bm y)$ and $\hat \bbq(\bm y)$ denote the distributions $\bbq^\star$ and $\hat \bbq$ of the proof of \Cref{thm:out-of-sample} with $\bm y$ in place of $\hat{\bm y}$, and let $R(\bm y)$ denote the constant $R$ with $\bm y$ in place of $\hat{\bm y}$, so that $R(\bm y) \le \bar R$ for all $\bm y \in \cy$. The derivation of inequality~\eqref{eq: bound qstar and qhat} uses no property of $\hat{\bm y}$ other than $\hat{\bm y} \in \cy$; repeating it with $\bm y$ in place of $\hat{\bm y}$ shows that, on $\mathcal{E}_2$,
\begin{equation}\label{eq:unif_wass}
\Delta\big(\bbq^\star(\bm y), \hat \bbq(\bm y)\big) \le R(\bm y)\Gamma_1 + D \Gamma_2 \le \bar R \Gamma_1 + D \Gamma_2 \qquad \forall \bm y \in \cy .
\end{equation}
By the definitions of $\bbq^\star(\bm y)$ and $\hat \bbq(\bm y)$ and the identity $\pi(\bm x, \bm y, \bm \xi, \bmu d, \bm r) = \Pi(\bm \Theta, \bmu v, \bm v, \bm x, \bm y, \bm \xi, \bmu d - \mfu A \bmu v - \mfu B(\bm y \odot \bm \xi), \bm r - \mf A \bm v - \mf B(\bm y \odot \bm \xi))$, together with its counterpart for $\bmh \Theta$, we have for every $(\bm x, \bm y) \in \cx \times \cy$
$$
\epstar{\Pi(\bm \Theta, \bmu v, \bm v, \bm x, \bm y, \bmt z)} = \mathbb{E}_{\bbq^\star(\bm y)}\big[\pi(\bm x, \bm y, \bmt \xi, \bmtu d, \bmt r)\big], \qquad
\ephat{\Pi(\bmh \Theta, \bmu v, \bm v, \bm x, \bm y, \bmt z)} = \mathbb{E}_{\hat \bbq(\bm y)}\big[\pi(\bm x, \bm y, \bmt \xi, \bmtu d, \bmt r)\big].
$$
Since $\pi(\bm x, \bm y, \cdot)$ is $\bar L$-Lipschitz with respect to $\|\cdot\|_\oplus$, the Kantorovich--Rubinstein representation used in the proof of \Cref{prop:bound} and inequality~\eqref{eq:unif_wass} give, on $\mathcal{E}_2$,
\begin{equation}\label{eq:unif_approx}
\Big| \epstar{\Pi(\bm \Theta, \bmu v, \bm v, \bm x, \bm y, \bmt z)} - \ephat{\Pi(\bmh \Theta, \bmu v, \bm v, \bm x, \bm y, \bmt z)} \Big| \le \bar L \big(\bar R \Gamma_1 + D \Gamma_2\big) =: \eta \qquad \forall (\bm x, \bm y) \in \cx \times \cy .
\end{equation}
If two functions on $\cx \times \cy$ differ by at most $\eta$ pointwise, their suprema differ by at most $\eta$. Taking suprema in~\eqref{eq:unif_approx} therefore yields, on $\mathcal{E}_2$,
\begin{equation}\label{eq:eto_value}
\big| \pi^\star_{\code{CSO}} - \hat{\pi}^\star_{\code{ETO}} \big| \le \eta .
\end{equation}
Here $\pi^\star_{\code{CSO}}$ is finite: since $(\bm x - \bmt d)^+ \ge \bm x - \bmt d$ componentwise and $\bm h + \bm \alpha \ge \bm 0$, the profit satisfies $\Pi(\bm \Theta, \bmu v, \bm v, \bm x, \bm y, \bmt z) \le (\bm p - \bm \alpha)^\top \bmtu d - (\bm \alpha + \bm \gamma)^\top \bmt r + (\bm h + \bm \alpha)^\top \bmt d - (\bm h + \bm c)^\top \bm x + \|\bm k\|_2 \|\bmt \xi\|_2$, the term $-(\bm h + \bm c)^\top \bm x$ is nonpositive because $\bm x \ge \bm 0$ and $\bm h + \bm c > \bm 0$, and $\bmt z$ has finite first moments under \Cref{assm: light tail}, so the expectation of the right-hand side is bounded by a finite constant that does not depend on $(\bm x, \bm y)$, because $\|\mfu B(\bm y \odot \bmt \xi)\|_2 \le \|\mfu B\|_{\rm F}\|\bmt \xi\|_2$ and $\|\mf B(\bm y \odot \bmt \xi)\|_2 \le \|\mf B\|_{\rm F}\|\bmt \xi\|_2$ for binary $\bm y$.

We now bound the fragility of the~\ref{model::ETR} solution. Let $(\bmh x_{\code{ETO}}, \bmh y_{\code{ETO}})$ be an optimal solution of~\ref{model::ETO}. By the first part of the proof, $\bar C(\mfhu B, \mfh B)$ is a Lipschitz constant of $\Pi(\bmh \Theta, \bmu v, \bm v, \bmh x_{\code{ETO}}, \bmh y_{\code{ETO}}, \cdot)$ with respect to $\|\cdot\|_\oplus$, and \Cref{prop:bound} allows any such constant in place of $C_{\code{ETO}}$. Hence, whenever $\tau \le \hat{\pi}^\star_{\code{ETO}}$, \Cref{prop:bound} applied with $C_{\code{ETO}} = \bar C(\mfhu B, \mfh B)$ gives
\begin{equation}\label{eq:kappa_unif}
\hat{\kappa}^\star_\tau \le \bar C(\mfhu B, \mfh B) .
\end{equation}

We now bound the suboptimality of the~\ref{model::ETR} solution for fixed $S$. On $\mathcal{E}_2$, \Cref{thm:out-of-sample} gives
$$
\epstar{\Pi(\bm \Theta, \bmu v, \bm v, \hat{\bm x}, \hat{\bm y}, \bmt z)} \ge \tau - \hat{\kappa}^\star_\tau \Big(D + \big(\sqrt{\|\bmu \Sigma^{-1}\|} + \sqrt{\|\bm \Sigma^{-1}\|}\big)\Gamma_1\Big)\big(R \Gamma_1 + D \Gamma_2\big) .
$$
Since $R = R(\hat{\bm y}) \le \bar R$ and, by~\eqref{eq:kappa_unif}, $\hat{\kappa}^\star_\tau \le \bar C(\mfhu B, \mfh B)$, the subtracted term is at most $\bar C(\mfhu B, \mfh B)\big(D + (\sqrt{\|\bmu \Sigma^{-1}\|} + \sqrt{\|\bm \Sigma^{-1}\|})\Gamma_1\big)(\bar R \Gamma_1 + D \Gamma_2)$. Moreover, writing $\tau = \hat{\pi}^\star_{\code{ETO}} - (\hat{\pi}^\star_{\code{ETO}} - \tau)$ and applying~\eqref{eq:eto_value},
$$
\tau \ge \pi^\star_{\code{CSO}} - \bar L\big(\bar R \Gamma_1 + D \Gamma_2\big) - \big(\hat{\pi}^\star_{\code{ETO}} - \tau\big) .
$$
Combining the last two displays, and using that $\pi^\star_{\code{CSO}}$ is a supremum for the left inequality, we obtain on $\mathcal{E}_2$
\begin{equation}\label{eq:subopt_bound}
\begin{aligned}
0 \le\;& \pi^\star_{\code{CSO}} - \epstar{\Pi(\bm \Theta, \bmu v, \bm v, \hat{\bm x}, \hat{\bm y}, \bmt z)} \\
\le\;& \big(\hat{\pi}^\star_{\code{ETO}} - \tau\big) + \bar L \big(\bar R \Gamma_1 + D \Gamma_2\big)
+ \bar C(\mfhu B, \mfh B) \Big(D + \big(\sqrt{\|\bmu \Sigma^{-1}\|} + \sqrt{\|\bm \Sigma^{-1}\|}\big)\Gamma_1\Big)\big(\bar R \Gamma_1 + D \Gamma_2\big).
\end{aligned}
\end{equation}

To prove the proposition, we restore the index $S$ and choose the confidence level $\delta_S := 1/S$. Let $\Gamma_{2,S}$ denote $\Gamma_2$ evaluated at $\delta = \delta_S$. For all $S$ large enough that $S \ge \log(2\beta_1 S)/\beta_2$, the first branch in the definition of $\Gamma_2$ applies and
$$
\Gamma_{2,S} = \left(\frac{\log(2\beta_1 S)}{\beta_2 S}\right)^{1/(3N)} \longrightarrow 0 \qquad \text{as } S \to \infty,
$$
deterministically. Let $\mathcal{E}_{2,S}$ denote the concentration event at level $\delta_S$, so that $\mathbb{P}^S[\mathcal{E}_{2,S}] \ge 1 - 1/(2S) \to 1$, and let $\mathcal{E}_{3,S}$ denote the event in \Cref{assm:design}, so that $\mathbb{P}^S[\mathcal{E}_{3,S}] \to 1$. On $\mathcal{E}_{3,S}$ we have $\|\bmu \Sigma_S^{-1}\| \le 1/\lambda$, $\|\bm \Sigma_S^{-1}\| \le 1/\lambda$, and $\bar R_S \le 1/\lambda$. On $\mathcal{E}_{2,S} \cap \mathcal{E}_{3,S}$, inequality~\eqref{eq:eto_value} gives $|\hat{\pi}^\star_{\code{ETO},S} - \pi^\star_{\code{CSO}}| \le \bar L(\Gamma_{1,S}/\lambda + D \Gamma_{2,S})$, so for any $\varepsilon > 0$
$$
\mathbb{P}^S\big[ |\hat{\pi}^\star_{\code{ETO},S} - \pi^\star_{\code{CSO}}| > \varepsilon \big]
\le \mathbb{P}^S[\mathcal{E}_{2,S}^c] + \mathbb{P}^S[\mathcal{E}_{3,S}^c] + \mathbb{P}^S\big[ \bar L(\Gamma_{1,S}/\lambda + D \Gamma_{2,S}) > \varepsilon \big] .
$$
The first two terms tend to zero, and the third tends to zero because $\Gamma_{1,S} \to 0$ in probability and $\Gamma_{2,S} \to 0$ deterministically. Hence $\hat{\pi}^\star_{\code{ETO},S} \to \pi^\star_{\code{CSO}}$ in probability. For the out-of-sample profit, on $\mathcal{E}_{2,S} \cap \mathcal{E}_{3,S}$ the right-hand side of~\eqref{eq:subopt_bound} at level $\delta_S$, combined with~\eqref{eq:Cbar_bound}, is at most
$$
\bar G_S := \big(\hat{\pi}^\star_{\code{ETO},S} - \tau_S\big) + \bar L\big(\Gamma_{1,S}/\lambda + D \Gamma_{2,S}\big) + \Big(\bar C(\mfu B, \mf B) + \tfrac{2\bar L}{\sqrt{\lambda}} \Gamma_{1,S}\Big)\Big(D + \tfrac{2}{\sqrt{\lambda}} \Gamma_{1,S}\Big)\big(\Gamma_{1,S}/\lambda + D \Gamma_{2,S}\big) .
$$
Each of $\hat{\pi}^\star_{\code{ETO},S} - \tau_S$ and $\Gamma_{1,S}$ converges to zero in probability by assumption, $\Gamma_{2,S} \to 0$ deterministically, and the remaining quantities are finite constants; since sums and products of sequences converging in probability converge in probability to the corresponding sums and products of their limits, $\bar G_S \to 0$ in probability. Therefore, for any $\varepsilon > 0$,
$$
\mathbb{P}^S\Big[ \pi^\star_{\code{CSO}} - \epstar{\Pi(\bm \Theta, \bmu v, \bm v, \hat{\bm x}_S, \hat{\bm y}_S, \bmt z)} > \varepsilon \Big]
\le \mathbb{P}^S[\mathcal{E}_{2,S}^c] + \mathbb{P}^S[\mathcal{E}_{3,S}^c] + \mathbb{P}^S[\bar G_S > \varepsilon] \longrightarrow 0 .
$$
Together with the left inequality in~\eqref{eq:subopt_bound}, which holds surely, this proves that the out-of-sample profit of $(\hat{\bm x}_S, \hat{\bm y}_S)$ converges to $\pi^\star_{\code{CSO}}$ in probability.
$\hfill\Halmos$

\subsection{Proof of \texorpdfstring{\Cref{prop::P_rsr_refo}}{Proposition~\ref{prop::P_rsr_refo}} \label{app:P_rsr_refo}}

We first verify the equivalence of~\ref{model::ETR} and Problem~\eqref{model::rs_lift}. Write $\bmh z_s := (\bmh \xi_s, \bmhu \epsilon_s, \bmh \epsilon_s)$. For any $\bbq \in \bar{\cq}$ with marginal $\bbp$ of $(\bmt\xi, \bmtu\epsilon, \bmt\epsilon)$, the pair $(\bmt z, \bmh z_{\tilde s})$ is a coupling of $\bbp$ and $\hat{\bbp}$, so $\Delta(\bbp, \hat{\bbp}) \le \mathbb{E}_{\bbq}[\|\bmt z - \bmh z_{\tilde s}\|_\oplus] \le \mathbb{E}_{\bbq}[\tilde w]$, and, since $\kappa \geq 0$, the constraint of~\eqref{model::ETR} implies that of~\eqref{model::rs_lift}. Conversely, for any $\bbp \in \mathcal{M}(\cz)$ and any coupling of $\bbp$ and $\hat{\bbp}$, setting $\tilde w := \|\bmt z - \bmh z_{\tilde s}\|_\oplus$ defines an element of $\bar{\cq}$, and taking the infimum over couplings shows that the constraint of~\ref{model::rs_lift} implies that of~\eqref{model::ETR}. 

It then suffices to reformulate the following constraint in Problem~\eqref{model::rs_lift}: $$\mathbb{E}_{\mathbb{Q}}\Big[\Pi(\mfhu A, \mfhu B, \mfh A, \mfh B, \bmu v, \bm v, \bm x, \bm y, \bmt \xi, \bmtu \epsilon, \bmt \epsilon) + \kappa \tilde{w}\Big] \ge \tau   \quad \forall \mathbb{Q} \in \bar{\cq}.$$  For this purpose, we first condition on the lifting variable $\tilde{s}$ to write the above constraint as 			
\begin{equation}\label{original_formulation_n}
    \displaystyle \frac{1}{S}\sum_{s\in [S]} \eps{\Pi(\mfhu A, \mfhu B, \mfh A, \mfh B, \bmu v, \bm v, \bm x, \bm y, \bmt \xi, \bmtu \epsilon, \bmt \epsilon) + \kappa\tilde{w}} \ge \tau  \quad \forall \bbp_s \in \bar{\cq}_s, s \in [S] 
\end{equation}
where we denote by $\bbp_s$ the conditional probability distribution of $(\bmt \xi, \bmtu \epsilon, \bmt \epsilon, \tilde w)$ given $\tilde s = s$, and \begin{equation*}
    \Bar{\cq}_s := \left\{ 
    \bbp \in \mathcal{M}(\mathbb{R}^{3N+1}) \left| 
    \begin{array}{ll}
        (\bmt\xi, \bmtu\epsilon, \bmt \epsilon, \tilde{w}) \sim \bbp &  \\
        \bbp[(\bmt\xi, \bmtu\epsilon, \bmt \epsilon, \tilde{w}) \in \widebar{\mathcal{Z}}_s  ]= 1 \\
    \end{array}
    \right.
    \right\}
\end{equation*} By introducing auxiliary variable $t_s \in \bbr$ for each $s \in [S]$, we rewrite the constraint \eqref{original_formulation_n} as 
\begin{equation}
    \begin{aligned}
        &\frac{1}{S}\sum_{s\in [S]} t_s \ge \tau \\
        &   \eps{\Pi(\mfhu A, \mfhu B, \mfh A, \mfh B, \bmu v, \bm v, \bm x, \bm y, \bmt \xi, \bmtu \epsilon, \bmt \epsilon) + \kappa\tilde{w}} \ge t_s & \forall s \in [S], \bbp_s \in \bar{\cq}_s
    \end{aligned}
\end{equation}
The second constraint can be further reformulated as 
$$
\begin{aligned}
    & \inf_{\bbp_s \in \bar{\cq}_s}  \eps{\Pi(\mfhu A, \mfhu B, \mfh A, \mfh B, \bmu v, \bm v, \bm x, \bm y, \bmt \xi, \bmtu \epsilon, \bmt \epsilon) + \kappa\tilde{w}} \ge t_s & \forall s \in [S] \\
    \Longleftrightarrow & \inf_{(\bm\xi, \bmu\epsilon, \bm \epsilon, w) \in \widebar{\mathcal{Z}}_s} \left\{ \Pi(\mfhu A, \mfhu B, \mfh A, \mfh B, \bmu v, \bm v, \bm x, \bm y, \bm \xi, \bmu \epsilon, \bm \epsilon) + \kappa w \right\} \ge t_s & \forall s \in [S] \\
    \Longleftrightarrow  &  \inf_{(\bm\xi, \bmu\epsilon, \bm \epsilon, w) \in \widebar{\mathcal{Z}}_s} \left\{ (\bm p - \bm \alpha)^\top \big(\mfhu A\, \bmu v + \mfhu B ( \bm y \odot \bm \xi) + \bmu{\epsilon}\big) - (\bm \alpha + \bm \gamma)^\top \big(\mfh A\, \bm v + \mfh B ( \bm y \odot \bm \xi) + \bm{\epsilon}\big) \right.  \\
    &~~~~~~ - (\bm h + \bm \alpha)^\top (\bm x - \mfhu A\, \bmu v - \mfh A\, \bm v - (\mfhu B + \mfh B) ( \bm y \odot \bm \xi) - \bmu{\epsilon} - \bm{\epsilon})^+ + (\bm \alpha - \bm c)^\top \bm x \left. - \bm k^\top(\bm y \odot \bm \xi)+ \kappa w\right\} \ge t_s& \forall s \in [S] \\
\end{aligned}
$$
where the first equivalence holds because $\bar{\cq}_s$ consists of all distributions supported on $\widebar{\mathcal{Z}}_s$, including the Dirac distributions, so the infimum of the expectation equals the pointwise infimum over $\widebar{\mathcal{Z}}_s$, and the second is by definition of $\Pi$.  In summary, Problem~\eqref{model::rs_lift} can be reformulated as Problem~\eqref{model::robust_opt_m1}. 
\hfill\Halmos

\subsection{Proof of \texorpdfstring{\Cref{thm::P_refo}}{Theorem~\ref{thm::P_refo}} \label{app:P_refo}}
Model~\eqref{model::robust_opt_m1} can be equivalently reformulated as the following mixed integer conic optimization problem. 
\begin{align}\label{model::exact}
    \min~& \kappa \notag \\
    {\rm s.t.}~& \bm t^\top \bm 1 \ge \tau S \notag\\
    & \displaystyle \Big(H(\ci)+ (\bm \eta_{s, \ci})^\top \bmh \xi_s + (\bm\theta_{s, \ci})^\top \bmhu\epsilon_{s} + (\bm\rho_{s, \ci})^\top \hat{\bm\epsilon}_{s} + (\bm\chi_{s, \ci}^{b})^\top \bmhu\epsilon^{b} \notag \\ 
    &~~~~  - (\bm\chi_{s, \ci}^{u})^\top \bmhu\epsilon^{u} + (\bm\varsigma_{s, \ci}^{b})^\top\bmh\epsilon^{b} - (\bm\varsigma_{s, \ci}^{u})^\top\bmh\epsilon^{u}
    \Big) \ge t_s &&\; \forall s \in [S], \ci \subseteq [N]  \notag\\ 
    &\displaystyle \sum_{n \in \ci}\Big( (p_n+h_n)(\bm y \odot \bmhu b_n) + (h_n-\gamma_n)(\bm y \odot \bmh b_n) \Big) &&&\tag{$\code{EXACT}$} \\
    &~~~~ + \displaystyle \sum_{n \in [N]\setminus \ci}\Big( (p_n-\alpha_n)(\bm y \odot \bmhu b_n)  - (\alpha_n+\gamma_n)(\bm y \odot \bmh b_n) \Big) - \bm y \odot \bm k \ge \bm \eta_{s, \ci} && \; \forall s \in [S], \ci \subseteq [N] \notag \\
    &\|\bm \eta_{s, \ci} \|_* \le \kappa,\;
    \|\bm \theta_{s, \ci} \|_* \le \kappa, \;
    \|\bm \rho_{s, \ci} \|_* \le \kappa && \; \forall s \in [S], \ci \subseteq [N] \notag \\
    & \bm \eta_{s, \ci} \in \mathbb{R}^N, \; (\bm\theta_{s,\ci}, \bm{\rho}_{s,\ci}, \bm\chi_{s, \ci}^{b}, \bm\chi_{s, \ci}^{u}, \bm\varsigma_{s, \ci}^{b}, \bm\varsigma_{s, \ci}^{u}) \in \Psi_{s, \ci} && \;  \forall s \in [S], \ci \subseteq [N] \notag \\
    & \bm x \in \cx, \; \bm y \in \cy, \; \kappa \in \bbr_+, \; \bm t \in \mathbb{R}^S. \notag
\end{align}
Here, for each  $s\in[S], \; \ci \subseteq [N] $, we denote 
$$
\begin{aligned}
    H(\ci):=  & \sum_{n \in \ci}\Big( (p_n+h_n)(\bmhu a_n)^\top \bmu v + (h_n-\gamma_n)(\bmh a_n)^\top \bm v - (h_n+c_n)x_n \Big) \notag\\
    & + \displaystyle \sum_{n \in [N]\setminus \ci}\Big( (p_n-\alpha_n)(\bmhu a_n)^\top \bmu v  - (\alpha_n+\gamma_n)(\bmh a_n)^\top \bm v + (\alpha_n-c_n)x_n \Big) \notag\\
\end{aligned}
$$
and 
\begin{equation*}
    \Psi_{s, \ci} \coloneq \left\{
    \begin{array}{c@{\;}l@{\;}l@{\;}l@{\;}}
        &(\bm\theta_{s,\ci}, \bm{\rho}_{s,\ci}, \bm\chi_{s, \ci}^{b}, \bm\chi_{s, \ci}^{u}, \bm\varsigma_{s, \ci}^{b}, \bm\varsigma_{s, \ci}^{u}) \in \bbr^{2N} \times \bbr_+^{4N} : \\
        &p_n+h_n  = \theta_{s, \ci, n}-\chi_{s, \ci, n}^{u} + \chi_{s, \ci, n}^{b} &  \forall n \in \ci\notag\\
        &h_n-\gamma_n  = \rho_{s, \ci, n} -\varsigma_{s, \ci, n}^{u} + \varsigma_{s, \ci, n}^{b}&  \forall n \in \ci \notag\\
        &p_n-\alpha_n  = \theta_{s, \ci, n} -\chi_{s, \ci, n}^{u} + \chi_{s, \ci, n}^{b} & \forall n \in [N]\setminus \ci\notag\\
        &-\alpha_n-\gamma_n  =  \rho_{s, \ci, n} -\varsigma_{s, \ci, n}^{u} + \varsigma_{s, \ci, n}^{b} &  \forall n \in [N]\setminus \ci \notag\\
    \end{array}
    \right\}.
\end{equation*}

To prove the theorem, we focus on the reformulation of the robust constraint
$$
\begin{aligned}
    & (\bm p - \bm \alpha)^\top \big(\mfhu A\, \bmu v + \mfhu B ( \bm y \odot \bm \xi) + \bmu{\epsilon}\big) - (\bm \alpha + \bm \gamma)^\top \big(\mfh A\, \bm v + \mfh B ( \bm y \odot \bm \xi) + \bm{\epsilon}\big)   \\
    & - (\bm h + \bm \alpha)^\top (\bm x - \mfhu A\, \bmu v - \mfh A\, \bm v - (\mfhu B + \mfh B) ( \bm y \odot \bm \xi) - \bmu{\epsilon} - \bm{\epsilon})^+ + (\bm \alpha - \bm c)^\top \bm x 
    - \bm k^\top(\bm y \odot \bm \xi) + \kappa w \ge t_s  \forall (\bm \xi, \bmu \epsilon, \bm \epsilon, w) \in \widebar{\mathcal{Z}}_s, s \in [S].
\end{aligned} 
$$
By the expansion in Equation~\eqref{eq: expand exponential pieces}, we can rewrite the above constraints as 
$$
\begin{aligned}
    &\inf_{(\bm \xi, \bmu \epsilon, \bm \epsilon, w) \in \widebar{\mathcal{Z}}_s} \left\{ (\bm p - \bm \alpha)^\top \big(\mfhu A\, \bmu v + \mfhu B ( \bm y \odot \bm \xi) + \bmu{\epsilon}\big) - (\bm \alpha + \bm \gamma)^\top \big(\mfh A\, \bm v + \mfh B ( \bm y \odot \bm \xi) + \bm{\epsilon}\big) \right.  \\
    & -  \sum\limits_{n \in \ci} \; (h_n + \alpha_n) (x_n - \bmhu a_n^\top \bmu v - \bmh a_n^\top \bm v - (\bmhu b_n + \bmh b_n)^\top (\bm y \odot \bm \xi) - \uline{\epsilon}_n - \epsilon_n) \left. + (\bm \alpha - \bm c)^\top \bm x 
    - \bm k^\top(\bm y \odot \bm \xi) + \kappa w \right\} \ge t_s  \forall  s \in [S], \ci \subseteq [N].
\end{aligned}
$$
Next, we first write down the minimization problem explicitly for the left-hand side of the above constraint $s \in [S]$ and $\ci \subseteq [N]$, 
\begin{align*}
        \inf \Bigg\{ &\displaystyle\sum_{n \in \ci}\Big( (p_n+h_n)(\bmhu a_n)^\top \bmu v + (h_n-\gamma_n)(\bmh a_n)^\top \bm v - (h_n+c_n)x_n \Big) \\& ~~~~~ + \displaystyle\sum_{n \in [N]\setminus \ci}\Big( (p_n-\alpha_n)(\bmhu a_n)^\top \bmu v  - (\alpha_n+\gamma_n)(\bmh a_n)^\top \bm v + (\alpha_n-c_n)x_n \Big) \\& ~~~~~ + \Bigg( \displaystyle\sum_{n \in \ci}\Big( (p_n+h_n)(\bm y \odot \bmhu b_n) + (h_n-\gamma_n)(\bm y \odot \bmh b_n) \Big) \\& ~~~~~ + \displaystyle\sum_{n \in [N]\setminus \ci}\Big( (p_n-\alpha_n)(\bm y \odot \bmhu b_n)  - (\alpha_n+\gamma_n)(\bm y \odot \bmh b_n) \Big) - \bm y \odot \bm k \Bigg)^\top \bm \xi  \\& ~~~~~ + \displaystyle\sum_{n \in \ci}\Big( (p_n+h_n)\uline\epsilon_n + (h_n-\gamma_n)\epsilon_n \Big) + \displaystyle\sum_{n \in [N]\setminus \ci}\Big( (p_n-\alpha_n)\uline\epsilon_n - (\alpha_n+\gamma_n)\epsilon_n \Big)  + \kappa w \Bigg\}   \\
        {\rm s.t.}~&\|\bm \varphi_s^0\| + \| \bm \varphi_s^1 \| + \| \bm \varphi_s^2 \| \le w\\  
        &\bm \varphi_{s}^0 = \bm \xi - \bmh \xi_s;\;  
        \bm \varphi_{s}^1 = \bmu \epsilon - \bmhu \epsilon_s;\;  
        \bm \varphi_{s}^2 = \bm \epsilon - \bmh \epsilon_s \\ 
        & \bmu \epsilon \ge \bmhu\epsilon^{b};\;  
         \bmu \epsilon \le \bmhu\epsilon^{u};\;  
         \bm \epsilon \ge \bmh\epsilon^{b};\;  
         \bm \epsilon \le \bmh\epsilon^{u}\\ 
        & \bm \xi \ge \bm 0,  
         \bm \xi, \bmu \epsilon, \bm \epsilon, \bm \varphi_s^0, \bm \varphi_s^1, \bm \varphi_s^2 \in \mathbb{R}^N, w \in \bbr
\end{align*}
The problem is a convex conic program whose only nonlinear constraint, the norm constraint, admits a strictly feasible point, namely any $(\bm \xi, \bmu \epsilon, \bm \epsilon) \in \cz$ together with $w$ larger than the sum of the three norms; hence strong duality holds, and the dual optimum is attained whenever the primal value is finite, while an unbounded primal corresponds to an infeasible dual, so the equivalence stated below holds in both cases. Correspondingly, the dual problem is given by:
\begin{align}\label{model::dual problem_n}
        \sup \Big\{ &\displaystyle \sum_{n \in \ci}\Big( (p_n+h_n)(\bmhu a_n)^\top \bmu v + (h_n-\gamma_n)(\bmh a_n)^\top \bm v - (h_n+c_n)x_n \Big) \notag\\
        &~~~~ + \displaystyle \sum_{n \in [N]\setminus \ci}\Big( (p_n-\alpha_n)(\bmhu a_n)^\top \bmu v  - (\alpha_n+\gamma_n)(\bmh a_n)^\top \bm v + (\alpha_n-c_n)x_n \Big) \notag\\
        &~~~~  + \bm \eta_{s, \ci}^\top \bmh \xi_s + (\bm\theta_{s,\ci})^\top \hat{\bmu\epsilon}_{s} + (\bm\rho_{s, \ci})^\top \hat{\bm\epsilon}_{s} + (\bm\chi_{s, \ci}^{b})^\top \bmhu\epsilon^{b} - (\bm\chi_{s, \ci}^{u})^\top \bmhu\epsilon^{u} + (\bm\varsigma_{s, \ci}^{b})^\top\bmh\epsilon^{b} - (\bm\varsigma_{s, \ci}^{u})^\top\bmh\epsilon^{u} \Big\} \notag\\
        {\rm s.t.}~& \displaystyle \sum_{n \in \ci}\Big( (p_n+h_n)(\bm y \odot \bmhu b_n) + (h_n-\gamma_n)(\bm y \odot \bmh b_n) \Big) \notag\\
        &~~~~ + \displaystyle \sum_{n \in [N]\setminus \ci}\Big( (p_n-\alpha_n)(\bm y \odot \bmhu b_n)  - (\alpha_n+\gamma_n)(\bm y \odot \bmh b_n) \Big) - \bm y \odot \bm k \ge \bm \eta_{s, \ci}  \\
        &p_n+h_n  = \theta_{s, \ci, n}-\chi_{s, \ci, n}^{u} + \chi_{s, \ci, n}^{b};\; h_n-\gamma_n  = \rho_{s, \ci, n} -\varsigma_{s, \ci, n}^{u} + \varsigma_{s, \ci, n}^{b}& \hspace{-1.9cm} \forall n \in \ci \notag\\
        &p_n-\alpha_n  = \theta_{s, \ci, n} -\chi_{s, \ci, n}^{u} + \chi_{s, \ci, n}^{b};\; -\alpha_n-\gamma_n  =  \rho_{s, \ci, n} -\varsigma_{s, \ci, n}^{u} + \varsigma_{s, \ci, n}^{b} & \hspace{-1.9cm} \forall n \in [N]\setminus \ci \notag\\
        &\|\bm \eta_{s, \ci} \|_* \le \kappa;\; \|\bm \theta_{s, \ci} \|_* \le \kappa;\; \|\bm \rho_{s, \ci} \|_* \le \kappa  \notag\\
        &\bm\eta_{s, \ci}, \bm\theta_{s, \ci}, \bm\rho_{s, \ci} \in \mathbb{R}^N; \bm\chi_{s, \ci}^{b}, \bm\chi_{s, \ci}^{u}, \bm\varsigma_{s, \ci}^{b}, \bm\varsigma_{s, \ci}^{u} \in \bbr_+^N, \kappa \in \bbr_+. \notag
\end{align}
By strong duality, the constraint that the infimum be at least $t_s$ holds if and only if there exist dual variables feasible for~\eqref{model::dual problem_n} whose objective value is at least $t_s$. Collecting these conditions over $s \in [S]$ and $\ci \subseteq [N]$ yields~\eqref{model::exact}.
$\hfill\Halmos$

\subsection{Proof of Theorem \ref{Theo::LDR_refo}} \label{app:LDR_refo}

First, the infimum over $\widebar{\cz}_s$ of the left-hand side of the second constraint in problem \eqref{model::bar_ETR} equals the optimal value of the following problem:
\begin{align*}
        \inf ~& \Big\{(\bm p-\bm \alpha)^\top (\mfhu A \, \bmu v) - (\bm \alpha+\bm \gamma)^\top (\mfh A \, \bm v) + (\bm \alpha-\bm c)^\top \bm x -(\bm h+\bm \alpha)^\top \bm{\ell}_{s}^0 \\
        &~~ + \left( ( \mfhu B^\top (\bm p-\bm \alpha) -  \mfh B^\top (\bm \alpha+\bm \gamma) - \bm k)  \odot \bm y - \mf{L}_{s}^{1\top} (\bm h+\bm \alpha)  \right)^\top \bm \xi \\
        &~~ + \left(\bm p - \bm \alpha - \mf{L}_{s}^{2\top}(\bm h+\bm \alpha)\right)^\top \bmu \epsilon - \left(\bm \alpha + \bm \gamma + \mf{L}_{s}^{3\top}(\bm h+\bm \alpha)\right)^\top \bm \epsilon + (\kappa - (\bm h+\bm \alpha)^\top \bm{\ell}_{s}^4) w \Big\}  \\
        {\rm s.t.}~&\|\bm \varphi_s^0\| + \| \bm \varphi_s^1 \| + \| \bm \varphi_s^2 \| \le w \\ 
        &\bm \varphi_s^0 = \bm \xi - \bmh \xi_s ;\; 
        \bm \varphi_s^1 = \bmu \epsilon - \bmhu \epsilon_s;\; 
        \bm \varphi_s^2 = \bm \epsilon - \bmh \epsilon_s\\ 
        & \bmu \epsilon \ge \bmhu\epsilon^{b};\; 
         \bmu \epsilon \le \bmhu\epsilon^{u};\;  
         \bm \epsilon \ge \bmh\epsilon^{b};\;  
         \bm \epsilon \le \bmh\epsilon^{u}\\ 
        &\bm \xi \ge \bm 0;\; 
        \bmu \epsilon,\bm \epsilon, \bm \varphi_s^0,\bm \varphi_s^1,\bm \varphi_s^2 \in \mathbb{R}^N, w \in \bbr. 
\end{align*}

As in the proof of \Cref{thm::P_refo}, the norm constraint admits a strictly feasible point and all other constraints are linear, so strong duality holds. Let $(\bm\eta^1_s, \bm\theta^1_s, \bm\rho^1_s) \in \bbr^{3N}$ denote the multipliers of the three equality constraints and $(\bm\chi^{b1}_s, \bm\chi^{u1}_s, \bm\varsigma^{b1}_s, \bm\varsigma^{u1}_s) \in \bbr_+^{4N}$ those of the four bound constraints; stationarity in $w$ gives $\kappa - (\bm h + \bm\alpha)^\top \bm\ell^4_s$ for the multiplier of the norm constraint, and stationarity in $\bm\xi$ eliminates the multiplier of $\bm\xi \ge \bm 0$. The dual problem is:
\begin{align*}
        \sup ~& \Big\{ (\bm p-\bm \alpha)^\top (\mfhu A \, \bmu v) - (\bm \alpha+\bm \gamma)^\top (\mfh A \, \bm v) + (\bm \alpha-\bm c)^\top \bm x -(\bm h+\bm \alpha)^\top \bm{\ell}_{s}^0 \\
        &~~~~  + (\bm \eta_s^1)^\top \bmh \xi_s+ (\bm \theta_s^1)^\top \bmhu \epsilon_s + (\bm \rho_s^1)^\top \bmh \epsilon_s + (\bm\chi_s^{b1})^\top \bmhu\epsilon^{b} - (\bm\chi_s^{u1})^\top \bmhu\epsilon^{u} + (\bm\varsigma_s^{b1})^\top \bmh\epsilon^{b} - (\bm\varsigma_s^{u1})^\top \bmh\epsilon^{u} \Big\}  \\ 
        {\rm s.t.}~& \Big( \mfhu B^\top (\bm p-\bm \alpha) -  \mfh B^\top (\bm \alpha+\bm \gamma) - \bm k\Big)  \odot \bm y - \mf{L}_{s}^{1\top} (\bm h+\bm \alpha) \ge \bm \eta_s^1 \\ 
        &\bm p - \bm \alpha - \mf{L}_{s}^{2\top}(\bm h+\bm \alpha) = \bm \theta_s^1 + \bm\chi_s^{b1} - \bm\chi_s^{u1}  \\ 
        &- \left(\bm \alpha + \bm \gamma + \mf{L}_{s}^{3\top}(\bm h+\bm \alpha)\right)  = \bm \rho_s^1 + \bm\varsigma_s^{b1} - \bm\varsigma_s^{u1}\\ 
        &\|\bm \eta_s^1 \|_* \le \kappa - (\bm h+\bm \alpha)^\top \bm{\ell}_{s}^4; \;
        \|\bm \theta_s^1 \|_* \le \kappa - (\bm h+\bm \alpha)^\top \bm{\ell}_{s}^4 \\ 
        &\|\bm \rho_s^1 \|_* \le \kappa - (\bm h+\bm \alpha)^\top \bm{\ell}_{s}^4; \;
         \bm\eta_s^1, \bm\theta_s^1, \bm\rho_s^1\in \mathbb{R}^N;\bm\chi_s^{b1}, \bm\chi_s^{u1}, \bm\varsigma_s^{b1}, \bm\varsigma_s^{u1} \in \bbr_+^N.
\end{align*}

Similarly, for each $n \in [N]$, dualizing the $n$-th component of the third and of the fourth constraint of~\ref{model::bar_ETR} over $\widebar{\mathcal{Z}}_s$, with multipliers $(\bm\eta^2_{s,n}, \bm\theta^2_{s,n}, \bm\rho^2_{s,n}) \in \bbr^{3N}$, $(\bm\chi^{b2}_{s,n}, \bm\chi^{u2}_{s,n}, \bm\varsigma^{b2}_{s,n}, \bm\varsigma^{u2}_{s,n}) \in \bbr_+^{4N}$ and $(\bm\eta^3_{s,n}, \bm\theta^3_{s,n}, \bm\rho^3_{s,n}) \in \bbr^{3N}$, $(\bm\chi^{b3}_{s,n}, \bm\chi^{u3}_{s,n}, \bm\varsigma^{b3}_{s,n}, \bm\varsigma^{u3}_{s,n}) \in \bbr_+^{4N}$ for the equality and bound constraints, and with the multipliers of the norm constraint and of $\bm\xi \ge \bm 0$ eliminated by stationarity as above (the former equals $\ell^4_{s,n}$), yields the remaining dual constraints.

Combining the dual reformulation and eliminating the variables $\bm \theta_{s}^1$, $\bm \rho_s^1$, $\bm \theta_{s,n}^2$, $\bm \rho_{s,n}^2$, $\bm \theta_{s,n}^3$, $\bm \rho_{s,n}^3$ by equality constraints, the robust counterpart of the \ref{model::bar_ETR} model can be formulated as follows,
\begin{align}\label{model::RS_L_CP}
    \min~& \kappa & \notag\\
    {\rm s.t.}~& \bm t^\top \bm 1 \ge \tau S \notag\\
    & (\bm p-\bm \alpha)^\top (\mfhu A \, \bmu v) - (\bm \alpha+\bm \gamma)^\top (\mfh A \, \bm v) + (\bm \alpha-\bm c)^\top \bm x -(\bm h+\bm \alpha)^\top \bm{\ell}_{s}^0   \notag \\
    &~~~~ + (\bm \eta_s^1)^\top \bmh \xi_s + (\bm p - \bm \alpha - \mf{L}_{s}^{2\top}(\bm h+\bm \alpha) - \bm\chi_s^{b1} + \bm\chi_s^{u1})^\top \bmhu \epsilon_s \notag\\
    &~~~~ - (\bm \alpha + \bm \gamma + \mf{L}_{s}^{3\top}(\bm h+\bm \alpha) + \bm\varsigma_s^{b1} - \bm\varsigma_s^{u1})^\top \bmh \epsilon_s \notag\\
    & ~~~~+ (\bm\chi_s^{b1})^\top \bmhu\epsilon^{b} - (\bm\chi_s^{u1})^\top \bmhu\epsilon^{u} + (\bm\varsigma_s^{b1})^\top \bmh\epsilon^{b} - (\bm\varsigma_s^{u1})^\top \bmh\epsilon^{u} \ge t_s &&\; \forall s \in [S] \notag \\
    & \Big( \mfhu B^\top (\bm p-\bm \alpha) -  \mfh B^\top (\bm \alpha+\bm \gamma) - \bm k\Big)  \odot \bm y - \mf{L}_{s}^{1\top} (\bm h+\bm \alpha) \ge \bm \eta_s^1&& \; \forall s \in [S] \notag\\
    &\|\bm \eta_s^1 \|_* \le \kappa - (\bm h+\bm \alpha)^\top \bm{\ell}_{s}^4 && \; \forall s \in [S] \notag\\
    &\|\bm p - \bm \alpha - \mf{L}_{s}^{2\top}(\bm h+\bm \alpha) - \bm\chi_s^{b1} + \bm\chi_s^{u1} \|_* \le \kappa - (\bm h+\bm \alpha)^\top \bm{\ell}_{s}^4 && \; \forall s \in [S] \notag\\
    &\|\bm \alpha + \bm \gamma + \mf{L}_{s}^{3\top}(\bm h+\bm \alpha) + \bm\varsigma_s^{b1} - \bm\varsigma_s^{u1} \|_* \le \kappa - (\bm h+\bm \alpha)^\top \bm{\ell}_{s}^4 && \; \forall s \in [S] \notag\\
    & \bmhu a_n^\top \bmu v + \bmh a_n^\top \bm v + \ell_{s,n}^0 - x_n + (\bm \eta_{s,n}^2)^\top \bmh \xi_s + (\bm{\ell}_{s,n}^{2}+\bm e_n - \bm\chi_{s,n}^{b2} + \bm\chi_{s,n}^{u2})^\top \bmhu \epsilon_s \notag\\
    &~~~~ + (\bm{\ell}_{s,n}^{3}+\bm e_n - \bm\varsigma_{s,n}^{b2} + \bm\varsigma_{s,n}^{u2})^\top \bmh \epsilon_s + (\bm\chi_{s,n}^{b2})^\top \bmhu\epsilon^{b} - (\bm\chi_{s,n}^{u2})^\top \bmhu\epsilon^{u} \notag \\
    &~~~~   + (\bm\varsigma_{s,n}^{b2})^\top \bmh\epsilon^{b} - (\bm\varsigma_{s,n}^{u2})^\top \bmh\epsilon^{u} \ge 0  && \; \forall n \in [N], s \in [S]  \\
    & \bm{\ell}_{s,n}^{1} + \bm y \odot(\bmhu b_n + \bmh b_n) \ge \bm \eta_{s,n}^2 && \; \forall n \in [N], s \in [S] \notag\\
    &\|\bm \eta_{s,n}^2 \|_* \le \ell_{s,n}^4 && \; \forall n \in [N], s \in [S] \notag\\
    &\| \bm{\ell}_{s,n}^{2}+\bm e_n - \bm\chi_{s,n}^{b2} + \bm\chi_{s,n}^{u2} \|_* \le \ell_{s,n}^4 && \; \forall n \in [N], s \in [S] \notag\\
    &\|\bm{\ell}_{s,n}^{3}+\bm e_n - \bm\varsigma_{s,n}^{b2} + \bm\varsigma_{s,n}^{u2} \|_* \le \ell_{s,n}^4 && \; \forall n \in [N], s \in [S] \notag\\
    &\ell_{s,n}^0 + (\bm \eta_{s,n}^3)^\top \bmh \xi_s + (\bm{\ell}_{s,n}^{2}  - \bm\chi_{s,n}^{b3} + \bm\chi_{s,n}^{u3})^\top \bmhu \epsilon_s + (\bm{\ell}_{s,n}^{3} - \bm\varsigma_{s,n}^{b3} + \bm\varsigma_{s,n}^{u3})^\top \bmh \epsilon_s \notag \\
    &~~~~ + (\bm\chi_{s,n}^{b3})^\top \bmhu\epsilon^{b} - (\bm\chi_{s,n}^{u3})^\top \bmhu\epsilon^{u} + (\bm\varsigma_{s,n}^{b3})^\top \bmh\epsilon^{b} - (\bm\varsigma_{s,n}^{u3})^\top \bmh\epsilon^{u} \ge 0  && \; \forall n \in [N], s \in [S] \notag\\ 
    & \bm{\ell}_{s,n}^{1} \ge \bm \eta_{s,n}^3;\; 
    \|\bm \eta_{s,n}^3 \|_* \le \ell_{s,n}^4 && \; \forall n \in [N], s \in [S] \notag\\  
    &\|\bm{\ell}_{s,n}^{2}  - \bm\chi_{s,n}^{b3} + \bm\chi_{s,n}^{u3} \|_* \le \ell_{s,n}^4 && \; \forall n \in [N], s \in [S] \notag\\  
    &\|\bm{\ell}_{s,n}^{3} - \bm\varsigma_{s,n}^{b3} + \bm\varsigma_{s,n}^{u3}  \|_* \le \ell_{s,n}^4 && \; \forall n \in [N], s \in [S] \notag\\
    &\bm\chi_{s,n}^{b2}, \bm\chi_{s,n}^{u2}, \bm\varsigma_{s,n}^{b2}, \bm\varsigma_{s,n}^{u2},\bm\chi_{s,n}^{b3}, \bm\chi_{s,n}^{u3}, \bm\varsigma_{s,n}^{b3}, \bm\varsigma_{s,n}^{u3} \in \bbr_+^N ; \bm \eta_{s,n}^2, \bm \eta_{s,n}^3\in \mathbb{R}^N  && \; \forall n \in [N], s \in [S] \notag \\
    &\bm{\ell}_{s}^0, \bm{\ell}_{s}^4, \bm \eta_{s}^1 \in \mathbb{R}^N;
    \mf{L}_{s}^1, \mf{L}_{s}^2, \mf{L}_{s}^3 \in \mathbb{R}^{N\times N}; \bm\chi_{s}^{b1}, \bm\chi_{s}^{u1}, \bm\varsigma_{s}^{b1}, \bm\varsigma_{s}^{u1} \in \mathbb{R}^N_+ && \; \forall s \in [S] \notag\\
    &\bm x \in \cx, \; \bm y \in \cy, \; \kappa \in \bbr_+,\; \bm t \in \bbr^S. \notag 
\end{align}
Here $\bm\ell^{k}_{s,n} \in \bbr^N$ denotes the transpose of the $n$-th row of $\mf L^{k}_s$, $k=1,2,3$, so that $\ell_{s,n}(\bm\xi, \bmu\epsilon, \bm\epsilon, w) = \ell^0_{s,n} + (\bm\ell^1_{s,n})^\top \bm\xi + (\bm\ell^2_{s,n})^\top \bmu\epsilon + (\bm\ell^3_{s,n})^\top \bm\epsilon + \ell^4_{s,n} w$. Under the representability assumptions on $\cx$ and $\cy$, the reformulation~\eqref{model::RS_L_CP} is a mixed-integer conic optimization problem
with $\mathcal{O}(SN^2)$ variables and constraints, in addition to those required in $\mathcal{X}$ and $\mathcal{Y}$.  $\hfill\Halmos$

\subsection{Proof of \texorpdfstring{\Cref{prop::optimality of single product}}{Proposition~\ref{prop::optimality of single product}}}  \label{app:P_LDR}

Let $\hat{\kappa}^\star_\tau$ and $\bar{\kappa}^\star_\tau$ be the optimal value of Problem \eqref{model::robust_opt_m1} and Problem \eqref{model::bar_ETR}, respectively. Because Problem~\eqref{model::bar_ETR} is a conservative approximation of Problem~\eqref{model::robust_opt_m1} ($\cl \subseteq \mathcal{R}^{3N+1,N}$), we have $\hat{\kappa}^\star_\tau \le \bar{\kappa}^\star_\tau$. 

When $N=1$, we show $\hat{\kappa}^\star_\tau \ge \bar{\kappa}^\star_\tau$ by proving that any feasible solution $(\ddot{x},\ddot{y},\ddot{\kappa}, \ddot{\bm t})$ of Problem~\eqref{model::robust_opt_m1} is also feasible for Problem~\eqref{model::bar_ETR}; taking the infimum over feasible $\ddot\kappa$ then gives $\bar{\kappa}^\star_\tau \le \hat{\kappa}^\star_\tau$.

Given any such solution, let $\ddot{\ell}_s(\xi, \uline{\epsilon}, \epsilon, w) := (\ddot{x} - \uline{d}(\xi, \uline{\epsilon}) - r(\xi, \epsilon))^+ \in \mathcal{R}^{4,1}$ for each $s \in [S]$, with $\uline d$ and $r$ defined below; feasibility for Problem~\eqref{model::robust_opt_m1} then states that
\begin{equation}\label{eq: feasibility of full recourse}
    \begin{array}{c@{\;}l@{\;}l@{\;}l@{\;}}
        & \ddot{\bm t}^\top \bm 1 \ge \tau S \\
        & ( p- \alpha)\cdot  \uline{d}(\xi, \uline{\epsilon}) - ( \alpha+ \gamma)\cdot   r(\xi, \epsilon) + ( \alpha- c)\cdot   \ddot x \\
        &~~~~ -( h+ \alpha)\cdot   \ddot{\ell}_s(\xi, \uline{\epsilon}, \epsilon, w) -  k\cdot  (\ddot{y} \odot \xi) + \ddot{\kappa} \cdot w \ge \ddot{t}_s &&\; \forall ( \xi, \uline{\epsilon},  \epsilon, w) \in \widebar{\cz}_s, s \in [S]\\
        & \ddot{\ell}_s( \xi, \uline{\epsilon},  \epsilon, w) \ge  \ddot{x}-\uline{d}(\xi, \uline{\epsilon})- r(\xi, \epsilon) && \;  \forall ( \xi, \uline{\epsilon},  \epsilon, w) \in \widebar{\cz}_s, s \in [S] \\
        &  \ddot{\ell}_s( \xi, \uline{\epsilon},  \epsilon, w) \ge  0 && \;  \forall ( \xi, \uline{\epsilon},  \epsilon, w) \in \widebar{\cz}_s, s \in [S] \\
    \end{array}
\end{equation} where $\uline{d}(\xi, \uline{\epsilon}) := \mfhu A \bmu v + \mfhu B(\ddot y \odot \xi) + \uline{\epsilon}$ and $r(\xi, \epsilon) := \mfh A \bm v + \mfh B(\ddot y \odot \xi) + \epsilon$ are the estimated net-demand and return mappings of Equations~\eqref{model:demand vector} and \eqref{model:return vector} evaluated at the advertising decision $\ddot y$. 

Since $h + \alpha > 0$, we may consider the affine mapping $\ring{\ell}_s \in \cl$ defined by 
$$
\ring{\ell}_s ( \xi, \uline{\epsilon},  \epsilon, w)  :=  \frac{( p- \alpha)\cdot \uline{d}(\xi, \uline{\epsilon}) - ( \alpha+ \gamma)\cdot   r(\xi, \epsilon) + ( \alpha- c)\cdot   \ddot x  -  k\cdot  (\ddot{y} \odot \xi) + \ddot{\kappa} \cdot w - \ddot{t}_s }{h+\alpha}
$$
for each $s \in [S]$. Then we have 
\begin{equation}
    \begin{array}{c@{\;}l@{\;}l@{\;}l@{\;}}
        & \ddot{\bm t}^\top \bm 1 \ge \tau S \\
        & ( p- \alpha)\cdot  \uline{d}(\xi, \uline{\epsilon}) - ( \alpha+ \gamma)\cdot   r(\xi, \epsilon) + ( \alpha- c)\cdot   \ddot x \\
        &~~~~ -( h+ \alpha)\cdot   \ring{\ell}_s(\xi, \uline{\epsilon}, \epsilon, w) -  k\cdot  (\ddot{y} \odot \xi) + \ddot{\kappa} \cdot w \ge \ddot{t}_s &&\; \forall ( \xi, \uline{\epsilon},  \epsilon, w) \in \widebar{\cz}_s, s \in [S]\\
        & \ring{\ell}_s( \xi, \uline{\epsilon},  \epsilon, w) \ge  \ddot{x}-\uline{d}(\xi, \uline{\epsilon})- r(\xi, \epsilon) && \;  \forall ( \xi, \uline{\epsilon},  \epsilon, w) \in \widebar{\cz}_s, s \in [S] \\
        &  \ring{\ell}_s( \xi, \uline{\epsilon},  \epsilon, w) \ge  0 && \;  \forall ( \xi, \uline{\epsilon},  \epsilon, w) \in \widebar{\cz}_s, s \in [S] \\
    \end{array}
\end{equation} 
where the first inequality follows from \eqref{eq: feasibility of full recourse}, the second follows from the definition of $\ring{\ell}_s$, the last two follow from both  the definition of $\ring{\ell}_s$ and \eqref{eq: feasibility of full recourse}:
$$
\begin{array}{c@{\;}l@{\;}l@{\;}l@{\;}}
    \displaystyle \ring{\ell}_s( \xi, \uline{\epsilon},  \epsilon, w) & \displaystyle = \frac{( p- \alpha)\cdot \uline{d}(\xi, \uline{\epsilon}) - ( \alpha+ \gamma)\cdot   r(\xi, \epsilon) + ( \alpha- c)\cdot   \ddot x  -  k\cdot  (\ddot{y} \odot \xi) + \ddot{\kappa} \cdot w - \ddot{t}_s }{h+\alpha} \\
    & \ge \ddot{\ell}_s( \xi, \uline{\epsilon},  \epsilon, w) \\
    & \ge  \max\left\{\ddot{x}-\uline{d}(\xi, \uline{\epsilon})- r(\xi, \epsilon) , \, 0\right\} &\hspace{-2.1cm}   \forall ( \xi, \uline{\epsilon},  \epsilon, w) \in \widebar{\cz}_s, s \in [S] 
\end{array}
$$
Hence, we have shown that there also exist feasible affine recourse mappings $\ring{\ell}_s \in \cl$, which implies that  $(\ddot{x},\ddot{y},\ddot{\kappa}, \ddot{\bm t})$ is also feasible to Problem~\eqref{model::bar_ETR}. $\hfill\Halmos$

\section{Additional Information or Results} \label{appendix:addi_comp_resu}
\subsection{Descriptive Analytics of the Data}
Table~\ref{table:weekly_stat_main} reports the descriptive statistics for the weekly data from our industrial partner. The statistics are computed on the 695 product-week observations of products P1--P10 retained after the screening described in Section~\ref{subsec::empirical evidences}, which are also the observations used in Appendices~\ref{app:prediction_models} and~\ref{app:protocol}. Campaigns are active, with positive ad spend, in 571 of these product-weeks; own ad clicks are zero in the remaining 124.
\begin{table}[htbp]
	\centering
	\caption{Descriptive statistics for the weekly data}
	\label{table:weekly_stat_main}
	\begin{threeparttable}
		\renewcommand{\arraystretch}{1}
		\setlength{\tabcolsep}{5pt}
		\fontsize{9}{9}\selectfont{
			\begin{tabular}{lrrrrrrr}
				\toprule
				Variable & Mean & S.D. & Min & 25\% & Median & 75\% & Max \\
				\midrule
				Sales quantity & 14.50 & 21.89 & 0.00 & 2.00 & 7.00 & 17.00 & 132.00 \\
				Return quantity & 3.49 & 4.74 & 0.00 & 0.00 & 2.00 & 5.00 & 27.00 \\
				Net demand & 11.01 & 18.32 & 0.00 & 1.00 & 4.00 & 11.00 & 114.00 \\
				Own ad clicks & 190.15 & 197.85 & 0.00 & 44.50 & 147.00 & 277.50 & 1492.00 \\
				Other-product ad clicks & 1704.88 & 873.77 & 236.00 & 1082.00 & 1432.00 & 1984.00 & 5268.00 \\
				Retail price & 25.71 & 17.19 & 5.69 & 17.98 & 19.99 & 26.99 & 69.99 \\
				\bottomrule
			\end{tabular}
		}
		\begin{tablenotes}\footnotesize
			\item Notes. 695 product-week observations. Own ad clicks are the product's recorded weekly clicks, zero in weeks in which it is not advertised; other-product ad clicks are the total weekly clicks of the other nine focal products. 
            Retail price is in dollars.
		\end{tablenotes}
	\end{threeparttable}
\end{table}

\subsection{Testing for Click Cannibalization}\label{app:click_exogeneity}

\Cref{assm:demand and return model} treats the distribution of ad clicks as invariant to the advertising decision $\bm y$. The competing hypothesis is attention cannibalization: if the retailer's campaigns competed for a fixed pool of shopper attention, advertising additional products would divert impressions from a focal product and depress its clicks, so a product's clicks would be negatively associated with the number of other focal products advertised in the same week. We examine this descriptively.

The sample consists of the 571 advertised product-weeks (positive ad spend) of the cleaned panel of Section~\ref{subsec::empirical evidences}, with mean weekly clicks of 231.45. We regress a product's weekly clicks on the number of other focal products advertised in the same week, which ranges from three to nine with a mean of 6.73, controlling for product fixed effects, retail price, a linear week trend, and month fixed effects; specifications (3)--(5) also control for the product's own ad spend, so that the coefficient is identified holding own campaign intensity fixed. \Cref{table:click_exogeneity} reports the results.

\begin{table}[htbp]
	\centering
	\begin{threeparttable}
		\caption{Click volume and concurrent advertising of other focal products}
		\label{table:click_exogeneity}
		\renewcommand{\arraystretch}{1}
		\setlength{\tabcolsep}{5pt}
		\fontsize{9}{9}\selectfont{
			\begin{tabular}{clllrrrcr}
				\toprule
				& Regressor & Outcome & Own spend & Coef. & S.E. & $p$ & 95\% CI & $R^2$ \\
				\midrule
				(1) & Count            & Clicks      & No  & $-1.73$   & $5.50$   & $0.753$ & $[-12.51,\ 9.04]$  & $0.539$ \\
				(2) & Count             & Clicks      & No  & $-1.73$   & $5.03$   & $0.731$ & $[-11.59,\ 8.13]$  & $0.539$ \\
				(3) & Count            & Clicks      & Yes & $3.08$    & $3.16$   & $0.331$ & $[-3.12,\ 9.28]$   & $0.843$ \\
				(4) & Count            & $\log$ clicks & Yes & $-0.0050$ & $0.0096$ & $0.603$ & $[-0.0237,\ 0.0138]$ & $0.952$ \\
				(5) & Spend (\$100s)   & Clicks      & Yes & $7.40$    & $3.12$   & $0.018$ & $[1.28,\ 13.53]$   & $0.845$ \\
				(6) & Spend (\$100s)   & Clicks      & No  & $11.53$   & $5.70$   & $0.043$ & $[0.37,\ 22.70]$   & $0.543$ \\
				\bottomrule
			\end{tabular}
		}
		\begin{tablenotes}[flushleft]\footnotesize
			\item Count is the number of other focal products advertised in the same week and Spend their total ad spend in hundreds of dollars. All specifications include product fixed effects, retail price, a linear week trend, and month fixed effects. 
            Standard errors clustered by product.
		\end{tablenotes}
	\end{threeparttable}
\end{table}

No specification supports cannibalization. In the baseline, one additional advertised product changes weekly clicks of $-1.73$ with a robust standard error of $5.50$, and clustering by product leaves the estimate unchanged. The confidence intervals are more informative than the $p$-values, because a null result could reflect low power whereas an interval bounds what the data permit: specification (4), the tightest, places the effect of one additional advertised product between $-2.3\%$ and $+1.4\%$ of own clicks. A cannibalization effect large enough to matter for the model is thus ruled out rather than merely undetected.

The two spend specifications, (5) and (6), return positive and significant coefficients. The sign is the opposite of what cannibalization implies, since weeks of heavier advertising across the retailer's products coincide with more own clicks, and most plausibly reflects common week-level variation in shopper traffic and in the retailer's advertising intensity that the month dummies and linear trend absorb only partly. We caution against a causal reading of either set of estimates: the number of advertised products is constant within a week across advertised products and hence collinear with week fixed effects, so week-level confounders cannot be differenced out. The evidence therefore bounds any click diversion among the retailer's own products to within a few percent of own clicks, which supports the invariance that \Cref{assm:demand and return model} posits; it does not establish the absence of week-level co-movement between advertising and traffic.

\subsection{Pooled IV Evidence and Instrument Relevance}\label{app:pooled_iv_evidence}

To examine the empirical implications of instrumenting advertising exposure, we estimate pooled OLS and 2SLS models for sales, returns, and net demand on a common sample. Own clicks and total clicks of the other products are the two endogenous regressors; their coefficients are pooled across products. Each outcome and click regressor is normalized by its product-specific sample mean. Both specifications include product fixed effects, retail price, a linear week trend, and month fixed effects. The excluded instruments are lagged own clicks and total advertising-market expenditure excluding the focal product. In this pooled calculation, the lag is the preceding available observation for the same product in the cleaned weekly panel. Removing observations without this lag leaves 685 product-weeks.

The heteroskedasticity-robust first-stage $F$ statistics for joint exclusion of the two instruments are $72.0$ for own clicks and $1{,}015.8$ for other-product clicks. These statistics provide evidence of instrument relevance in the pooled specification. Exclusion additionally requires that lagged clicks and market expenditure affect current outcomes only through the included advertising exposures, conditional on the controls. The first-stage statistics do not test these restrictions or establish identification of the larger product-specific system.

\begin{table}[htbp]
    \vspace{-0.4cm}
  \centering
  \caption{Pooled OLS and IV estimates of direct and spillover advertising effects}
  \label{table:ols_iv_effects}
  \begin{threeparttable}
    \renewcommand{\arraystretch}{1}
    \setlength{\tabcolsep}{5pt}
    \fontsize{9}{9}\selectfont{
    \begin{tabular}{lcccc}
      \toprule
      & \multicolumn{2}{c}{OLS} & \multicolumn{2}{c}{IV} \\
      \cmidrule(lr){2-3} \cmidrule(lr){4-5}
      Outcome & Direct & Spillover & Direct & Spillover \\
      \midrule
      Sales quantity & 0.789$^{***}$ & $-$0.360$^{***}$ & 0.985$^{***}$ & $-$0.530$^{***}$ \\
      Return quantity & 0.482$^{***}$ & $-$0.407$^{***}$ & 0.917$^{***}$ & $-$0.723$^{***}$ \\
      Net demand & 0.886$^{***}$ & $-$0.324$^{**}$ & 0.999$^{***}$ & $-$0.438$^{***}$ \\
      \bottomrule
    \end{tabular}
    }
    \begin{tablenotes}
     \item Notes. 685 product-week observations; heteroskedasticity-robust inference. $^\text{*}p<0.1$, $^{\text{**}}p<0.05$, $^{\text{***}}p<0.01$.
    \end{tablenotes}
  \end{threeparttable}
  \vspace{-0.2cm}
\end{table}

\Cref{table:ols_iv_effects} shows that all six IV advertising coefficients retain the signs of their OLS counterparts and are significant at the 1\% level. The direct effects rise from 0.789 to 0.985 for sales, from 0.482 to 0.917 for returns, and from 0.886 to 0.999 for net demand; the negative spillover effects also increase in magnitude. The change is particularly pronounced for returns, indicating that accounting for the potentially endogenous component of clicks can materially affect the estimated demand--return relationship. This comparison motivates the IV policy experiment in Section~\ref{sec:endo_extension}. The pooled estimates summarize that empirical comparison; the policy experiment retains unrestricted product-pair effects and expands the instrument set accordingly.

\subsection{Prediction Specifications and Validation}\label{app:prediction_models}

This subsection specifies the five prediction specifications compared in \Cref{table:product_specific_prediction} and the measures used to evaluate them, in the notation of Section~\ref{subsec::empirical evidences}: products are indexed by $n,m\in[N]$ and weeks by $t$; $\uline{d}_{nt}$, $r_{nt}$, and $d_{nt}=\uline{d}_{nt}+r_{nt}$ are the net demand, returns, and sales (realized demand) of product $n$ in week $t$; $\xi_{mt}$ is the recorded weekly ad clicks of product $m$, equal to zero in weeks in which $m$ is not advertised, and $\bm\xi_t=(\xi_{1t},\ldots,\xi_{Nt})^{\top}$; and $\bmu v_{nt}=\bm v_{nt}$ collects the retail price, the week index, and the eleven month indicators. Clicks are taken from the original weekly records, so a product's clicks enter the equations of the other products also in weeks in which its own outcome observation was removed by the screening of Section~\ref{subsec::empirical evidences}; an absent click record contributes zero. Throughout, $(\uline{\hat{d}}_{nt},\hat r_{nt})$ denotes the prediction of the specification under evaluation, the superscripts $\mathrm{lin}$ and $\mathrm{log}$ identify fitted values of the linear and log-linear equations, and all prediction errors are measured in original quantity units.

\paragraph{Linear specification.}
The linear specification consists of the net-demand and return equations~\eqref{eq:ridge_nd}--\eqref{eq:ridge_r}, with fitted values $\uline{\hat{d}}^{\,\mathrm{lin}}_{nt}$ and $\hat r^{\,\mathrm{lin}}_{nt}$. The return-rate specification below additionally requires a sales equation of the same form,
\begin{equation}
d_{nt} = a^{\mathrm d}_{n0} + (\bm a^{\mathrm d})^{\top}\bm v_{nt} + \sum_{m\in[N]} b^{\mathrm d}_{nm}\,\xi_{mt} + \epsilon^{\mathrm d}_{nt},
\label{eq:prediction_sales}
\end{equation}
where the superscript marks the sales equation, with fitted value $\hat d_{nt}$. In each equation the term $m=n$ is the direct effect and the terms $m\neq n$ are the spillover effects, one for every ordered pair of advertised and receiving products; the intercepts are product-specific, and the price, trend, and month coefficients are shared across products.

\paragraph{Ridge estimation.}
Each equation is estimated by ridge regression on a training set $\mathcal T$ of product-weeks, with the penalty applied to the advertising coefficients only. For the return equation~\eqref{eq:ridge_r},
\begin{equation*}
\big(\{\hat a_{n0}\}_{n\in[N]},\hat{\bm a},\mfh B\big)\in\argmin_{\{a_{n0}\}_{n\in[N]},\,\bm a,\,\mf B}\left\{\frac{1}{|\mathcal T|}\sum_{(n,t)\in\mathcal T}\big(r_{nt}-a_{n0}-\bm a^{\top}\bm v_{nt}-\bm b_n^{\top}\bm\xi_t\big)^2+\varrho\sum_{n\in[N]}\sum_{m\in[N]} s_{nm}^2\,b_{nm}^2\right\},
\end{equation*}
where $\bm b_n^{\top}$ is the $n$-th row of $\mf B$ and $s_{nm}$ is the root-mean-square value over $\mathcal T$ of the regressor multiplying $b_{nm}$, set to one when that regressor is identically zero, so that the penalty acts on standardized regressors while the coefficients are reported in original click units. Redundant intercept and control columns are removed by pivoted QR within each fit; because $\varrho>0$ and the retained controls have full column rank, the solution is unique even when the advertising regressors are rank deficient. The net-demand equation~\eqref{eq:ridge_nd} and the sales equation~\eqref{eq:prediction_sales} are estimated in the same way. The penalty $\varrho$ is selected separately for each equation from the 13-point grid $\{10^{-4},10^{-3.5},\ldots,10^{2}\}$ by five-fold cross-validation with folds grouped by week: each candidate is scored by its total validation squared error divided by the number of validation observations, exact ties favor the larger penalty, and standardization, control selection, and estimation are recomputed within each fold. The equation is then refitted at the selected penalty on all of $\mathcal T$.

\paragraph{Log-linear specification.}
The log-linear specification applies $\log(1+\cdot)$ to the outcome and to every advertising regressor:
\begin{align}
\log(1+\uline{d}_{nt}) &= \uline{a}^{\,\mathrm{log}}_{n0} + (\bmu a^{\mathrm{log}})^{\top}\bmu v_{nt} + \sum_{m\in[N]} \uline{b}^{\,\mathrm{log}}_{nm}\,\log(1+\xi_{mt}) + \uline{\epsilon}^{\,\mathrm{log}}_{nt},\label{eq:prediction_loglinear_nd}\\
\log(1+r_{nt}) &= a^{\mathrm{log}}_{n0} + (\bm a^{\mathrm{log}})^{\top}\bm v_{nt} + \sum_{m\in[N]} b^{\mathrm{log}}_{nm}\,\log(1+\xi_{mt}) + \epsilon^{\mathrm{log}}_{nt}.\label{eq:prediction_loglinear_r}
\end{align}
Both equations are estimated by the ridge procedure above, with penalties selected on the log scale. The transform retains zero outcomes and zero clicks. Predictions are returned to the original scale as $\uline{\hat{d}}^{\,\mathrm{log}}_{nt}=\exp\big(\widehat{\log(1+\uline{d}_{nt})}\big)-1$ and $\hat r^{\,\mathrm{log}}_{nt}=\exp\big(\widehat{\log(1+r_{nt})}\big)-1$, without smearing or other bias correction.

\paragraph{Monthly return rates.}
Let $M(t)$ denote the calendar year--month containing the starting date of week $t$, so that March 2022 and March 2023 are distinct periods. For a training set $\mathcal T$, the return rate of product $n$ in year--month $M$ is the sales-weighted ratio
\begin{equation}
\widehat\omega_{nM}(\mathcal T)
 =\frac{\sum_{(n,t)\in\mathcal T:\,M(t)=M}r_{nt}}
        {\sum_{(n,t)\in\mathcal T:\,M(t)=M}d_{nt}},
\label{eq:prediction_monthly_rate}
\end{equation}
when the denominator is positive. When the product-month is empty or has no training sales, the product's ratio of total returns to total sales over all its training observations is used, with the training-wide ratio as the final fallback. Rates are estimated from training observations only and held fixed when scoring test observations.

\paragraph{The five specifications.}
With $\hat d^{\,+}_{nt}=\max\{0,\hat d_{nt}\}$ the truncated sales prediction, the predictions $(\uline{\hat{d}}_{nt},\hat r_{nt})$ of the five specifications compared in \Cref{table:product_specific_prediction} are
\begin{align*}
\text{Linear:}\qquad &(\uline{\hat{d}}^{\,\mathrm{lin}}_{nt},\hat r^{\,\mathrm{lin}}_{nt}),\\
\text{Log-linear:}\qquad &(\uline{\hat{d}}^{\,\mathrm{log}}_{nt},\hat r^{\,\mathrm{log}}_{nt}),\\
\text{Return rate:}\qquad &\big((1-\widehat\omega_{n,M(t)})\,\hat d^{\,+}_{nt},\ \widehat\omega_{n,M(t)}\,\hat d^{\,+}_{nt}\big),\\
\text{Linear net/log return:}\qquad &(\uline{\hat{d}}^{\,\mathrm{lin}}_{nt},\hat r^{\,\mathrm{log}}_{nt}),\\
\text{Log net/linear return:}\qquad &(\uline{\hat{d}}^{\,\mathrm{log}}_{nt},\hat r^{\,\mathrm{lin}}_{nt}).
\end{align*}
Truncation at zero applies only to the sales input of the return-rate specification; predictions from the net-demand and return equations are not truncated. The two mixed specifications combine separately fitted equations and introduce no additional parameters. The return-rate specification is a prediction benchmark only; the numerical experiments use the linear net-demand and return equations.

\paragraph{Splits and reported errors.}
All specifications use the 695 cleaned product-week observations of Section~\ref{subsec::empirical evidences}. The in-sample MSE fits and scores the full panel. The out-of-sample columns average over 200 random 80/20 splits with 556 training and 139 test observations; coefficients, penalties, and return rates are re-estimated on the training portion of each split. For an evaluation set $\mathcal H$ of product-weeks, we define
\begin{align*}
\operatorname{MSE}_{\uline{d}}(\mathcal H)&=\frac{1}{|\mathcal H|}\sum_{(n,t)\in\mathcal H}\big(\uline{d}_{nt}-\uline{\hat{d}}_{nt}\big)^2,\\
\operatorname{MSE}_{r}(\mathcal H)&=\frac{1}{|\mathcal H|}\sum_{(n,t)\in\mathcal H}\big(r_{nt}-\hat r_{nt}\big)^2,\\
\operatorname{MSE}_{\mathrm{joint}}(\mathcal H)&=\tfrac12\big(\operatorname{MSE}_{\uline{d}}(\mathcal H)+\operatorname{MSE}_{r}(\mathcal H)\big).
\end{align*}
The in-sample and out-of-sample columns of \Cref{table:product_specific_prediction} report joint MSE, the two outcome-specific columns report the means of $\operatorname{MSE}_{\uline{d}}$ and $\operatorname{MSE}_{r}$ on the test portions over the 200 splits, and the 90\% range is the 5th--95th percentile range of split-level joint test MSE, not a confidence interval for its mean. Bold values mark column minima, including ties; because each mixed specification reuses the fitted net-demand or return equation of the linear or log-linear specification on the same splits, its outcome-specific error coincides with that of the parent specification, which produces the ties in \Cref{table:product_specific_prediction}. These splits are within the observed panel rather than chronological, and they are distinct from the week-based splits of the numerical experiments in Section~\ref{subsec:experiments}.

\paragraph{Full-panel coefficients.}
\Cref{table:product_specific_coefficients} reports the advertising-effect matrices $\mfu B$ and $\mf B$ estimated by ridge on the full panel, in units of outcome quantity per 100 clicks. Direct effects are positive for both outcomes for all ten products, while spillover effects vary in sign and magnitude. These full-panel estimates define the replay environment of Section~\ref{subsec:experiments}; the decision models are re-estimated on each training sample.

\begin{table}[htbp]
    \centering
    \caption{Full-panel ridge estimates of the advertising-effect matrices}
    \label{table:product_specific_coefficients}
    \begin{threeparttable}
        \renewcommand{\arraystretch}{1}
        \setlength{\tabcolsep}{4pt}
        \fontsize{9}{10}\selectfont
        \begin{tabular}{l*{10}{r}}
        \toprule
        \multicolumn{11}{c}{Net demand, $\mfu B$: receiving product $n$ (row), advertised product $m$ (column)} \\
        \midrule
         & P1 & P2 & P3 & P4 & P5 & P6 & P7 & P8 & P9 & P10 \\
        \midrule
        P1 & \textbf{1.14} & 2.21 & -1.56 & 0.74 & 0.83 & 0.63 & -0.54 & 0.40 & -1.87 & -0.34 \\
        P2 & -0.81 & \textbf{2.64} & 0.09 & 0.36 & 1.64 & -1.02 & -0.20 & -0.07 & 0.46 & 0.11 \\
        P3 & -0.10 & 0.16 & \textbf{1.63} & -0.51 & 1.77 & -1.01 & -0.04 & -0.15 & 0.25 & -0.18 \\
        P4 & 0.33 & 0.24 & -1.27 & \textbf{3.69} & 1.38 & 0.28 & 0.31 & 0.16 & -0.46 & 0.38 \\
        P5 & -0.69 & -0.26 & -0.51 & -0.46 & \textbf{8.12} & 0.86 & 0.35 & 0.36 & 0.22 & -0.72 \\
        P6 & 0.30 & -0.96 & -0.08 & 3.06 & 1.14 & \textbf{2.79} & 0.09 & 0.13 & -0.58 & -0.99 \\
        P7 & -0.25 & -1.63 & -0.61 & 0.12 & 1.19 & 0.25 & \textbf{0.73} & -0.03 & 0.68 & -0.50 \\
        P8 & -2.54 & 2.79 & 1.48 & 4.07 & 12.60 & 1.59 & -1.37 & \textbf{2.98} & -4.78 & 0.17 \\
        P9 & 0.08 & 0.61 & -0.57 & 0.04 & 1.65 & -0.11 & -0.02 & 0.05 & \textbf{1.94} & -0.60 \\
        P10 & 0.00 & -0.09 & -0.79 & 0.53 & 0.88 & 0.52 & 0.05 & 0.25 & -0.02 & \textbf{2.29} \\
        \midrule
        \multicolumn{11}{c}{Returns, $\mf B$: receiving product $n$ (row), advertised product $m$ (column)} \\
        \midrule
         & P1 & P2 & P3 & P4 & P5 & P6 & P7 & P8 & P9 & P10 \\
        \midrule
        P1 & \textbf{0.30} & 0.50 & -0.24 & 1.32 & 0.16 & -0.56 & 0.07 & 0.15 & -1.02 & -0.72 \\
        P2 & -0.30 & \textbf{0.81} & -0.04 & 0.26 & 0.86 & 0.18 & -0.11 & 0.03 & -0.18 & -0.21 \\
        P3 & -0.39 & 0.73 & \textbf{0.53} & -0.90 & 0.51 & -0.38 & -0.25 & 0.12 & 0.34 & -0.24 \\
        P4 & -0.03 & -0.08 & -0.15 & \textbf{0.24} & 0.24 & -0.16 & 0.10 & 0.05 & -0.18 & 0.35 \\
        P5 & -0.20 & -0.29 & -0.37 & 0.07 & \textbf{1.14} & 0.36 & 0.01 & -0.04 & -0.09 & 0.52 \\
        P6 & 0.59 & -1.32 & -0.73 & 0.47 & -0.37 & \textbf{0.54} & 0.19 & 0.07 & 0.22 & 0.77 \\
        P7 & -0.36 & -0.94 & -0.50 & 0.17 & 0.49 & 0.50 & \textbf{0.47} & -0.14 & 0.13 & 0.84 \\
        P8 & -0.52 & -0.04 & 0.23 & 0.28 & 1.84 & -0.60 & -0.13 & \textbf{0.17} & -0.37 & 1.45 \\
        P9 & -0.30 & 0.23 & -0.05 & 0.39 & 0.49 & 0.08 & -0.01 & 0.07 & \textbf{0.22} & 0.00 \\
        P10 & -0.02 & -0.18 & -0.40 & 0.08 & 0.48 & 0.01 & 0.10 & 0.08 & -0.01 & \textbf{0.36} \\
        \bottomrule
        \end{tabular}
        \begin{tablenotes}
            \item Notes. Entry $(n,m)$ is $\uline{b}_{nm}$ (upper panel) or $b_{nm}$ (lower panel) in~\eqref{eq:ridge_nd}--\eqref{eq:ridge_r}. 
            Bold diagonal entries are direct effects; off-diagonal entries are spillover effects. The selected ridge penalty is $0.1$ for both equations. 
        \end{tablenotes}
    \end{threeparttable}
\end{table}

\subsection{Design of the Numerical Experiments}\label{app:protocol}

We assess policy profits through a data-calibrated counterfactual simulation anchored in our partner's historical weekly observations. Each historical week supplies the covariates and fitted net-demand and return residuals, together with recorded click exposures; a product-week with no recorded clicks contributes zero exposure. Holding these scenario inputs fixed across policies, we use the calibrated response models to evaluate outcomes under alternative advertising decisions and compute profits under the associated inventory decisions. This design compares policies under common scenarios anchored in historical conditions. The following paragraphs specify the economic parameters, sampling and policy construction, and counterfactual evaluation procedure.

\paragraph{Parameters.}\label{app:key_para}
Products follow the fixed order P1--P10, and the instances with $N$ products use the first $N$. Retail prices are the full-panel mean weekly unit prices, in dollars, 
$
\bm p=(21.88,\,19.67,\,26.53,\,14.46,\,7.99,\,32.60,\,69.58,\,19.08,\,19.53,\,14.92).
$
The procurement cost is $\bm c=0.25\bm p$, an approximation of realized procurement cost; the backorder penalty is $\bm\alpha=0.30\bm p$, reflecting the international express-shipping cost of delayed fulfillment; the return-processing fee is $\bm\gamma=0.20\bm p$, calibrated from historical handling records; and the holding cost is $\bm h=0.06\bm p$, reflecting Amazon's monthly warehousing charges. The cost per click $\bm k$ is 1.5 times each product's average historical cost per click, recomputed on each training sample; the full-panel averages, in dollars, are $(0.31,\,0.27,\,0.26,\,0.31,\,0.30,\,0.30,\,0.17,\,0.39,\,0.22,\,0.46)$, and the multiplier represents a planning scenario with higher paid-search costs. All parameters are identical for~\ref{model::ETO} and~\ref{model::bar_ETR} within an instance.

\paragraph{Sampling.}
For each $(N,S)$ cell, we evaluate 20 instances, each formed by a random permutation of the 76 weeks into a 30-week training pool and a disjoint 46-week test set. The first $S$ weeks of the pool supply the training observations; training samples overlap across instances. The robustification and spillover comparisons use the same 400 matched instances. 
The computational comparison uses 20 matched instances per $(N,S)$ cell under the same sampling design.

\paragraph{Decision models.}
On each training sample, the net-demand and return equations~\eqref{eq:ridge_nd}--\eqref{eq:ridge_r} are estimated by the ridge procedure of Appendix~\ref{app:prediction_models}, including its penalty selection, and supply the estimates $(\mfhu A,\mfhu B,\mfh A,\mfh B)$ and the residuals $(\bmhu\epsilon_s,\bmh\epsilon_s)$ of Section~\ref{subsec:ETO}; the empirical distribution $\hat{\bbp}$ places mass $1/S$ on each training week's realized clicks and residuals $(\bmh\xi_s,\bmhu\epsilon_s,\bmh\epsilon_s)$. The support set $\mathcal Z$ of \Cref{thm::P_refo} uses $\bm\xi\ge\bm 0$ and sets the bounds on the net-demand and return errors to the componentwise minima and maxima of the training residuals. The Wasserstein distance uses the block-separable norm $\|\cdot\|_\oplus$ of Section~\ref{sec:ETR}, so the dual-norm constraints of~\ref{model::exact} and~\eqref{model::RS_L_CP} are second-order cone (SOC) constraints and both models are solved as mixed-integer SOC programs.

\paragraph{Target selection.}
Within each training sample, Algorithm~\ref{alg::tau} selects $\vartheta$ from the eleven candidates $\{0.70,\allowbreak 0.75,\allowbreak 0.80,\allowbreak 0.85,\allowbreak 0.90,\allowbreak 0.95,\allowbreak 0.98,\allowbreak 0.99,\allowbreak 0.995,\allowbreak 0.999,\allowbreak 1.00\}$ by five-fold cross-validation with folds grouped by week. Each fold re-estimates the net-demand and return models on its own training weeks, including the ridge penalty selection of Appendix~\ref{app:prediction_models}, and forms the empirical distribution from the resulting residuals; validation profits are computed by residual replay under the fold's own calibration, so no test-week information enters the selection. After selection, penalties, coefficients, and residuals are refitted on all $S$ training weeks. The implemented target is $\tau=\hat{\pi}^{\star}_{\code{ETO}}-(1-\vartheta)|\hat{\pi}^{\star}_{\code{ETO}}|$, which equals $\vartheta\hat{\pi}^{\star}_{\code{ETO}}$ when the empirical optimal profit is nonnegative and $(2-\vartheta)\hat{\pi}^{\star}_{\code{ETO}}$ otherwise, so that lowering $\vartheta$ always relaxes the target and $\tau\le\hat{\pi}^{\star}_{\code{ETO}}$ holds. All~\ref{model::bar_ETR} solves, in target selection and for the final policies, reach Gurobi's optimal status within a relative MIP-gap tolerance of $10^{-4}$.

\begin{algorithm}[ht]
\caption{Cross-validation selection of the target parameter $\vartheta$}
\label{alg::tau}
\fontsize{10}{10}\selectfont
\begin{algorithmic}[1]
\Require training weeks $\mathcal S$, candidate set $\Lambda$, number of folds $K$.
\State Partition $\mathcal S$ into $K$ disjoint folds $\mathcal S_1,\ldots,\mathcal S_K$ of as nearly equal size as possible.
\For{$j = 1,\ldots,K$}
    \State Estimate the net-demand and return models on $\mathcal S\setminus\mathcal S_j$ and form the empirical distribution $\hat{\bbp}_{-j}$ from their residuals.
    \For{each validation week $s\in\mathcal S_j$ with covariates $(\bmhu v_s,\bmh v_s)$}
        \State Solve~\ref{model::ETO} to obtain $\hat{\pi}^\star_{\code{ETO},-j}(s)$.
        \For{$\vartheta\in\Lambda$}
            \State Solve~\ref{model::bar_ETR} with target $\tau=\hat{\pi}^\star_{\code{ETO},-j}(s) - (1-\vartheta)\,|\hat{\pi}^\star_{\code{ETO},-j}(s)|$ to obtain $(\bmh x_{\vartheta}(s),\bmh y_{\vartheta}(s))$.
            \State Evaluate $(\bmh x_{\vartheta}(s),\bmh y_{\vartheta}(s))$ on week $s$ by residual replay and record its profit $\pi_{\vartheta}(s)$.
        \EndFor
    \EndFor
\EndFor
\Ensure $\vartheta^\star\in\displaystyle\argmax_{\vartheta\in\Lambda}\ \frac{1}{|\mathcal S|}\sum_{j\in[K]}\sum_{s\in\mathcal S_j}\pi_{\vartheta}(s)$, with ties resolved in favor of the larger $\vartheta$.
\end{algorithmic}
\end{algorithm}

\paragraph{Residual replay.}
Policies are evaluated on the 46 test weeks under a common full-panel ridge calibration of the net-demand and return models for each $N$. Replay retains the covariates and baseline variation of each test week, replaces its ad exposure by $\bm y\odot\bmh\xi_s$ for the policy under evaluation, and truncates return quantities at zero. Click exposures are taken from the original weekly records, so a product-week with no recorded clicks contributes zero exposure under either policy. Because the full-panel calibration uses test-period observations, the outcomes are calibrated simulations rather than prospective policy evaluations, and because training samples overlap across instances, the comparisons are descriptive summaries rather than independent field trials.

\paragraph{Computational comparison.}
\Cref{table:time_ave} solves~\ref{model::exact} and~\ref{model::bar_ETR} on 400 matched training samples, with 20 per $(N,S)$ cell, with the same estimates and residuals and the common fixed target $\vartheta=0.98$. Total computation time covers model construction and optimization and excludes the regression fit and the~\ref{model::ETO} solve, which the two models share. The empirical profit gap of a cell is the mean over its 20 instances of $100(\hat\pi_i^{\widebar{\code{ETR}}}-\hat\pi_i^{\code{EXACT}})/\hat\pi_i^{\code{EXACT}}$, where $\hat\pi_i^{\widebar{\code{ETR}}}$ and $\hat\pi_i^{\code{EXACT}}$ are the expected profits of the two solutions under the empirical distribution $\hat{\bbp}$ of instance $i$; all instances, including those with a negative gap, are retained. Times and gaps are rounded to three and six decimals, so a reported gap of $0.000000$ does not assert exact equality. 

\paragraph{Spillover comparison.}
The diagonal-only specification of \Cref{fig:robustification_spillover_density} estimates~\eqref{eq:ridge_nd}--\eqref{eq:ridge_r} with the spillover coefficients $\uline b_{nm}$ and $b_{nm}$, $m\ne n$, restricted to zero, so that each equation retains the product-specific intercept, the shared controls, and the direct effect. It is refitted by the same ridge procedure on the same training weeks, and its~\ref{model::ETO} and~\ref{model::bar_ETR} policies are obtained and evaluated as those of the full specification in the same 400 instances. The gain of each specification is computed by~\eqref{eq:aggregate_profit_gain} against the~\ref{model::ETO} policy of the same specification, so the full-specification gains in \Cref{fig:robustification_spillover_density} are those of \Cref{fig:per_comp}. The 95\% percentile bootstrap bands use 20,000 resamples of instances within each $(N,S)$ cell. Each resample recomputes the total profit difference and total ETO profit, using the same draws for both specifications; paired draws also yield the interval for the difference between the two gains.


\newpage		
\end{document}